\documentclass{amsart}
\usepackage{amsthm}
\usepackage{amsmath}
\usepackage{amssymb}
\usepackage{graphics}
\usepackage[dvipdfm]{graphicx}
\usepackage{bmpsize}

\theoremstyle{plain}
\newtheorem{lemma}{Lemma}[section]
\newtheorem{theorem}[lemma]{Theorem}
\newtheorem{propn}[lemma]{Proposition}
\newtheorem{cor}[lemma]{Corollary}

\theoremstyle{definition}
\newtheorem{defn}[lemma]{Definition}

\theoremstyle{remark}
\newtheorem{remark}[lemma]{Remark}

\newcommand{\C}{\mathbb C}
\newcommand{\R}{\mathbb R}
\newcommand{\Z}{\mathbb Z}

\newcommand{\D}{\mathbb D}
\newcommand{\est}{e^{2\pi(s+it)}}
\newcommand{\util}{\tilde u}
\newcommand{\wtil}{\tilde w}
\newcommand{\vtil}{\tilde v}
\newcommand{\jtil}{\tilde J}
\newcommand{\J}{\mathcal J}
\newcommand{\V}{\mathcal V}
\newcommand{\M}{\mathcal M}

\newcommand{\E}{\mathcal E}
\newcommand{\DD}{\mathcal D}

\renewcommand{\O}{\mathcal O}
\renewcommand{\P}{\mathcal P}
\newcommand{\wind}{\text{wind}}
\newcommand{\link}{\text{link}}

\newcommand{\maslov}{\text{Maslov}}
\renewcommand{\sp}{\text{Sp}}

\newcommand{\norma}[1]{\left\|{#1}\right\|}

\newcommand{\cl}[1]{\overline{{#1}}}

\begin{document}

\title[Systems of global surfaces of section]{Systems of global surfaces of section for dynamically convex Reeb flows on the 3-sphere}

\author[U.L. Hryniewicz]{Umberto L. Hryniewicz}

\address{Departamento de Matem\'atica Aplicada -- Instituto de Matem\'atica\\
Universidade Federal do Rio de Janeiro\\ Rio de Janeiro, Brazil}
\email{umberto@labma.ufrj.br}


\begin{abstract}
We characterize which closed Reeb orbits of a dynamically convex contact form on the 3-sphere bound disk-like global surfaces of section for the Reeb flow, without any genericity assumptions. We show that these global surfaces of section come in families, organized as open book decompositions. As an application we obtain new global surfaces of section for the Hamiltonian dynamics on strictly convex 3-dimensional energy levels.
\end{abstract}

\maketitle

\section{Introduction}\label{intro}

Our goal is to study the existence question for global surfaces of section for Reeb dynamics associated to a dynamically convex contact form on $S^3$, initiated by Hofer, Wysocki and Zehnder~\cite{char1,char2,convex,fols}.

Recall that a 1-form $\lambda$ on a $3$-manifold $M$ is a contact form if $\lambda\wedge d\lambda$ never vanishes. Associated to $\lambda$ is a (co-oriented) contact structure given by the 2-plane distribution
\begin{equation}\label{contact_str}
\xi = \ker\lambda,
\end{equation}
and a vector field $R$ which is uniquely determined by the equations
\begin{equation}\label{reebvector}
\begin{array}{ccc}
  i_R d\lambda = 0 &  & i_R \lambda = 1.
\end{array}
\end{equation}
It is called the Reeb vector field and its flow $\phi_t$ is referred to as the Reeb flow. By a periodic Reeb orbit $P$ we mean a pair $(x,T)$ where $T>0$ and $x$ is a $T$-periodic trajectory of $\phi_t$, which is called prime or simply covered if $T$ is its minimal positive period. Pairs with the same geometric image and period are identified, and the set of equivalence classes will be denoted by $\P$. When $c_1(\xi)$ vanishes on $\pi_2(M)$ one can associate to any contractible $P\in \P$ its Conley-Zehnder index $\mu_{CZ}(P) \in \Z$ even when $P$ is degenerate, see~\ref{cz_geom_def} and~\ref{cz_operators_description} in section~\ref{cz_index_descriptions} for two equivalent definitions. This is an invariant of the linearized dynamics at $P$ which is lower semi-continuous for $C^1$ perturbations of the linearized Reeb flow, in particular lower bounds for $\mu_{CZ}$ are preserved under such perturbations. In~\cite{convex} we find the following important definition.

\begin{defn}[Hofer, Wysocki and Zehnder]
The contact form $\lambda$ is dynamically convex if $c_1(\xi)$ vanishes on $\pi_2(M)$ and $\mu_{CZ}(P) \geq 3$ for every contractible periodic Reeb orbit $P$.
\end{defn}

A global surface of section for a flow without rest points on a closed connected $3$-manifold is a compact connected embedded surface $\Sigma$ such that $\partial \Sigma$ consists of periodic trajectories, $\Sigma\setminus \partial\Sigma$ is transverse to the flow, and every trajectory not in $\partial \Sigma$ hits $\Sigma$ infinitely many times in the future and in the past. We shall also consider families of global surfaces of section organized as an open book decomposition. Recall that an open book decomposition of a closed 3-manifold $M$ is a pair $(K,\Pi)$, where $K\subset M$ is an oriented link and $\Pi: M\setminus K \to \R/\Z$ is a (smooth) fibration such that each $\Pi^{-1}(\vartheta)$ is the interior of a compact embedded oriented surface $S_\vartheta$ satisfying $\partial S_\vartheta = K$ (orientations included). $K$ is called the binding and $\Pi^{-1}(\vartheta)$ is called a page.

\begin{defn}\label{open_book_adapted}
An open book decomposition $(K,\Pi)$ of $M$ is adapted to the contact form $\lambda$ if the Reeb vector field is positively tangent to $K$ and the pages are global surfaces of section for the Reeb flow.
\end{defn}

The following remarkable statement is the main result of~\cite{convex}.

\begin{theorem}[Hofer, Wysocki and Zehnder]\label{theo_convex}
Let $\lambda$ be a dynamically convex contact form on $S^3$. Then there exists an embedded disk $D_0 \hookrightarrow S^3$ that is a global surface of section for the Reeb flow. Its boundary $\partial D_0 = P_0$ is a closed Reeb orbit satisfying $\mu_{CZ}(P_0) = 3$, and $D_0 \setminus \partial D_0$ is a page of an open book decomposition of $S^3$ adapted to $\lambda$.
\end{theorem}

Theorem~\ref{theo_convex} has very strong consequences for the dynamics. Since the Poincar\'e first return map to $D_0\setminus \partial D_0$ is well-defined and preserves the finite area form $\omega_0|_{D_0\setminus \partial D_0}$, results of J. Franks~\cite{franks} on area-preserving disk maps can be applied.

\begin{cor}[Hofer, Wysocki and Zehnder]
The Reeb dynamics of a dynamically convex contact form on $S^3$ admit either two or infinitely many geometrically distinct periodic orbits.
\end{cor}

It should be noted that dynamical convexity of a contact form also imposes restrictions on the contact manifold. Recall that an embedded disk $F$ on a contact 3-manifold $(M,\xi)$ is overtwisted if $T\partial F \subset \xi$ and $T_pF \not= \xi_p$, $\forall p\in\partial F$. The contact structure $\xi$ is called tight if there are no overtwisted disks and, when $M$ is closed Hofer, Wysocki and Zehnder show in~\cite{char1} that if $\xi = \ker\lambda$ for some dynamically convex contact form $\lambda$ then $\xi$ is tight and $\pi_2(M)$ vanishes.

Let us describe a well-known family of examples. Equipping $\R^4$ with coordinates $(q_1,p_1,q_2,p_2)$, the Liouville form $$ \lambda_0 = \frac{1}{2} (q_1dp_1 - p_1dq_1 + q_2dp_2 - p_2dq_2) $$ restricts to a contact form on the boundary $S$ of a smooth bounded domain of $\R^4$ which is star-shaped with respect to the origin. Writing $S = \{ \sqrt{f(x)}x : x\in S^3 \}$ for some smooth $f:S^3\to \R^+$ then $\Psi^*(\lambda_0|_S) = f\lambda_0|_{S^3}$ where $\Psi:S^3\to S$ is the diffeomorphism given by $\Psi(x) = \sqrt{f(x)}x$. In this case we say $f\lambda_0|_{S^3}$ comes from $S$. The integral leaves of the associated characteristic line bundle $$ \mathcal L = \bigcup_{p\in S} (T_pS)^{\omega_0}, \ \ \text{ where } (T_pS)^{\omega_0} = \{ v \in \R^4 \mid \omega_0(v,w)=0 \ \forall w\in T_pS \}, $$ coincide with the integral curves of the Reeb flow associated to $\lambda_0|_S$. In~\cite{convex} Hofer, Wysocki and Zehnder show that $\lambda_0|_S$ is dynamically convex when $S$ is strictly convex. By a theorem of Eliashberg every tight contact form on $S^3$ is diffeomorphic to a contact form coming from some star-shaped domain.

Motivated by Theorem~\ref{theo_convex}, one may ask which closed orbits of the Reeb dynamics associated to a tight contact form on $S^3$ bound global disk-like global surfaces of section. This was answered in~\cite{pedro} when the contact form is non-degenerate, that is, when the infinitesimal Poincar\'e return map associated to any orbit in $\mathcal P$ does not have $1$ as an eigenvalue. The conditions also depend on contact-topological properties of the closed Reeb orbit seen as a transverse knot, encoded in its self-linking number.

\begin{defn}
Let $L \hookrightarrow (M,\xi)$ be a transverse knot, that is $TL \pitchfork \xi$, and $\Sigma \hookrightarrow M$ be a Seifert surface\footnote{A compact connected orientable embedded surface satisfying $\partial\Sigma = L$.} for $L$. The vector bundle $\xi|_\Sigma$ admits a non-vanishing section $Z$. If $\exp$ is any exponential map on $M$ and $\epsilon>0$ is small then $p \in L \mapsto \exp (\epsilon Z_p)$ is a diffeomorphism between $L$ and an embedded loop $L^\prime \hookrightarrow M$ satisfying $L \cap L^\prime = \emptyset$. An orientation for $\Sigma$ also orients $L$ which, in turn, orients $L^\prime$. The self-linking number $sl(L,\Sigma) \in \Z$ is defined as the oriented intersection number of $L^\prime$ and $\Sigma$; it is independent of all choices.
\end{defn}

If $c_1(\xi)$ vanishes on $\pi_2(M)$, $P = (x,T)$ is a periodic Reeb orbit and $D$ is an embedded disk spanning $x(\R)$, then $sl(x(\R),D)$ is independent of the choice of the disk with these properties  and will be denoted by $sl(P)$.

\begin{theorem}[\cite{pedro}]\label{tight_reeb_flows}
Let $\lambda$ be a non-degenerate tight contact form on $S^3$. Then a prime closed Reeb orbit $P$ bounds a disk-like global section for the Reeb flow if, and only if, $P$ is unknotted, $\mu_{CZ}(P) \geq 3$, $sl(P)=-1$ and $P$ is linked to every orbit $P'$ satisfying $\mu_{CZ}(P')=2$. Moreover, one finds an open book decomposition with disk-like pages of $S^3$ adapted to $\lambda$ with binding~$P$.
\end{theorem}

In particular, when $\lambda$ is a non-degenerate dynamically convex contact form on $S^3$ then a prime closed Reeb orbit bounds a disk-like global surface of section if, and only if, it is unknotted and has self-linking number $-1$. This fact was first proved in~\cite{hry}. Our first result reads as follows.

\begin{theorem}\label{main}
Let $\lambda$ be any dynamically convex contact form on $S^3$. Then a prime periodic Reeb orbit $\bar P$ bounds a disk-like global surface of section if, and only if, it is unknotted and $sl(\bar P)=-1$. Moreover, one can find an open book decomposition with disk-like pages of $S^3$ adapted to $\lambda$ with binding $\bar P$.
\end{theorem}

In the case of a contact form on $S^3$ that comes from an ellipsoid, one easily sees that the axes are bindings of adapted open book decompositions as in Definition~\ref{open_book_adapted}. However, when this ellipsoid is not the round 3-sphere, only one open book is detected by Theorem~\ref{theo_convex}, namely, the one with the axis of smaller action as binding. The other prime orbit (with larger action) has Conley-Zehnder index $>3$. As an application of Theorem~\ref{main} we show that the same global picture holds in general.

\begin{theorem}\label{2nd_global}
Let $\lambda$ be any dynamically convex contact form on $S^3$ and let $D_0$ be any disk-like global surface of section for the Reeb flow of $\lambda$. Consider any Reeb orbit $P_1$ associated to a fixed point of the first return map to $D_0\setminus \partial D_0$. Then $P_1$ is unknotted, $sl(P_1)=-1$ and, consequently, bounds a disk-like global surface of section $D_1$ which is a page of an open book decomposition of $S^3$ adapted to $\lambda$.
\end{theorem}

Brouwer's translation theorem implies the existence of a fixed point of the first return map to the disk $D_0$ from Theorem~\ref{theo_convex}, so that Theorem~\ref{main} gives new disk-like global sections, geometrically distinct from $D_0$. Thus, the Hamiltonian flow is globally twisting in two ``independent'' directions: around $\partial D_0$ and around $\partial D_1$, see Figure 1.

\begin{figure}
\includegraphics[width=320\unitlength]{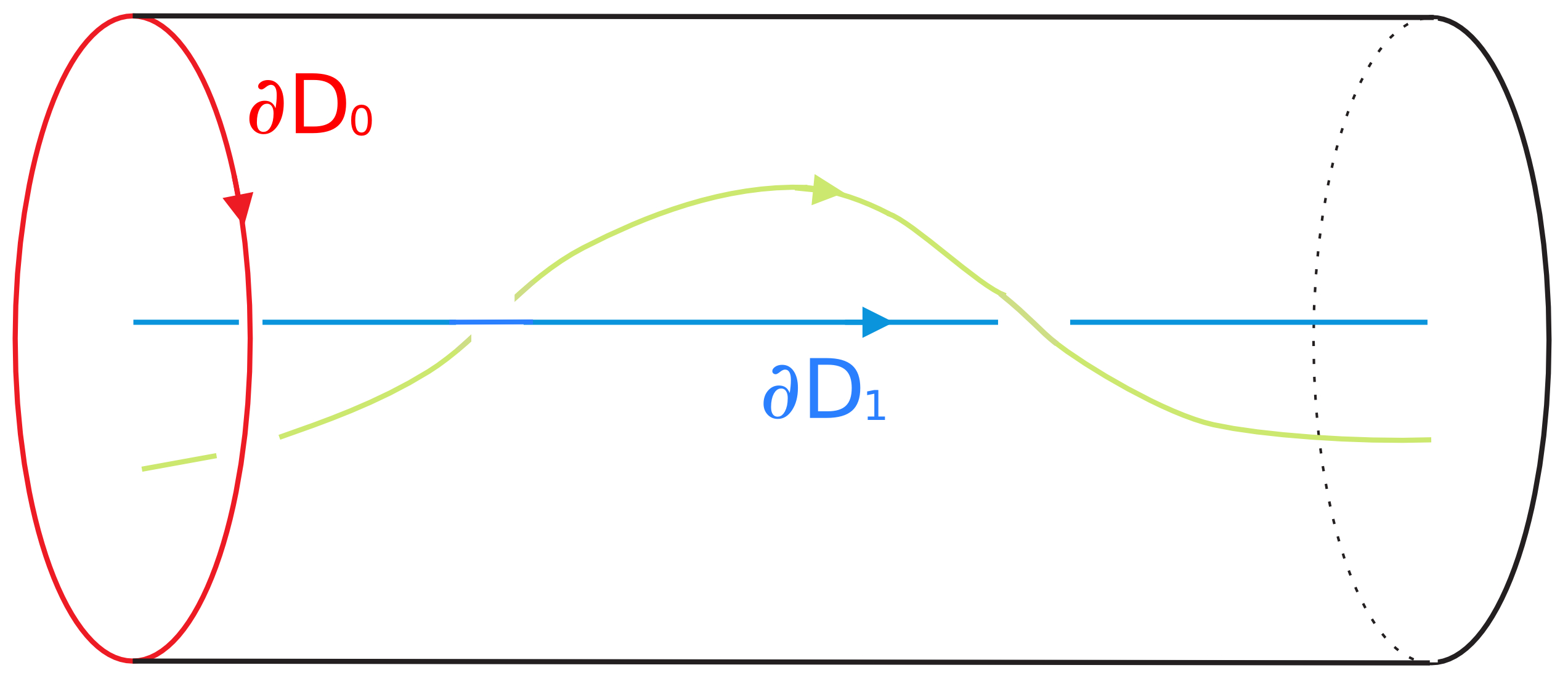}
\caption{$S \setminus \partial D_0$ gets identified with the interior of a solid torus, the flow is globally twisting with respect to $\partial D_0$ and $\partial D_1$ in a monotone fashion.}
\end{figure}

It should be noted that the Reeb dynamics of tight contact forms on $S^3$ are conjugate or semi-conjugate to the Hamiltonian dynamics on certain energy levels of many important classical Hamiltonian systems with two degrees of freedom. For example, in a recent paper~\cite{3body} Albers, Fish, Frauenfelder, Hofer and van Koert obtain disk-like global sections in the planar restricted 3-body problem. After a special transformation, a component of the levels slightly below the first Lagrange value is covered by a compact component of an energy level of the transformed Hamiltonian to which Theorem~\ref{theo_convex} can be applied. Theorem~\ref{main} can be applied as well to study the existence of new disk-like global sections in this case.

These results are examples of results in Symplectic Dynamics, as explained by Bramham and Hofer~\cite{BH}. Symplectic Dynamics has roots in the work of Poincar\'e, but recent techniques have proved to be extremely successful to uncover new global phenomena in Hamiltonian dynamics. Among these we would like to emphasize methods from holomorphic curve theory introduced in symplectic geometry by Gromov~\cite{gromov}. The proofs of the above mentioned results are based on a holomorphic curve theory in symplectizations introduced by Hofer in~\cite{93}, further developed by Hofer, Wysocki and Zehnder~\cite{props1,props2,props3,props4} and later by many other authors. \\

\noindent {\bf Outline of main arguments.} If the periodic orbit $\bar P$ bounds a disk-like global surface of section then, clearly, it must be unknotted. In Proposition 2.1 of~\cite{hry} we show that $sl(\bar P) = -1$; this follows since $\bar P$ bounds a disk transverse to the Reeb flow. Necessity in Theorem~\ref{main} follows, and sufficiency remains to be addressed. In view of Theorem~\ref{tight_reeb_flows} we only need to pass to the degenerate case, which will be done by following~\cite{convex,hry,pedro} closely. It is always possible to find non-degenerate contact forms $\lambda_k$ admitting $\bar P$ as a closed Reeb orbit and satisfying $\lambda_k \to \lambda$ in $C^\infty$. We will obtain finite-energy pseudo-holomorphic planes $\util_k$ (with respect to a suitable almost complex structure induced by $\lambda_k$) which project onto pages of an open book decomposition adapted to $\lambda_k$, and then pass to the limit as $k\to\infty$ to show that these planes converge to pages of the desired open book adapted to $\lambda$. However, Theorem~\ref{tight_reeb_flows} does not apply directly to $\lambda_k$ and $\bar P$ since there could be closed $\lambda_k$-Reeb orbits with very high action which have Conley-Zehnder index $= 2$ and are not linked to $\bar P$. So we need to carefully revisit the arguments from~\cite{convex,hry,pedro} to circumvent this difficulty. For details we refer to Section~\ref{deg_case}. Theorem~\ref{2nd_global} follows from the fact that any orbit $P_1$ given by a fixed point of the return map has self-linking number $-1$. This is the content of Proposition~\ref{propn_sl} which is proved in Section~\ref{section_sl}. \\

\noindent {\bf Acknowledgements.} This work has its origin in the author's Ph.D. thesis, and we thank Professor Helmut Hofer for early influential discussions on this topic. Special thanks to Pedro A. S. Salom\~ao for numerous conversations about our results and their future applications, for suggestions on how to generalize previous versions of this paper, and for his mathematical help. We would also like thank the referee for suggesting many improvements of the paper. This material is based upon work supported by the National Science Foundation under agreement No. DMS-0635607. Any opinions, findings and conclusions or recommendations expressed in this material are those of the authors and do not necessarily reflect the views of the National Science Foundation. We thank the IAS for its hospitality during part of the preparation of this paper.

\section{Preliminaries}\label{preliminaries}

\subsection{Descriptions of the Conley-Zehnder index in $3$-dimensions}\label{cz_index_descriptions}

Here we recall the basic facts about the Conley-Zehnder index, and from the theory of pseudo-holomorphic curves in symplectizations. The reader may skip this section on a first read, referring back only for the notation.

Consider the set $\Sigma^*$ of smooth paths $\varphi: [0,1] \to \sp(1)$ of symplectic $2\times 2$ matrices satisfying $\varphi(0) = I$ and $\det (\varphi(1)-I) \not= 0$. As explained in~\cite{fols}, the Conley-Zehnder index $\mu:\Sigma^* \to \Z$ is uniquely determined by the following axioms:
\begin{itemize}
  \item \textbf{Homotopy:} If $\{\varphi^s\}$ is a homotopy of paths in $\Sigma^*$ then $\mu(\varphi^s)\equiv \mu(\varphi^0)$.
  \item \textbf{Maslov Index:} If $\psi:([0,1],\{0,1\}) \to (\sp(1),I)$ is a smooth closed loop and $\varphi \in \Sigma^*$ then $\mu(\psi\varphi) = 2\maslov(\psi) + \mu(\varphi)$.
  \item \textbf{Inversion:} If $\varphi \in \Sigma^*$ then $\mu(\varphi^{-1}) = -\mu(\varphi)$.
  \item \textbf{Normalization:} $\mu(t\in[0,1]\mapsto e^{i\pi t}) = 1$.
\end{itemize}
Assuming $c_1(\xi)$ vanishes on $\pi_2(M)$, the Conley-Zehnder index of a contractible periodic Reeb orbit $P = (x,T)$ is defined as follows. Choose a disk map $F : \D \to M$ satisfying $F(e^{i2\pi t}) = x(Tt)$, fix a smooth $d\lambda$-symplectic trivialization $\Psi$ of $F^*\xi$ and consider $\varphi = \Psi_{e^{i2\pi t}} \cdot d\phi_{Tt}|_{x(0)} \cdot \Psi_{1}^{-1}$. Then $\varphi \in \Sigma^*$ if, and only if $P$ is non-degenerate. The integer $\mu(\varphi)$ is independent of $F$ and $\Psi$, allowing one to define
\begin{equation}\label{cz_index_defn}
 \mu_{CZ}(P) = \mu(\varphi).
\end{equation}
Below we discuss how it is possible to extend the Conley-Zehnder index to periodic orbits of arbitrary contact forms, following~\cite{convex}.

In our applications we need two concrete descriptions of the integer $\mu$ satisfying the axioms described in the introduction.

\subsubsection{The geometrical description}\label{cz_geom_def}

For any smooth path $\varphi:[0,1] \to \sp(1)$ satisfying $\varphi(0)=I$ there exist unique smooth functions $r,\theta : [0,1]\times [0,2\pi] \to \R$ satisfying $\varphi(t)e^{is} = r(t,s)e^{i\theta(t,s)}$, $r>0$, $\theta(0,s)=s$. The image of the map $$ \Delta: [0,2\pi] \to \R, \ \ s \mapsto \frac{\theta(1,s)-s}{2\pi} $$ is a closed interval $I(\varphi)$ with length\footnote{In fact, suppose $|\Delta(s_2)-\Delta(s_1)|=1/2$. We may assume $0\leq s_1<s_2<\pi$ without loss of generality. Defining $g(t) = \theta(t,s_2)-\theta(t,s_1)$ we have $g(1)=g(0) \pm \pi$. By continuity there must be a value $t^* \in (0,1)$ satisfying $g(t^*) \in \{0,\pi\}$, which implies $\varphi(t^*)e^{is_2} \in \R\varphi(t^*)e^{is_1}$. We conclude $\Delta(s_1)=\Delta(s_2)$ by linearity of the equation and uniqueness of solutions, a contradiction.} $<1/2$.

\begin{lemma}\label{deg_path}
$\partial I(\varphi) \cap \Z \neq \emptyset \Leftrightarrow \varphi \not\in \Sigma^*(1)$.
\end{lemma}

\begin{proof}
Assume $\Delta(\bar s) = \inf \Delta \in \Z$. Then $\theta(1,\bar s) \in \bar s +2\pi\Z$ and $0=2\pi\Delta'(\bar s) = \theta_s(1,\bar s)-1$. We must have $\varphi(1)e^{i\bar s} = r(1,\bar s)e^{i\bar s}$ and $\varphi(1)ie^{i\bar s} = \frac{d}{ds}|_{s=\bar s} \varphi(1)e^{i s} = r_s(1,\bar s)e^{i\bar s} + r(1,\bar s)\theta_s(1,\bar s)ie^{i\bar s}$. Thus $1=\det \varphi(1) = r(1,\bar s)^2\theta_s(1,\bar s) = r(1,\bar s)^2$ shows that $r(1,\bar s)=1$ is an eigenvalue of $\varphi(1)$. If $\sup \Delta \in \Z$ the argument is analogous. This proves $\partial I(\varphi) \cap \Z \neq \emptyset \Rightarrow \varphi \not\in \Sigma^*(1)$.

Now assume $\varphi \not\in \Sigma^*(1)$ and choose $s_0 \in [0,\pi)$ such that $\varphi(1)e^{is_0} = e^{is_0}$, or equivalently, $\theta(1,s_0) = s_0 + 2k\pi$ for some $k\in\Z$. If $\varphi(1)=I$ then $I(\varphi) \subset \Z$ and we are done, so we also assume $\varphi(1)\neq I$. The final matrix $\varphi(1)$ must be a symplectic shear satisfying $\varphi(1)  e^{i(s_0+\pi/2)} = e^{i(s_0+\pi/2)} + c e^{is_0}$ with $c\neq 0$. By linearity we have $\theta(t,s_0+\pi) = \theta(t,s_0)+\pi \ \forall t$ and $s\in(s_0,s_0+\pi) \Rightarrow \theta(t,s_0)<\theta(t,s)<\theta(t,s_0)+\pi \ \forall t$. In particular $s_0+2k\pi<\theta(1,s)<s_0+(2k+1)\pi \ \forall s\in(s_0,s_0+\pi)$. Moreover, $\theta(1,s) \neq s+2k\pi \ \forall s\in(s_0,s_0+\pi)$ since, otherwise, $c= 0$ and $\varphi(1)=I$. It is clear that $c>0 \Leftrightarrow \theta(1,s) < s+2k\pi \ \forall s\in(s_0,s_0+\pi)$ and $c<0 \Leftrightarrow \theta(1,s) > s+2k\pi \ \forall s\in(s_0,s_0+\pi)$. The first case implies $\sup I(\varphi)=k$, the second case implies $\inf I(\varphi)=k$.
\end{proof}

On the set of closed intervals $J$ of length $<1/2$ satisfying $\partial J \cap \Z = \emptyset$ we may consider the function
\begin{equation}\label{}
  \hat \mu(J) = \left\{ \begin{aligned} &2k \text{ if } k\in J \\ & 2k+1 \text{ if } J \subset (k,k+1). \end{aligned} \right.
\end{equation}
One checks that $\mu(\varphi) = \hat \mu(I(\varphi))$ satisfies the required axioms for the Conley-Zehnder index of paths in $\Sigma^*(1)$. The function $\hat \mu$ may be extended to the set of all closed intervals of length $<1/2$ by
\begin{equation}\label{}
  \hat \mu(J) = \lim_{\epsilon\to0^+} \hat\mu(J-\epsilon).
\end{equation}
This induces an extension $\mu(\varphi) = \hat \mu(I(\varphi))$ of the Conley-Zehnder index to the set of all smooth paths $\varphi:[0,1]\to \sp(1)$ satisfying $\varphi(0)=I$. Clearly $\mu$ is lower semi-continuous if the space of such paths is endowed with the $C^0$-topology, since small variations of $\varphi$ induce small variations of the end-points of $I(\varphi)$.

\subsubsection{A description via self-adjoint operators}\label{cz_operators_description}

A smooth path $\varphi:[0,1] \to \sp(1)$ satisfies a differential equation $-i\varphi'-S\varphi=0$, for some (unique) smooth path of symmetric matrices $S(t), \ t\in[0,1]$. We may consider the unbounded self-adjoint operator $L$ on $L^2(\R/\Z,\R^2)$ given by
\begin{equation}\label{}
  \begin{array}{ccc}
    Lv = -iv'-Sv & \text{ where } & i \simeq \begin{bmatrix} 0 & -1 \\ 1 & 0 \end{bmatrix}.
  \end{array}
\end{equation}
It has compact resolvent, its spectrum $\sigma(L)$ is a discrete set of real eigenvalues accumulating at $\pm\infty$, geometric and algebraic multiplicities coincide, and all eigenvalues have multiplicity at most $2$. A non-trivial eigenvector $v(t) = r(t)e^{i\theta(t)}$ associated to some $\lambda\in\sigma(L)$ never vanishes, so the total winding $\theta(1)-\theta(0)$ is well-defined. This number does not depend on the eigenvector for $\lambda$, and we have an integer $\wind(\lambda) := (\theta(1)-\theta(0))/2\pi$. It turns out that for each $k\in\Z$ there are precisely two (multiplicities counted) eigenvalues with $\wind=k$, and that $\lambda \leq \mu \Rightarrow \wind(\lambda)\leq \wind(\mu)$. For more details see section 3 from~\cite{props2}.

Let us denote $\lambda^- = \max \{\lambda \in \sigma(L) \mid \lambda < 0\}$, $\lambda^+ = \min \{\lambda \in \sigma(L) \mid \lambda \geq 0\}$. Set $p=0$ if $\wind(\lambda^-) = \wind(\lambda^+)$, or $p=1$ if $\wind(\lambda^-) = \wind(\lambda^+)-1$. In case $\varphi \in \Sigma^*(1)$ we have $0 \not\in \sigma(L)$, so $\lambda^+ > 0$ and one checks that
\begin{equation}\label{CZ_operators}
  \mu(\varphi) = 2\wind(\lambda^-)+p
\end{equation}
also satisfies the axioms for the Conley-Zehnder index, see~\cite{props2} for the proof. Even in case $\varphi \not\in \Sigma^*(1)$, Hofer, Wysocki and Zehnder define $\mu(\varphi)$ by~\eqref{CZ_operators}, which provides an extension to all smooth paths in $\sp(1)$ starting at $I$. When the space of such paths is endowed with the $C^1$-topology, small variations of $\varphi$ induce $C^0$-small variations of $S$, which in turn induce small variations of the spectrum, as explained in~\cite{kato}. It becomes clear that the extension of $\mu$ described above is lower-semicontinuous since $\mu(\varphi)$ either stays constant or jumps up (by $1$ or $2$) under small variations of the spectrum.

\subsubsection{Comparing both extensions}\label{both_extensions}

\begin{lemma}\label{lemma_both_extensions}
Both extensions of $\mu$ described in~\ref{cz_geom_def} and~\ref{cz_operators_description} coincide.
\end{lemma}

\begin{proof}
We need only to consider a path $\varphi\not\in \Sigma^*(1)$. We denote by $\mu_1(\varphi)$ the extension described in~\ref{cz_geom_def}, and by $\mu_2(\varphi)$ the one from~\ref{cz_operators_description}.

Let $I(\varphi)= [a,b]$ be the ``winding interval'' described in~\ref{cz_geom_def}. Consider also some pair of eigenvalues $\nu^+>0$ and $\nu^-<0$ of the operator $L = -i\partial_t - S$ ($S=-i\varphi'\varphi^{-1}$). We claim that $\wind(\nu^+) > a$ and $\wind(\nu^-) < b$. To prove the first inequality choose an eigenvector $v(t)$ satisfying $Lv = \nu^+v$ and consider $u(t) = \varphi(t)\cdot v(0)$. Then the vector $z(t) = v(t) \overline{u(t)}$ satisfies
\begin{equation}\label{}
  -i\dot z = (Sv)\bar u - v (\overline{Su}) + \nu^+ z
\end{equation}
Whenever $v \in \R^+ u$ we have $z\in \R$ and $(Sv)\bar u - v (\overline{Su}) \in i\R$. Then, writing $z(t) = \rho(t)e^{i\vartheta(t)}$, we must have $\vartheta \in 2\pi\Z \Rightarrow \dot\vartheta = \nu^+>0$. Thus the total angular variation $\vartheta(1)-\vartheta(0)$ of $z$ is strictly positive, in other words, the total angular variation of $v$ is strictly larger than that of $u$, as we wanted to show. The other inequality is proved analogously.

By Lemma~\ref{deg_path} one finds $k\in \Z$ such that $\{k\} = I(\varphi) \cap \Z = \partial I(\varphi) \cap \Z$. The eigenvalue $0$ of the operator $L$ must have winding precisely $k$. Note that $$ \varphi(1)=I \Leftrightarrow I(\varphi)=\{k\} \Leftrightarrow 0 \text{ has multiplicity two as an eigenvalue of }L. $$ The first equivalence is obvious, and the second follows from the fact that two non-colinear eigenvectors of $L$ for the eigenvalue $0$ are pointwise linearly independent. In this case, the windings of the largest negative eigenvalue and that of the smallest positive eigenvalue are $k-1$ and $k+1$, respectively. According to the definitions explained in~\ref{cz_geom_def} and~\ref{cz_operators_description} and the spectral properties of $L$ we have $\mu_1(\varphi) = \mu_2(\varphi) = 2k-1$.

It remains to handle the cases where $\varphi(1)\neq I$. If $a=k$ and $b>a$ then $\mu_1(\varphi) = 2k$. As proved above, the winding of the smallest positive eigenvalue of $L$ is $>k$. The spectral properties of $L$ explained in~\ref{cz_operators_description} imply that this winding is precisely $k+1$. Since $0$ is a simple eigenvalue, the winding of the largest negative eigenvalue must also be $k$, proving $\mu_2(\varphi)=2k$.

If $b=k$ and $a<b$ then $\mu_1(\varphi) = 2k-1$. As before, the winding of the largest negative eigenvalue of $L$ is $<k$ and, consequently, it must be precisely $k-1$. Thus $\mu_2(\varphi) = 2k-1$.
\end{proof}

\subsection{Pseudo-holomorphic curves in symplectizations}\label{hol_curve_theory}

Let us recall the basic definitions and facts of the theory as introduced by Hofer in~\cite{93}. Throughout we fix a contact form $\lambda$ on the closed $3$-manifold $M$, with Reeb vector $R$ and induced contact structure $\xi = \ker\lambda$. The projection $\R\times M \to M$ onto the second factor is denoted by $\pi_M$. Consider also
\begin{equation}\label{projection}
 \pi : TM \to \xi
\end{equation}
the projection along the Reeb direction.

\subsubsection{Almost complex structures}

A complex structure $J$ on $\xi$ is said to be $d\lambda$-compatible if $d\lambda(\cdot,J\cdot)$ is a metric. The space of these complex structures, which is well-known to be contractible, is denoted by $\J(\xi,d\lambda)$. Any $J \in \J(\xi,d\lambda)$ induces an almost complex-structure $\jtil$ on $\R\times M$ by
\begin{equation}\label{alm_cpx_str_def}
  \begin{array}{cc}
    \jtil \cdot \partial_a = R, & \jtil|_\xi = J.
  \end{array}
\end{equation}
Above we see $TM$ and $\xi$ as ($\R$-invariant) subbundles of $T(\R\times M)$, and denote by $a$ the $\R$-coordinate.

\subsubsection{Finite-energy curves}\label{fe_curves_section}

Let us consider $\jtil$ as in~\eqref{alm_cpx_str_def} induced by some $J \in \J(\xi,d\lambda)$.

\begin{defn}[Hofer]\label{fe_curve_def}
Let $(S,j)$ be a closed Riemann surface and $\Gamma\subset S$ a finite set. A finite-energy curve is a map $\util: S\setminus \Gamma \to \R\times M$ satisfying the Cauchy-Riemann equations
\begin{equation}\label{cr_eqns}
  \bar\partial_{\jtil}(\util) = \frac{1}{2} \left( d\util + \jtil(\util) \cdot d\util \cdot j \right) = 0
\end{equation}
and an energy condition $0<E(\util)<\infty$. The energy is defined as
\begin{equation*}\label{}
E(\util) = \sup_{\phi \in \Lambda} \int_{S\setminus \Gamma} \util^*d(\phi\lambda)
\end{equation*}
where $\Lambda = \{ \phi :\R\to\R \text{ smooth} \mid 0\leq \phi\leq 1, \ \phi'\geq 0 \}$.
\end{defn}

Solutions of~\eqref{cr_eqns} are called $\jtil$-holomorphic maps. Each integrand in the definition of the energy is non-negative, and a quick calculation shows $\util$ is constant when $E(\util)=0$. The elements of $\Gamma$ are the so-called punctures.

\begin{remark}[Cylindrical coordinates]\label{cyl_coord}
Fix $z\in\Gamma$ and choose a holomorphic chart $\psi:(U,z) \to (\psi(U),0)$, where $U$ is a neighborhood of~$z$. We identify $[s_0,+\infty) \times \R/\Z$ with a punctured neighborhood of $z$ via $(s,t) \simeq \psi^{-1}(e^{-2\pi(s+it)})$, for $s_0 \gg 1$, and call $(s,t)$ positive cylindrical coordinates centered at $z$. We may also identify $(s,t) \simeq \psi^{-1}(\est)$ where $s<-s_0$ and, in this case, $(s,t)$ are negative coordinates. In both cases we write $\util(s,t) = \util\circ\psi^{-1}(e^{-2\pi(s+it)})$ or $\util(s,t) = \psi^{-1}(e^{2\pi(s+it)})$.
\end{remark}

Let $(s,t)$ be positive cylindrical coordinates centered at some $z\in\Gamma$, and write $\util(s,t) = (a(s,t),u(s,t))$. $E(\util)<\infty$ implies
\begin{equation}\label{mass_defn}
  m = \lim_{s\to+\infty} \int_{\{s\}\times \R/\Z} u^*\lambda
\end{equation}
exists. This number is the mass of $\util$ at $z$, and does not depend on the choice of coordinates. The puncture $z$ is called positive, negative or removable when $m>0$, $m<0$ or $m=0$ respectively, and $\util$ can be smoothly extended to $(S\setminus\Gamma)\cup\{z\}$ when $z$ is removable. Moreover, $a(s,t) \to \epsilon\infty$ as $s\to+\infty$, where $\epsilon$ is the sign of $m$.

\subsubsection{Asymptotic operators and asymptotic behavior}\label{asymptotic_behavior_section}

Let $P=(x,T)$ be a closed Reeb orbit and consider the bundle $\xi_P = x_T^*\xi \to \R/\Z$, where $x_T(t) = x(Tt)$. Let $k\in \Z^+$ be its multiplicity, that is, $T = kT_{min}$ where $T_{min}$ is minimal positive period of $x$. A choice of $J \in \J(\xi,d\lambda)$ induces an inner-product $$ \left< \eta,\zeta \right> = \int_0^1 (d\lambda)_{x(Tt)}(\eta(t),J_{x(Tt)}\zeta(t)) dt $$ of pair of sections $\eta,\zeta$ of $\xi_P$. The corresponding space of square integrable sections is denoted by $L^2_J(\xi_P)$.

\begin{defn}\label{asymp_op_def}
The asymptotic operator at $P$ is the unbounded self-adjoint operator $A_P$ on $L^2_J(\xi_P)$ defined by
\begin{equation}
  A_P(\eta) = J(-\nabla_t\eta + T\nabla_\eta R)
\end{equation}
where $\nabla$ is any symmetric connection\footnote{$A_P$ does not depend on $\nabla$.} on $M$ and $\nabla_t$ denotes covariant derivative along the curve $x_T(t)$.
\end{defn}

\begin{remark}\label{rmk_asymp_op}
A $d\lambda$-symplectic frame for $\xi_P$ presents $A_P$ as an operator of the form $-J(t)\partial_t -S(t)$, where $J(t)$ is the representation of $J_{x(Tt)}$ and $S(t)$ is symmetric with respect to $\left<\cdot,-iJ(t)\cdot\right>$. If the frame is $(d\lambda,J)$-unitary then the  $A_P$ takes the form $-i\partial_t-S(t)$ where $S(t)$ is symmetric. It follows that $A_P$ has all the spectral properties described in~\ref{cz_operators_description}.
\end{remark}

Consider coordinates $(\theta,x,y) \in \R/\Z \times \R^2$ and the contact form $\alpha_0 = d\theta + xdy$.

\begin{defn}\label{defn_martinet}
A Martinet tube for a simply covered $T$-periodic orbit $P$ is a pair $(U, \Psi)$ where $U$ is a neighborhood of $x(\R)$ in $M$ and $\Psi:U \to \R/\Z \times B$ is a diffeomorphism ($B \subset \R^2$ is an open ball centered at the origin) satisfying
\begin{itemize}
  \item[(1)] $\Psi^*(f\alpha_0) = \lambda$ where the smooth function $f:\R/\Z \times B \to \R^+$ satisfies $f|_{\R/\Z \times 0} \equiv T$ and $df|_{\R/\Z \times 0} \equiv 0$.
  \item[(2)] $\Psi(x(Tt)) = (t,0,0)$.
\end{itemize}
\end{defn}

\begin{remark}\label{rmk_martinet}
According to~\cite{props1} Martinet tubes always exist. We note that if $P$ is simply covered and $\eta$ is any non-vanishing section of $\xi_P$ then $(U,\Psi)$ can be constructed so that $\Psi_*\eta = \partial_x$.
\end{remark}

Let $\util = (a,u)$ be as in Definition~\ref{fe_curve_def}, $z\in \Gamma$ be a non-removable puncture, and $(s,t)$ be positive cylindrical coordinates centered at $z$. Let $m$ be the mass of $\util$ at $z$ and $\epsilon$ be its sign.

\begin{defn}\label{non_deg_puncture}
We call $z$ a non-degenerate puncture of $\util$ if there exists a periodic Reeb orbit $P=(x,T)$ and constants $c,d\in \R$ such that
\begin{itemize}
  \item[(1)] $\sup_{t\in \R/\Z} |a(s,t) - \epsilon Ts-d| \to 0$ as $s\to+\infty$.
  \item[(2)] $u(s,t) \to x(\epsilon Tt+c)$ in $C^0(\R/\Z,M)$ as $s\to+\infty$.
  \item[(3)] If $\pi\cdot du$ does not vanish identically then $\pi\cdot du(s,t) \neq0$ when $s\gg1$.
  \item[(4)] If we write $u(s,t) = \exp_{x(\epsilon Tt+c)}(\zeta(s,t))$ for sufficiently large values of $s$ then $\sup_{t\in \R/\Z} e^{bs}|\zeta(s,t)| \to 0$ as $s\to+\infty$, for some $b>0$.
\end{itemize}
We also say that $\util$ has non-degenerate asymptotic behavior at $z$.
\end{defn}

The above definition is independent of the choice of $\psi$ and $\exp$. If $P$ is as above then we say $\util$ is asymptotic to $P$ at $z$. When every non-removable puncture is non-degenerate we simply say $\util$ has non-degenerate asymptotics.

Here is a partial description of the asymptotic behavior from~\cite{props1}.

\begin{theorem}[Hofer, Wysocki and Zehnder]\label{partial_asymptotics}
If $\lambda$ is non-degenerate then every finite-energy $\jtil$-holomorphic curve, for any $J \in \J(\xi,d\lambda)$, has non-degenerate asymptotics.
\end{theorem}

A much more precise description of the asymptotic behavior is given in~\cite{props1}. Let $\util=(a,u)$ be asymptotic to $P = (x,T)$ at the non-degenerate puncture $z$. Consider cylindrical holomorphic coordinates $(s,t)$ centered at $z$, positive if $z$ is positive or negative if $z$ is negative. Choose a Martinet tube $(U,\Psi)$ for the underlying prime orbit, as described above. The frame $\{e_1 \simeq \partial_x/\sqrt{f}, e_2 \simeq (-x\partial_\theta + \partial_y)/\sqrt{f}\}$ of $\xi$ on $U$ is $d\lambda$-symplectic and can be used to represent $A_P \simeq L=-J(t)\partial_t - S(t)$, as explained in Remark~\ref{rmk_asymp_op}. The functions
\[
\begin{array}{cccc}
  a(s,t) \in \R, & \theta(s,t)\in \R/\Z & \text{and} & z(s,t) =(x(s,t),y(s,t))\in \R^2
\end{array}
\]
given by $(id_\R\times\Psi)\circ \util = (a,\theta,x,y)$ are well-defined for $|s|\gg1$. Let $\epsilon=\pm1$ be the sign of the puncture $z$.

The following theorem gives a precise description of the asymptotic behavior of a finite-energy surface near a non-removable puncture when the contact form is non-degenerate. It incorporates a refinement due to R. Siefring~\cite{cpam} of the original asymptotic formula of Hofer, Wysocki and Zehnder~\cite{props1}; see also E.~Mora's dissertation~\cite{mora}.

\begin{theorem}\label{precise_asymptotics}
Assume that $\lambda$ is non-degenerate and that $\int u^*d\lambda>0$. One finds $b>0$, $\mu \in \sigma(L)$ such that $\epsilon\mu<0$, an eigenvector $v$ of $L$ satisfying $Lv=\mu v$, and constants $c,d\in \R$ such that
\[
\lim_{|s|\to\infty} \sup_t e^{b|s|} \left( |D^\beta[a- Ts-d]| + |D^\beta[\theta- kt-c]| \right) = 0
\]
for every multi-index $\beta$, and
\[
z(s,t) = e^{\mu s} \left( v(t) + \Delta(s,t) \right)
\]
where $\lim_{|s|\to\infty} \sup_t |D^\beta\Delta| =0 \ \forall \beta$. Here $k\geq 1$ is the multiplicity of $P$.
\end{theorem}

The eigenvalue $\mu \in \sigma(A_P)$ and the eigenvector of $A_P$ corresponding to $v$ as in the above statement will be loosely referred to as the asymptotic eigenvalue and asymptotic eigenvector of $\util$ at $z$, respectively.

In~\cite{props4} it is proved that the conclusions of Theorem~\ref{precise_asymptotics} also hold under the assumption that $\lambda$ is only Morse-Bott. In~\cite{hry} we prove the following lemma.

\begin{lemma}
Let $\lambda$ be any contact form on $M$, fix any $J \in \J(\xi,d\lambda)$, and let $\util$ be a finite-energy $\jtil$-holomorphic curve in $\R\times M$ asymptotic to $P$ at the non-degenerate puncture $z$, in the sense of Definition~\ref{non_deg_puncture}. Then the conclusions of Theorem~\ref{precise_asymptotics} are true.
\end{lemma}

\subsubsection{Some algebraic invariants}\label{section_alg_inv}

We need to recall a few definitions from~\cite{props2}. Let $J \in \J(\xi,d\lambda)$ induce $\jtil$ as in~\eqref{alm_cpx_str_def}, $(S,j)$ be a closed Riemann surface, and $\util=(a,u):S\setminus \Gamma \to \R\times M$ be a finite-energy $\jtil$-holomorphic curve, where $\Gamma \subset S$ is a finite set of non-removable punctures. Assume also that every $z\in\Gamma$ is a non-degenerate puncture as in Definition~\ref{non_deg_puncture}.

The section $\pi \cdot du$ of $\mathcal E = \wedge^{0,1}T^*(S\setminus\Gamma)\otimes_\C u^*\xi$ satisfies the Cauchy-Riemann type equation
\begin{equation*}\label{}
  \pi \cdot du \cdot j = J(u) \cdot \pi \cdot du.
\end{equation*}
This follows from~\eqref{cr_eqns}. Thus either $\pi \cdot du$ vanishes identically or its zeros are isolated. In the second case one defines
\begin{equation}\label{wind_pi_defn}
  \wind_\pi(u) = \text{ algebraic count of zeros of } \pi\cdot du
\end{equation}
where $\mathcal E$ is oriented by its natural complex structure.

\begin{remark}[Winding numbers]\label{rmk_winding_numbers}
Let $E\to \R/\Z$ be an oriented rank-2 real vector bundle, and consider two non-vanishing sections $Z$ and $W$ of $E$. A choice of complex structure $J$ inducing the orientation of $E$ gives unique functions $a,b: \R/\Z\to\R$ satisfying $W = aZ + bJZ$. The function $f=a+ib: \R/\Z\to\C$ does not vanish and we define
\begin{equation}\label{}
  \wind(W,Z) = \deg \frac{f}{|f|} \in \Z.
\end{equation}
This integer depends only on the homotopy class of non-vanishing sections of $Z$ and $W$, and on the orientation of $E$. When $E$ is symplectic, we use the induced orientation.
\end{remark}

Let $Z$ be a non-vanishing section of $u^*\xi$ and assume $\pi\cdot du$ does not vanish identically. If $(s,t)$ are positive cylindrical coordinates centered at some $z\in\Gamma$ then set $$ \wind_\infty(\util,z,Z) := \lim_{s\to+\infty} \wind(t\mapsto \pi\cdot \partial_s u(s,\epsilon t), t\mapsto Z(s,\epsilon t)) $$ where $\epsilon = \pm 1$ is the sign of the puncture $z$. Splitting $\Gamma = \Gamma^+ \sqcup \Gamma^-$ into positive and negative punctures, Hofer, Wysocki and Zehnder define in~\cite{props2}
\begin{equation}\label{wind_infty_defn}
  \wind_\infty(\util) = \sum_{z\in\Gamma^+} \wind_\infty(\util,z,Z) - \sum_{z\in\Gamma^-} \wind_\infty(\util,z,Z).
\end{equation}
Clearly this number does not depend on $Z$.

\begin{defn}\label{fast_planes}
A finite-energy plane $\util:\C \to \R \times M$ will be called fast if $\infty$ is a non-degenerate puncture, the asymptotic limit $P$ of $\util$ is a simply covered (prime) Reeb orbit and $\wind_\pi(\util)=0$.
\end{defn}

We do not make any non-degeneracy assumptions on $\lambda$ or $P$ in the above definition.

An application of the argument principle proves the following important identity.

\begin{lemma}[Hofer, Wysocki and Zehnder]\label{lemma_wind_pi_infty}
$\wind_\pi(\util) = \wind_\infty(\util) - \chi(S) + \#\Gamma$.
\end{lemma}

The following statement will be left without proof.

\begin{lemma}\label{winding_asym_technical}
Let $\util=(a,u)$ be a $\jtil$-holomorphic finite-energy surface, $z$ be a non-degenerate puncture of $\util$ with sign $\epsilon=\pm1$, and $P=(x,T)$ be the asymptotic limit of $\util$ at $z$. Fix holomorphic cylindrical coordinates $(s,t)$ centered at $z$ with sign $\epsilon$. For any smooth non-vanishing section $Z$ of $x_T^*\xi$ there is a smooth section $\overline Z$ of $u^*\xi$ defined near $z$ such that $\overline Z(s,t) \to Z(t)$ in $C^0(\R/\Z,\xi)$ as $\epsilon s\to+\infty$.
\end{lemma}

\begin{remark}\label{asymp_evalue_wind_infty}
Suppose $\lambda$ is non-degenerate and that the finite-energy plane $\util$ is asymptotic to the closed Reeb orbit $P=(x,T)$. It follows from Theorem~\ref{precise_asymptotics}, Lemma~\ref{winding_asym_technical} and~\eqref{wind_infty_defn} that the winding of the asymptotic eigenvector with respect to a trivialization of $x_T^*\xi$ induced by the capping disk given by $\util$ is precisely $\wind_\infty(\util)$.
\end{remark}

\subsubsection{Curves with vanishing $d\lambda$-energy}

Let $\Gamma \subset S$ be finite, $J\in\J(\xi,d\lambda)$ and $\util=(a,u):S\setminus \Gamma \to \R\times M$ be a finite-energy $\jtil$-holomorphic surface. The integral
\begin{equation}\label{dlambda-energy}
  \int_{S\setminus \Gamma} u^*d\lambda
\end{equation}
is non-negative, bounded by $E(\util)$ and vanishes if, and only if, $\pi\cdot du\equiv0$.

\begin{theorem}[Hofer, Wysocki and Zehnder]\label{zera_dlambda_theorem}
If $\Gamma \subset \C$ is finite, $\util:\C\setminus\Gamma \to \R\times M$ is as above, $\Gamma$ consists of negative punctures and $\int_{\C\setminus \Gamma} u^*d\lambda=0$, then one finds a non-constant polynomial $p:\C\to \C$ and a closed Reeb orbit $P=(x,T)$ such that $p^{-1}(0) = \Gamma$ and $\util = F_P \circ p$. Here $F_P:\C\setminus\{0\} \to \R\times M$ is defined by $F_P(\est) = (Ts,x(Tt))$.
\end{theorem}

\begin{cor}
If $\Gamma = \emptyset$ then $\int_\C u^*d\lambda>0$.
\end{cor}

\subsubsection{Bubbling-off points}

The basic tool for the bubbling-off analysis is the following lemma, where norms are taken with respect to any $\R$-invariant metric on $\R\times M$ and the euclidean metric on $\C$.

\begin{lemma}\label{beforeclaim}
Let $\Gamma \subset \C$ be finite and $U_n \subset \C\setminus\Gamma$ be an increasing sequence of open sets such
that $\cup_n U_n = \C\setminus\Gamma$. Let $\util_n = (a_n,u_n) : (U_n,i) \rightarrow (\R \times M,\jtil)$ be $\jtil$-holomorphic maps satisfying $\sup_n E(\util_n) = C < \infty$, and $z_n \in U_n$ be a sequence
such that $|d\util_n(z_n)| \rightarrow +\infty$. If $z_n$ stays bounded away from $\Gamma \sqcup \{\infty\}$,
or if some $U_m$ is a neighborhood of $\infty$ and $z_n$ stays bounded away from
$\Gamma$, then there exist subsequences $\{\util_{n_j}\}$ and $\{z_{n_j}\}$, sequences $z^\prime_j \in \C$ and $r_j \in\R$, and a contractible periodic Reeb orbit $\hat P = (\hat x , \hat T)$ such that $|z_{n_j}-z_j^\prime|\rightarrow 0$, $r_j \rightarrow 0^+$, $\hat T \leq C$ and
\[
 \limsup_{j\rightarrow+\infty} \int_{|z-z_j^\prime|\leq r_j} u_{n_j}^*d\lambda \geq \hat T.
\]
\end{lemma}

The proof is standard and will be omitted.

\section{Passing to the degenerate case}\label{deg_case}

Our goal here is to prove Theorem~\ref{main}. We assume, without loss of generality, that the dynamically convex contact form $\lambda$ on $S^3$ is of the form $f\lambda_0|_{S^3}$, where $f:S^3 \to (0,+\infty)$ is smooth and $\lambda_0=\frac{1}{2} \sum qdp - p dq$ is the standard Liouville form on $\R^4$. As before, the Reeb vector is denoted by $R$, its flow by $\phi_t$ and the standard contact structure on $S^3$ by $\xi = \ker\lambda_0$. From now on we assume $\bar P=(\bar x,\bar T)$ is a closed Reeb orbit of $\lambda$ as in the statement of Theorem~\ref{main}.

\subsection{A suitable spanning disk for $\bar P$}\label{perturbing}

Recall that on any embedded oriented surface $S\subset S^3$ there is a singular distribution
\begin{equation}\label{}
  (\xi \cap TS)^\bot
\end{equation}
called the characteristic distribution of $S$, where here $\bot$ denotes the $d\lambda$-symplectic orthogonal (with respect to any defining contact form $\lambda$ for $\xi$). It equals $\xi \cap TS$ except at the singular points where $\xi = TS$ and $(\xi \cap TS)^\bot = \{0\}$. If $S$ is in a regular level set of some function $H$ then equations
\begin{equation}\label{}
  \begin{array}{ccc}
    i_V\lambda = 0, &  & i_Vd\lambda = dH - (i_RdH)\lambda
  \end{array}
\end{equation}
define a vector field $V$ on $S$ that parametrizes the characteristic distribution. A singular point $p\in \mathcal D$ is called non-degenerate when the linearization $DV_p$ is an isomorphism. Then $p$ is called elliptic or hyperbolic if the determinant of $DV_p$ is positive or negative, respectively. The space $T_pS = \xi|_p$ has two orientations: $o_p$ induced by the given orientation of $S$, and $o'_p$ induced by $d\lambda$. The singular point $p$ is positive if $o_p=o'_p$, and negative otherwise. Following Hofer~\cite{93}, we call $p$ nicely elliptic if it is elliptic and the eigenvalues of $DV_p$ are real.

Below we only consider the case where $\partial S\neq \emptyset$ is a knot transverse to $\xi$ and, in this case, $S$ will always be oriented by $\lambda_{T\partial S} >0$. Applying arguments from~\cite{char1,char2} (which use Giroux's elimination lemma from~\cite{giroux1}) to our situation one proves

\begin{theorem}[Hofer, Wysocki and Zehnder]\label{special_disk_0}
Let $L \subset (S^3,\xi)$ be a transverse unknot spanned by an embedded disk $\mathcal D_0 \subset S^3$ satisfying $sl(L) = -1$. Then there exists another embedded disk $\mathcal D_1$ spanning $L$ such that the characteristic distribution of $\mathcal D_1$ has precisely one singularity, which is a positive nicely elliptic point. Moreover, the disk $\mathcal D_1$ can be obtained by a smooth and arbitrarily $C^0$-small perturbation of $\mathcal D_0$ supported away from $L$.
\end{theorem}

In order to consider special Bishop families of pseudo-holomorphic disks we need suitable boundary conditions provided by the following statement.

\begin{lemma}\label{special_disk}
Consider a sequence of smooth functions $h_k:S^3 \to \R^+$ satisfying $h_k \equiv 1$, $dh_k \equiv 0$ on $\bar x(\R)$ and $h_k \to 1$ in $C^\infty$. There exists an embedded disk $\mathcal D \subset S^3$ spanning $\bar x(\R)$ for which we can find $k_0$ and a neighborhood $\O$ of $\partial\mathcal D = \bar x(\R)$ in $\mathcal D$ such that $\R R_k|_p \cap T_p\mathcal D = 0$, $\forall \ p \in \O\setminus \bar x(\R)$ and $k\geq k_0$. Here $R_k$ denotes the Reeb vector of $h_k\lambda$.
\end{lemma}

Note that, by the assumptions on $h_k$, $t\mapsto \bar x(t)$ is $\bar T$-periodic trajectory of the vector fields $R_k$.

\begin{proof}
Let $\mathcal D_0 \subset S^3$ be an embedded disk spanning $\bar x(\R)$ and $$ t \in \R/\Z \mapsto W(t) \in (T\mathcal D_0 \cap \xi)|_{\bar x(\bar Tt)} $$ be a non-vanishing vector along the curve ${\bar x}_{\bar T}$. Let us denote $\lambda_k = h_k\lambda$. We claim that there exists an embedding
\begin{equation}\label{embedding_psi}
  \psi: (1-\epsilon,1] \times \R/\Z \to S^3
\end{equation}
and an integer $k_0\geq 1$ satisfying the following properties:
\begin{itemize}
  \item[(i)] $\psi(1,t) = \bar x(\bar Tt)$.
  \item[(ii)] Denoting $S = \psi((1-\epsilon,1] \times \R/\Z)$, let $t\mapsto N(t)$ be a non-vanishing vector field  satisfying $N(t) \in (TS \cap \xi)|_{\bar x(\bar T t)}, \ \forall t$. Then $\wind(N(t),W(t)) = 0$.
  \item[(iii)] If $k\geq k_0$ then $\{\partial_r\psi,\partial_t\psi,R_k\circ\psi\}$ is a basis of $T_{\psi(r,t)}S^3$ when $1-\epsilon < r< 1$.
\end{itemize}

Let $U$ be a small neighborhood of $\bar x(\R)$ in $S^3$ and $\Phi : U \to \R/\Z\times B$ be a diffeomorphism, where $B \subset \R^2$ is the unit ball centered at the origin, satisfying $\Phi(\bar x(\bar T\theta)) = (\theta,0,0)$ $\forall \theta \in \R/\Z$, $\Phi_*\lambda|_{\R/\Z\times (0,0)} = \bar Td\theta$ and $\Phi_*(d\lambda)|_{\R/\Z\times (0,0)} \equiv dx \wedge dy$. Here $(\theta,x,y)$ are standard coordinates in $\R/\Z\times \R^2$. Note that $\Phi_*W(t) \in 0\times \R^2$ for every $t$. The map $\Phi$ can also be arranged so that $\Phi_*W(t)$ and $\partial_x|_{(t,0,0)}$ do not wind with respect to each other.

In these coordinates the linearized $\lambda$-Reeb flow along $\bar x(\R) \simeq \R/\Z\times (0,0)$ is represented as
\begin{equation*}\label{}
  d\phi_{\bar Tt}|_{\bar x(0)} \simeq \begin{pmatrix} 1 & \\ & \varphi(t) \end{pmatrix}
\end{equation*}
where $\varphi:[0,1] \to Sp(1)$ is a path of symplectic matrices satisfying $\varphi(0) = I$. Note that $$ \mu(\varphi) = \mu_{CZ}(\bar P) + 2 sl(\bar P) \geq 1 $$ where $sl(\bar P)=-1$ is the self-linking number of $\bar P$. Let $I(\varphi)$ be the winding interval of $\varphi$ defined as in~\ref{cz_geom_def}. In view of the geometric description of the $\mu$-index we know $I(\varphi) \subset \R^+$ and $\varphi(1) = I \Leftrightarrow I(\varphi) \subset \Z^+$. One of the following cases hold:
\begin{enumerate}
  \item[(a)] $\sigma(\varphi(1)) = \{a,a^{-1}\}$ with $a>0$, $a\neq 1$. In this case we set $$ Y = \begin{pmatrix} \ln a & 0 \\ 0 & -\ln a \end{pmatrix} $$ and find $T\in Sp(1)$ such that $\varphi(1) = T^{-1} e^Y T$.
  \item[(b)] $\sigma(\varphi(1)) = \{1\}$ and $\varphi(1) \neq I$. In this case we find $a\in \R\setminus \{0\}$ and $T\in Sp(1)$ such that $\varphi(1) = T^{-1}e^YT$ where $$ Y = \begin{pmatrix} 0 & a \\ 0 & 0 \end{pmatrix}. $$
  \item[(c)] $\sigma(\varphi(1)) = \{e^{\pm i\gamma}\}$ for some $\gamma \in (0,2\pi)\setminus \{\pi\}$, or $\varphi(1)=-I$ with $\gamma = \pi$, or $\varphi(1)=I$ with $\gamma = 0$. In this case we set $$ Y = \begin{pmatrix} 0 & -\gamma \\ \gamma & 0 \end{pmatrix} $$ and, perhaps after changing $\gamma$ by $2\pi-\gamma$, we find $T\in Sp(1)$ such that $\varphi(1) = T^{-1} e^Y T$.
  \item[(d)] $\sigma(\varphi(1)) = \{a,a^{-1}\}$ with $a<0$, $a\neq -1$. In this case we set $$ Y = \begin{pmatrix} \ln (-a) & 0 \\ 0 & -\ln (-a) \end{pmatrix} $$ and find $T\in Sp(1)$ such that $-\varphi(1) = T^{-1} e^Y T$.
  \item[(e)] $\sigma(\varphi(1)) = \{-1\}$ and $\varphi(1) \neq -I$. In this case we find $a$ and $Y$ as in case (b) such that $-\varphi(1) = T^{-1} e^Y T$.
\end{enumerate}

After composing $\Phi$ with the linear diffeomorphism $$ \begin{pmatrix} \theta \\ x \\ y \end{pmatrix} \mapsto \begin{pmatrix} \theta \\ T\begin{pmatrix} x \\ y \end{pmatrix} \end{pmatrix} $$ in each case, we could have assumed that $\Phi$ satisfies all the properties mentioned before and, moreover, that $\varphi(1) = e^Y$ in cases (a)-(c) or $\varphi(1) = -e^Y$ in case (d)-(e). In cases (a), (b) and (c) we set $K(t) = e^{tY}$, and in case (d)-(e) we set\footnote{We may denote the matrix $\begin{pmatrix} \cos y & -\sin y \\ \sin y & \cos y \end{pmatrix}$ by $e^{iy}$.} $K(t) = e^{i\pi t}e^{tY}$, so that $K(1) = \varphi(1)$ in all cases.

We need to understand in detail the index $\mu(K)$ and the winding interval $I(K)$ of the path $t\mapsto K(t)$.
\begin{itemize}
  \item In case (a): $I(K)$ contains $0$ in its interior and $\mu(K) = 0$.
  \item In case (b): if $a<0$ then $I(K) = [0,c]$ for some $0<c<1/2$ and $\mu(K) = 0$, if $a>0$ then $I(K) = [c,0]$ for some $-1/2<c<0$ and $\mu(K) = -1$.
  \item In case (c): $I(K) = \{\gamma/2\pi\} \subset [0,1)$. If $\gamma>0$ then $\mu(K) = 1$, if $\gamma=0$ then $\mu(K) = -1$.
  \item In case (d): $I(K)$ contains $1/2$ in its interior and $\mu(K) = 1$.
  \item In case (e): $I(K)$ contains $1/2$ in its boundary and $\mu(K) = 1$.
  \end{itemize}
We shall now construct an embedding~\eqref{embedding_psi} satisfying conditions (i) and (ii) above which is transverse to $R$ in $\psi((1-\epsilon,1)\times \R/\Z)$ in each case separately. Then, after this is done, we will check condition (iii). For simplicity we will assume that $\bar T=1$ in the following, without loss of generality. \\

\noindent {\bf Case (a).} The loop
\begin{equation}\label{loop_M}
  M(t) = K(t)\varphi^{-1}(t)
\end{equation}
satisfies $\maslov(M) = -k$ where $\mu(\varphi) = 2k$ for some $k\geq 1$. Note that $t\in \R/\Z \mapsto M(t) \in \sp(1)$ is a smooth map (since so is $\varphi'\varphi^{-1}$). We still write $(\theta,x,y)$ for the new coordinates obtained by composing $\Phi$ with the diffeomorphism
\begin{equation}\label{diffeo_H}
  H: \begin{pmatrix} \theta \\ x \\ y \end{pmatrix} \mapsto \begin{pmatrix} \theta \\ M(\theta)\begin{pmatrix} x \\ y \end{pmatrix} \end{pmatrix}.
\end{equation}
The linearized $\lambda$-Reeb flow $d\phi_{\bar Tt}: \xi|_{\bar x(0)} \to \xi|_{\bar x(\bar Tt)}$ along $\bar P$ is now represented by
\begin{equation*}\label{}
   \begin{pmatrix} 1 & \\ & K(t) \end{pmatrix}
\end{equation*}
and the differential of the $\lambda$-Reeb vector by
$$ DR(\theta,0,0) = \begin{pmatrix} 0 & \\ & K'K^{-1} = Y \end{pmatrix} = \begin{pmatrix} 0 & 0 & 0 \\ 0 & \ln a & 0 \\ 0 & 0 & -\ln a \end{pmatrix}. $$

We consider $\epsilon>0$ small and a map $F:(1-\epsilon,1] \times \R/\Z \to \R/\Z\times B$ of the form
\begin{equation}\label{strip_F}
  F(r,t) = (t,(1-r)\cos b(t),(1-r)\sin b(t))
\end{equation}
for some smooth function $b:\R\to \R$ satisfying $b(t+1) = b(t) - 2\pi k$. Before describing $b(t)$ in detail we make some \emph{a priori} computations. First, the Reeb vector can be written as
\begin{equation}\label{aprox_R}
  \begin{aligned}
    R\circ F(r,t) &= R(t,0,0) + (r-1) DR(t,0,0) \partial_rF(1,t) + O(|1-r|^2) \\
    &= \begin{pmatrix} 1 \\ (1-r)Y\begin{pmatrix} \cos b(t) \\ \sin b(t) \end{pmatrix} \end{pmatrix} + O(|1-r|^2)
  \end{aligned}
\end{equation}
for $r\sim 1$. Consider the function
\begin{equation}\label{function_d}
  d(r,t) = \det (\partial_r F(r,t),\partial_t F(r,t),R\circ F(r,t)).
\end{equation}
Then
\begin{equation*}\label{}
  \begin{aligned}
    d(r,t) &= \det \begin{pmatrix} 0 & 1 & 1 \\ -\cos b & (r-1)b'\sin b & (1-r)(\ln a) \cos b \\ -\sin b & (1-r)b'\cos b & (r-1)(\ln a) \sin b \end{pmatrix} + O(|1-r|^2) \\
    &= -(1-r)(b'+(\ln a)\sin 2b) + O(|1-r|^2)
  \end{aligned}
\end{equation*}
for $r\sim 1$. We found
\begin{equation}\label{}
  \partial_r d(1,t) = b'+(\ln a)\sin 2b.
\end{equation}
Since $-2\pi k\leq -2\pi$ the function $b\in[-2\pi k,0] \mapsto \sin 2b$ changes sign, and we find an non-empty open interval $J\subset (-2\pi k,0)$ such that $b\in J \Rightarrow (\ln a)\sin 2b < 0$. If $\alpha<0$ satisfies $|\alpha| > |\ln a|$ then we may take $b(t)$ satisfying $b'<0$ and $b' \geq \alpha \Leftrightarrow b \in J$. Thus $b'(t) +(\ln a)\sin 2b(t) < 0$ for every $t\in \R$, and $d(r,t) \neq 0$ if $1-r$ is small and positive. It follows that the Reeb vector is transverse to the embedded (open) strip $F((1-\epsilon,1)\times \R/\Z)$ when $\epsilon$ is small. Setting
\begin{equation}\label{map_h}
  h(r,t) = H^{-1} \circ F(r,t)
\end{equation}
then $$ \partial_r h(1,t) = \begin{pmatrix} 0 \\ M(t)^{-1} \begin{pmatrix} -\cos b(t) \\ -\sin b(t) \end{pmatrix} \end{pmatrix} $$ and the non-vanishing vector $t\mapsto d\pi_0 \cdot \partial_r h(1,t)$ has winding number equal to $0$ around the origin in $\R^2$. Here $\pi_0 : \R/\Z\times \R^2 \to \R^2$ is the projection onto the second factor. This follows since $\maslov(M^{-1}) = k$ and $b(1)-b(0) = -2\pi k$. Setting $\psi(r,t) = \Phi^{-1} \circ h(r,t)$, $\psi$ satisfies the required conditions. \\

\noindent {\bf Case (b).} Assuming $a>0$ we have $\mu(\varphi) = 2k+1$ for some $k\geq 0$ and the loop $$ M(t) = e^{i2\pi t}K(t)\varphi^{-1}(t) $$ satisfies $\maslov(M) = -k$. As before, $t\in \R/\Z \mapsto M(t) \in \sp(1)$ is smooth. Composing with the diffeomorphism $H$~\eqref{diffeo_H} we obtain new coordinates, still denoted $(\theta,x,y)$, where the linearized $\lambda$-Reeb flow $d\phi_{\bar Tt}: \xi|_{\bar x(0)} \to \xi|_{\bar x(\bar Tt)}$ along $\bar P$ is represented by $e^{i2\pi t}K(t)$ and the differential of the $\lambda$-Reeb vector by $$ DR(\theta,0,0) = \begin{pmatrix} 0 & \\ & i2\pi + e^{i2\pi t}Ye^{-i2\pi t} \end{pmatrix} $$ where $$ i2\pi + e^{i2\pi t}Ye^{-i2\pi t} = \begin{pmatrix} -a\sin(2\pi t)\cos(2\pi t) & -2\pi + a\cos^2(2\pi t) \\  2\pi - a\sin^2(2\pi t) & a\sin(2\pi t)\cos(2\pi t) \end{pmatrix} = S(t). $$
Set $F(r,t)$ as in~\eqref{strip_F} where the smooth function $b:\R\to \R$ satisfying $b(t+1) = b(t) - 2\pi k$ and b(0)=0 will be defined \emph{a posteriori}. We compute similarly to~\eqref{aprox_R}
\[
  \begin{aligned}
    R\circ F(r,t) &= R(t,0,0) + (r-1) DR(t,0,0) \partial_rF(1,t) + O(|1-r|^2) \\
    &= \begin{pmatrix} 1 \\ (1-r)S(t)\begin{pmatrix} \cos b(t) \\ \sin b(t) \end{pmatrix} \end{pmatrix} + O(|1-r|^2).
  \end{aligned}
\]
If $d(r,t)$ is the function defined as in~\eqref{function_d} then
\begin{equation*}\label{}
  \begin{aligned}
    d(r,t) &= \det \begin{pmatrix} 0 & 1 & 1 \\ -\cos b & (r-1)b'\sin b & (1-r)(Se^{ib})_{11} \\ -\sin b & (1-r)b'\cos b & (1-r)(Se^{ib})_{21} \end{pmatrix} + O(|1-r|^2) \\
    &= \det \begin{pmatrix} 0 & 1 & 1 \\ -1 & 0 & (1-r)(e^{-ib}Se^{ib})_{11} \\ 0 & (1-r)b' & (1-r)(e^{-ib}Se^{ib})_{21} \end{pmatrix} + O(|1-r|^2) \\
    &= -(1-r)[b' - (e^{-ib}Se^{ib})_{21}] + O(|1-r|^2) \\
    &= -(1-r)[b'-2\pi + a\sin^2(2\pi t-b(t))] + O(|1-r|^2) \\
    &= -(1-r)[\beta' + a\sin^2 \beta(t)] + O(|1-r|^2)
  \end{aligned}
\end{equation*}
where $\beta(t) = b(t)-2\pi t$. Let $\alpha<0$ satisfy $|\alpha|\gg a$. In view of the condition $\beta(1) = -2(k+1)\pi \leq -2\pi$ we may choose $\delta>0$ small and define $\beta(t)$ satisfying
\begin{itemize}
  \item $\beta'<0$ and $\beta'(t) \neq \alpha \Leftrightarrow \beta(t) \in [-\pi-\delta,-\pi]$.
  \item $\beta' < -a\sin^2 \beta$ when $\beta(t) \in [-\pi-\delta,-\pi]$.
\end{itemize}
Here we strongly use that $\sin^2 \beta = (\beta+\pi)^2 + O(|\beta+\pi|^4)$ when $\beta \to -\pi$. In fact, consider $C>0$ and $\delta_0>0$ small so that $|x^2-\sin^2(-\pi-x)| \leq  Cx^4$ if $0\leq x \leq\delta_0$. Take $0<\delta<\delta_0$ satisfying $C\delta^4 < \delta^2/2$ and $-\delta+3a\delta^2<0$. Let $0<t_0<t_1<1$ be defined by $\alpha t_0 =-\pi$ and $-2\pi (k+1) + \alpha(t_1-1) = -\pi-\delta$. Note that $t_0\to 0^+$ and $t_1\to 1^-$ as $\alpha \to -\infty$. There exists a smooth function $\beta:[t_0,t_1] \to \R$ satisfying
\begin{itemize}
  \item $\beta(t_0)=-\pi$ and $\beta(t_1) =-\pi -\delta$,
  \item $\beta'(t_0) = \beta'(t_1) = \alpha$ and $\beta^{(j)}(t_0) = \beta^{(j)}(t_1) = 0 \ \forall j\geq 2$ and
  \item $\beta'(t) \leq -\frac{\delta}{2(t_1-t_0)} \ \forall t\in[t_0,t_1]$.
\end{itemize}
If we extend $\beta$ to $[0,1]$ by $\beta(t) = \alpha t$ on $[0,t_0]$ and $\beta(t) = -2\pi (k+1) + \alpha(t-1)$ on $[t_1,1]$ then $\beta$ is $C^\infty$ and we can estimate for $t\in[t_0,t_1]$:
\[
  \begin{aligned}
    \beta'(t) + a\sin^2 \beta(t) & \leq -\frac{\delta}{2(t_1-t_0)} + a (\beta(t) + \pi)^2 + aC(\beta(t)+\pi)^4 \\
    & \leq -\frac{\delta}{2} + a\delta^2  + aC\delta^4 \leq -\frac{\delta}{2} + \frac{3a\delta^2}{2} <0.
  \end{aligned}
\]
If $\alpha$ is close to $-\infty$ then the same inequality is satisfied for all $t\in[0,1]$.

Thus $b'-2\pi + a\sin^2(2\pi t-b(t)) \neq 0$ and $\partial_r d(1,t) \neq 0$ for every $t$. As before we define $h(r,t)$ as in~\eqref{map_h} and $\psi = \Phi^{-1} \circ h(r,t)$.

The case $a<0$ is much simpler since $\mu(\varphi) = 2k$ for some $k\geq 1$ and we may consider the closed loop $M(t) = K(t)\varphi^{-1}(t)$ which has Maslov index $-k$. The argument is entirely analogous. \\

\noindent {\bf Case (c).} We assume $\gamma\in(0,2\pi)\setminus \{\pi\}$, the cases $\varphi(1)=\pm I$ are left to the reader. The loop $M(t) = K(t)\varphi^{-1}(t)$ has Maslov index $-k \leq 0$ in this case, where $k$ is given by $\mu(\varphi) = 2k+1 \geq 1$. This is so since $\mu(K) = 1$. Again we change coordinates by composing with the diffeomorphism~\eqref{diffeo_H}, and the linearized $\lambda$-Reeb flow $d\phi_{\bar Tt} : \xi|_{\bar x(0)} \to \xi|_{\bar x(\bar Tt)}$ becomes represented by $t\mapsto K(t)$. Consequently, the differential of the Reeb vector field is $$ DR(t,0,0) = \begin{pmatrix} 0 \\ & K'K^{-1} = i\gamma \end{pmatrix}. $$

We define $F(r,t)$ as in~\eqref{strip_F} with some smooth function $b:\R\to \R$ satisfying $b(t+1) = b(t) - 2\pi k$ to be constructed \emph{a posteriori}, and $d(r,t)$ by~\eqref{function_d}. As before we have
\[
  \begin{aligned}
    R\circ F(r,t) &= R(t,0,0) + (r-1) DR(t,0,0) \partial_rF(1,t) + O(|1-r|^2) \\
    &= \begin{pmatrix} 1 \\ (r-1)\gamma\sin b(t) \\ (1-r)\gamma\cos b(t) \end{pmatrix} + O(|1-r|^2).
  \end{aligned}
\]
Thus
\begin{equation*}\label{}
  \begin{aligned}
    d(r,t) &= \det \begin{pmatrix} 0 & 1 & 1 \\ -\cos b & (r-1)b'\sin b & (r-1)\gamma\sin b(t) \\ -\sin b & (1-r)b'\cos b & (1-r)\gamma\cos b(t) \end{pmatrix} + O(|1-r|^2) \\
    &= \det \begin{pmatrix} 0 & 1 & 1 \\ -1 & 0 & 0 \\ 0 & (1-r)b' & (1-r)\gamma \end{pmatrix} + O(|1-r|^2) \\
    &= -(1-r)(b' - \gamma) + O(|r-1|^2).
  \end{aligned}
\end{equation*}
In this case we simply set $b(t) = -2\pi kt$, so that $b'-\gamma =-2\pi k-\gamma<0$ for every $t$. Again we achieved $\partial_r d(1,t) \neq 0, \ \forall t$. As in cases (a)-(b) we set $h(r,t)$ as in~\eqref{map_h} and $\psi = \Phi^{-1} \circ h(r,t)$. \\

\noindent {\bf Case (d).} The loop $M(t) = K(t)\varphi^{-1}(t)$ has Maslov index $-k \leq 0$ where $k$ is given by $\mu(\varphi) = 2k+1 \geq 1$. This is so since $\mu(K) = 1$. Composing with the diffeomorphism~\eqref{diffeo_H} we obtain new coordinates and the linearized $\lambda$-Reeb flow $d\phi_{\bar Tt} : \xi|_{\bar x(0)} \to \xi|_{\bar x(\bar Tt)}$ becomes represented by $t\mapsto K(t)$. Consequently, the differential of the Reeb vector field is $$ DR(t,0,0) = \begin{pmatrix} 0 \\ & K'K^{-1} \end{pmatrix}. $$

We define $F(r,t)$ as in~\eqref{strip_F} with some smooth function $b:\R\to \R$ satisfying $b(t+1) = b(t) - 2\pi k$ and b(0)=0 to be constructed \emph{a posteriori}, and $d(r,t)$ by~\eqref{function_d}. As before we have
\[
  \begin{aligned}
    R\circ F(r,t) &= R(t,0,0) + (r-1) DR(t,0,0) \partial_rF(1,t) + O(|1-r|^2) \\
    &= \begin{pmatrix} 1 \\ (1-r)K'K^{-1}\begin{pmatrix} \cos b(t) \\ \sin b(t) \end{pmatrix} \end{pmatrix} + O(|1-r|^2).
  \end{aligned}
\]
Thus
\begin{equation*}\label{}
  \begin{aligned}
    d(r,t) &= \det \begin{pmatrix} 0 & 1 & 1 \\ -\cos b & (r-1)b'\sin b & (1-r)(K'K^{-1}e^{ib})_{11} \\ -\sin b & (1-r)b'\cos b & (1-r)(K'K^{-1}e^{ib})_{21} \end{pmatrix} + O(|1-r|^2) \\
    &= \det \begin{pmatrix} 0 & 1 & 1 \\ -1 & 0 & (1-r)(e^{-ib}K'K^{-1}e^{ib})_{11} \\ 0 & (1-r)b' & (1-r)(e^{-ib}K'K^{-1}e^{ib})_{21} \end{pmatrix} + O(|1-r|^2)
  \end{aligned}
\end{equation*}
Substituting
\[
  \begin{aligned}
    e^{-ib}K'K^{-1}e^{ib} &= e^{-ib}(i\pi + e^{i\pi t}Ye^{-i\pi t})e^{ib} \\
    &= \begin{pmatrix} \ln(-a) \cos (2\pi t-2b) & -\pi + \ln(-a)\sin (2\pi t-2b) \\ \pi+ \ln(-a)\sin (2\pi t-2b) & -\ln(-a)\cos (2\pi t-2b) \end{pmatrix}
  \end{aligned}
\]
we get
\begin{equation*}
  \begin{aligned}
    d(r,t) &= -(1-r)[b' - \pi - \ln(-a) \sin (2\pi t-2b)] + O(|1-r|^2) \\
    &= -(1-r)[\beta' + \ln(-a)\sin 2\beta(t)] + O(|1-r|^2)
  \end{aligned}
\end{equation*}
where $\beta(t) = b(t) - \pi t$. The condition $b(1) = -2\pi k$ forces $\beta(1) = -(2k+1)\pi \leq -\pi$, so $\sin 2\beta$ is forced to change sign. Let $J \subset (-(2k+1)\pi,0)$ be an non-empty open interval so that $\beta \in J \Rightarrow \ln(-a)\sin 2\beta < 0$, and pick a number $\alpha < 0$ satisfying $|\alpha| > |\ln(-a)|$. We can find a function $\beta(t)$ satisfying $\beta' < 0$, $\beta' \neq \alpha \Leftrightarrow \beta(t) \in J$, $\beta(0)=0$ and $\beta(1) = -(2k+1)\pi$. It follows that $\beta'(t) + \ln(-a)\sin 2\beta(t)<0$ for every $t\in [0,1]$, and that $b(t)$ can be smoothly extended to $\R$ by $b(t+1) = b(t)-2\pi k$. Finally we set $h(r,t)$ as in~\eqref{map_h} and $\psi = \Phi^{-1} \circ h(r,t)$. \\

\noindent {\bf Case (e):} This case is similar to case (b). \\

In all cases the embedding $\psi:(1-\epsilon,1]\times \R/\Z \to S^3$ was obtained by the formula $\psi(r,t) = \Phi^{-1} \circ H^{-1} \circ F(r,t)$ and we checked that $d(r,t) = \det (\partial_rF, \partial_tF, R\circ F)$ satisfies $\partial_rd(1,t) \neq 0, \ \forall t$ (here $R$ is the representation of the $\lambda$-Reeb vector in local coordinates). Now consider, for each $k$, the function $$ d_k(r,t) = \det (\partial_rF, \partial_tF, R_k\circ F). $$ Then $d_k \to d$ in $C^\infty_{loc}$ so that, in each case, we find $k_0$ large and $c>0$ satisfying $k\geq k_0 \Rightarrow |\partial_rd_k(1,t)| \geq c, \ \forall t$. Thus, possibly after taking $\epsilon$ smaller, the embedding $\psi$ is transverse to $R_k$ if $1-\epsilon<r<1$ and $k\geq k_0$, as required.

Finally note that since $\wind(N,W) = 0$ the strip $S$ can be glued with the disk $\mathcal D_0$ (away from their common boundaries) to obtain an embedded disk $\mathcal D$ with all the required properties.
\end{proof}

\begin{lemma}\label{special_disk_2}
Let $h_k:S^3 \to \R$ be a sequence of smooth functions satisfying $h_k \equiv 1$, $dh_k \equiv 0$ on $\bar x(\R)$ and $h_k \to 1$ in $C^\infty$. There exists a spanning disk $\mathcal D \subset (S^3,\xi)$ for $\bar x(\R)$ such that
\begin{itemize}
  \item The characteristic distribution of $\mathcal D$ has precisely one singularity. This singularity is a positive nicely elliptic point.
  \item There is a neighborhood $\mathcal O \subset \mathcal D$ of $\bar x(\R)$ such that $R_k|_p \not\in T_p\mathcal D$ for every $p\in \mathcal O \setminus \bar x(\R)$, when $k$ is large enough.
\end{itemize}
Here $R_k$ is the Reeb vector of $h_k\lambda = h_kf\lambda_0$. If we further assume that the $h_k\lambda$ are non-degenerate contact forms, then the disk $\DD$ may be arranged to satisfy the above properties, and also $y(\R) \not\subset \DD$ whenever $y$ is a periodic trajectory of $R_k$ satisfying $y(\R) \neq \bar x(\R)$, for $k$ large enough.
\end{lemma}

\begin{proof}
Combining Theorem~\ref{special_disk_0} with Lemma~\ref{special_disk} we get a disk $\DD$ satisfying the first two conditions. Suppose that the $h_k\lambda$ are non-degenerate contact forms and let $f:\D\to S^3$ be a smooth embedding satisfying $\DD = f(\D)$ and $f(0)=e$, and fix $\delta>0$ small. Then $$ X = \{ g \in C^\infty(\D,S^3) \mid g(z)=f(z) \text{ if } |z|\leq \delta \text{ or } 1-\delta\leq |z|\leq1 \} $$ is closed in the complete metric space $C^\infty(\D,S^3)$. Taking $\delta$ small and $k_0$ large we may assume all $R_k$ are transverse to $f(\{z\in\D:|z|\leq\delta \text{ or } 1-\delta\leq |z|< 1 \})$ if $k\geq k_0$. Here we used the properties of $\DD$ near the boundary given by Lemma~\ref{special_disk}. For fixed $k\geq k_0$ and periodic trajectory $y$ of $R_k$ satisfying $y(\R) \neq \bar x(\R)$ consider $X_{k,y} = \{ g\in X \mid y(\R) \subset g(\D) \}$. The $R_k$ have only countably many geometrically distinct periodic orbits since the $h_k\lambda$ are non-degenerate, and one easily checks that each $(X_{k,y})^c$ is open and dense in $X$ since $y(\R) \neq \bar x(\R)$. By Baire's category theorem $(\bigcup_{k,y} X_{k,y})^c$ is dense in $X$. Consequently, after a further $C^\infty$-small perturbation of $\DD$ supported away from $\{e\} \cup \partial\DD$, we may assume $\DD$ satisfies all the required properties.
\end{proof}

\subsection{Bishop families and fast planes}\label{bishop_families}

Consider the set
\begin{equation}\label{}
  \mathcal F = \{ h\in C^\infty(S^3,\R^+) \mid h\lambda \text{ is non-degenerate} \}.
\end{equation}
$\mathcal F$ is residual in $C^\infty(S^3,\R^+)$. Well-known arguments show that there exists a sequence $h_k: S^3 \to \R$ of smooth functions satisfying
\begin{equation}\label{functions_h_k}
  \begin{array}{cccc}
    h_k \in \mathcal F, & h_k \equiv1,\ dh_k \equiv 0 \text{ on } \bar x(\R) & \text{and} & h_k \to 1 \text{ in } C^\infty.
  \end{array}
\end{equation}
In the remainder of Section~\ref{deg_case} $R_k$ denotes the Reeb vector of $\lambda_k = h_k \lambda$. As noted before, $t\mapsto \bar x(t)$ is a $\bar T$-periodic orbit of $R_k$.

We now recall how arguments from~\cite{char1,char2} and~\cite{hry,pedro} prove the following existence result for fast planes asymptotic to $\bar P$. The sets $\J(\xi,d\lambda_k)$ and $\J(\xi,d\lambda)$ all coincide with $\J(\xi,d\lambda_0)$ and will be simply denoted by $\J$.

\begin{theorem}\label{existence_fast_non_deg_1}
For every $k$ large enough there exists some $J \in \J$ and an embedded fast $\jtil_k$-holomorphic finite-energy plane $\util_k : \C \to \R\times S^3$ asymptotic to $\bar P$.
\end{theorem}

The almost complex structure $\jtil_k$ is defined by~\eqref{alm_cpx_str_def} using $J$ and the Reeb vector field $R_k$. The remainder of this subsection is devoted to the proof of Theorem~\ref{existence_fast_non_deg_1}.

Let $\DD$ be the special spanning disk for $\bar x(\R)$ given by applying Lemma~\ref{special_disk_2} to the sequence $h_k$~\eqref{functions_h_k}, which has a unique singular point $e$. For each $k$ define
\begin{equation}\label{}
  \text{area}_{d\lambda_k}(\DD) = \sup \left\{ \int_U d\lambda_k : U \subset \DD \text{ is open} \right\}
\end{equation}
Since $\lambda_k \to \lambda$ we have $\text{area}_{d\lambda_k}(\DD) \to \text{area}_{d\lambda}(\DD)$. Applying Darboux's theorem we find a fixed small neighborhood $V$ of $e$ and embeddings $\Psi_k : V \to \R^3$ satisfying $\Psi_k(e) = (0,0,0)$, $\Psi_k^*(dz+xdy) = \lambda_k$. By arguments from~\cite{93} the disk $\DD$ may be perturbed to another spanning disk $\DD_k$ for $\bar x(\R)$ so that
\[
  \begin{array}{cccc}
    e \in \DD_k, & \DD_k \setminus V = \DD \setminus V & \text{and} & \DD_k \subset \{ z=-\frac{1}{2}xy \} \text{ near }e.
  \end{array}
\]
Moreover, this perturbation can be constructed so that $e$ is the only singularity of the characteristic foliation of $\DD_k$ and
\begin{equation}\label{energy_bound_C}
  \sup_k \text{area}_{d\lambda_k}(\DD_k) = C < \infty.
\end{equation}
Since $V$ can be chosen arbitrarily small, $\DD_k$ still satisfies all the properties obtained by Lemma~\ref{special_disk_2}.

For each $J\in \J$ and $k$ consider the set $\M_k(J)$ of solutions of the following boundary-value problem:
\begin{equation}\label{bishop_disks_defn}
  \left\{
  \begin{aligned}
    &\util= (a,u): \D \to \R\times S^3 \text{ satisfies } \bar\partial_{\jtil_k}(\util) = 0 \\
    &\util \text{ is an embedding, } a\equiv0 \text{ on } \partial\D \text{ and } u(\partial \D) \subset \DD_k\setminus \{e\} \\
    &u(\partial\D) \text{ winds once and positively around } e.
  \end{aligned}
  \right.
\end{equation}
$\D$ is equipped with its standard complex structure $i$ and its usual orientation, while $\DD$ is oriented by $\lambda_0|_{T\partial\DD}>0$. Here $\jtil_k$ is determined by $J$ and $R_k$ as in~\eqref{alm_cpx_str_def}.

\begin{lemma}[Hofer]\label{existence_bishop_disks}
For every $k$, the set of $J \in \J$ such that $\M_k(J) \neq\emptyset$ contains an non-empty open subset of $\J$.
\end{lemma}

\begin{proof}
There exists some $J$ that satisfies $J \cdot \partial_x = -x\partial_z + \partial_y$ near $e$ in the coordinates $(x,y,z)$ given by the Darboux chart $\Psi_k$. Then the maps $\util_\tau = (a_\tau,u_\tau) : \D \to \R\times S^3$ defined by
\[
  \begin{array}{ccc}
    a_\tau(s+it) = \frac{\tau^2}{4}(s^2+t^2-1), &  & \Psi_k \circ u_\tau(s+it) = \left( \tau s, \tau t, -\frac{\tau^2}{2}st \right)
  \end{array}
\]
form a 1-parameter family of elements of $\M_k(J)$, which converges to the constant map $(0,e)$ as $\tau \to 0^+$. It is shown in~\cite{93} that the linearization of the Cauchy-Riemann operator $\bar\partial_{\jtil_k}$ at a given element of $\M_k(J)$ is automatically surjective, with respect to a Fredholm theory of disks with boundary on $\DD_k$. Thus if $\M_k(J) \neq\emptyset$ for some $J$ then $\M_k(J') \neq\emptyset$ for every $J'$ in a neighborhood of $J$ in $\J$.
\end{proof}

Since $\lambda$ is dynamically convex we find $k_0$ such that if $k\geq k_0$ then all closed $\lambda_k$-Reeb orbits satisfying $\mu_{CZ} \leq 2$ have action strictly larger than the constant $C$ in~\eqref{energy_bound_C}. Consider for given $k\geq k_0$ and $J \in \J$ the set $\M^*_k(J)$ of non-constant finite-energy $\jtil_k$-holomorphic embeddings $$ \util = (a,u) : \D\setminus \Gamma \to \R\times S^3, $$ where $\Gamma \subset \D \setminus \partial\D$ is some non-empty finite set, satisfying
\begin{equation}\label{}
  \left\{
  \begin{aligned}
    & a\equiv0 \text{ on } \partial\D \text{ and } u(\partial\D) \subset \DD_k\setminus \{e\}, \\
    & u(\partial\D) \text{ winds once and positively around } e, \\
    & \text{every } z\in\Gamma \text{ is a negative puncture and } \\
    & \int_{\D\setminus \Gamma} u^*d\lambda_k > 0.
  \end{aligned}
  \right.
\end{equation}

\begin{theorem}[Hofer, Wysocki and Zehnder]\label{J_gen}
If $k$ is large enough then there exists a dense set $\J_{gen}(k)$ such that if $J \in \J_{gen}(k)$ then $\M^*_k(J) = \emptyset$.
\end{theorem}

The above theorem, which is proved in~\cite{props3}, roughly follows from the fact that the elements of $\M_k^*(J)$ are solutions of a (finite) number of Fredholm problems for which the expected dimension of the solution set is negative. Here we strongly use that, by Stokes theorem, any given $\util \in \M^*_k(J)$ with $k\geq k_0$ is asymptotic at a negative puncture $z\in \Gamma$ to a closed $\lambda_k$-Reeb orbit $P_z$ with action $\leq C$ and, consequently, satisfies $\mu_{CZ}(P_z) \geq 3$.

The fundamental result regarding the above defined Bishop families of disks is

\begin{theorem}[Hofer, Wysocki and Zehnder]\label{consequence_bishop_family}
For each $k$ large enough there exists $J \in \J$ and a sequence $\{\util_n = (a_n,u_n)\} \subset \M_k(J)$ satisfying $\util_n(\D) \cap (\R\times\bar x(\R)) = \emptyset \ \forall n$ and $\util_n \to F_{\bar P}$ in $C^\infty_{loc}(\D\setminus \{0\},\R\times S^3)$. Here $F_{\bar P} : \D\setminus \{0\} \to \R\times S^3$ is the map $F_{\bar P}(\est) = (\bar Ts, \bar x(\bar Tt))$.
\end{theorem}

\begin{proof}[Sketch of Proof]
We fix three leaves $l_1,l_i,l_{-1}$ of the characteristic foliation of $\DD_k$. They have finite length and connect $e$ to the boundary $\partial \DD_k$. By Lemma~\ref{existence_bishop_disks} and Lemma~\ref{J_gen} we find $J\in \J_{gen}(k)$ such that $\M_k(J) \neq \emptyset$. For each $\util=(a,u) \in \M_k(J)$ the curve $u(\partial \D) \subset \DD_k\setminus\{e\}$ hits every leaf of the characteristic foliation transversely and once. Thus we may define $\tau(\util)\in \R^+$ to be the length of the piece of $l_1$ connecting $e$ to $u(\partial \D)$ with respect to some metric. Let $\util_0 \in \M_k(J)$ be close to the point $(0,e)$, as described in Lemma~\ref{existence_bishop_disks}. The connected component $\mathcal Y \subset \M_k(J)$ containing $\util_0$ is a trivial principle $\text{M\"ob}(\D)$-bundle over the open real interval $I=\tau(\mathcal Y)$ with projection $\tau:\mathcal Y\to I$, where $\text{M\"ob}(\D)$ denotes the group of holomorphic self-diffeomorphisms of $\D$. This is proved using the implicit function theorem, see~\cite{93} for details. There is a unique global section $t\in I \mapsto \util^t = (a^t,u^t)$ satisfying $\tau(\util^t) = t$, $u(z) \in l_z$ for $z\in\{1,i,-1\}$. Let $t^* = \sup I$ and $\bar t$ be the length of $l_1$. We follow~\cite{char1,char2} closely. A non-trivial intersection argument shows that $\util(\D) \cap \R\times \bar x(\R) = \emptyset$ for every $\util \in \mathcal Y$, and a bubbling-off analysis proves that the family $\{\util^t\}$ has uniform $C^\infty$-bounds on a fixed neighborhood of $\partial \D$. Take $\{t_n\} \subset I$, $t_n \to t^*$. If the sequence $\{\util^{t_n}\}$ is $C^1$-bounded then, up to a subsequence, we may assume $\util^{t_n} \to \util$. It is possible to show that $\util \in \mathcal Y$ and, by the implicit function theorem, $\tau$ takes values larger than $t^*$ on $\mathcal Y$, a contradiction. Thus bubbling-off occurs and one finds a non-empty finite set $\Gamma \subset \D\setminus \partial \D$ and a non-constant finite-energy $\jtil_k$-holomorphic map $\vtil =(b,v) : \D\setminus \Gamma \to \R\times S^3$ such that, up to selection of a subsequence, $\util^{t_n} \to \vtil$ in $C^\infty_{loc}(\D\setminus \Gamma)$ as $n\to\infty$. Note that $E(\vtil) \leq \sup_n E(\util^{t_n}) \leq C$. Moreover, $v(\partial \D) \subset \DD_k \setminus\{e\}$ winds positively and once around $e$ and $\vtil$ is an embedding. Then $\int v^*d\lambda_k = 0$ since $J \in \J_{gen}(k)$ and $k$ was fixed large enough. It follows that $v(\D\setminus \Gamma) \subset \hat x(\R)$ for some periodic $\lambda_k$-Reeb orbit and, by the properties of $\DD_k$, we must have $\hat x (\R)= x(\R)$ and $\#\Gamma = 1$. Up to reparametrization by an element of $\text{M\"ob}(\D)$ we may assume $\Gamma = \{0\}$. Hence $t^* = \bar t$ and the argument is complete.
\end{proof}

From now on we follow closely the arguments from~\cite{hry,pedro}. We fix $k$ large and consider $J$ and the sequence $\util_n=(a_n,u_n) \in \M_k(J)$ given by Theorem~\ref{consequence_bishop_family}. Let $\sigma(C)>0$ be small enough so that
\begin{itemize}
  \item Every closed $\lambda_k$-Reeb orbit has action strictly larger than $\sigma(C)$.
  \item If $P'=(x',T')$ and $P''=(x'',T'')$ are closed $\lambda_k$ Reeb orbits satisfying $T',T''\leq C$ and $T'\neq T''$ then $|T'-T''| > \sigma(C)$.
\end{itemize}
Fix $\epsilon>0$ small enough so that $$ \left| \int_{B_\epsilon(0)} u_n^*d\lambda_k -\bar T \right| \leq \sigma(C)/2, \ \text{ for } n\gg 1. $$ Note that $\epsilon$ exists since $m_r = \lim_{n\to\infty} \int_{B_r(0)} u_n^*d\lambda_k$ exists for every $r>0$ and $m_r \to \bar T$ as $r\to 0$. Following~\cite{fols}, we take $z_n \in \D$ satisfying $a_n(z_n) = \inf_\D a_n$ and $\delta_n>0$ such that $u_n^*d\lambda_k$ integrates to $\sigma(C)$ over $B_\epsilon(0)\setminus B_{\delta_n}(0)$. Then $z_n \to 0$ and $\delta_n\to0$. Choose $R_n\to\infty$ such that $\delta_nR_n \to 0$ and define
\begin{equation}\label{}
  \vtil_n = (b_n,v_n) : B_{R_n}(0) \to \R\times S^3
\end{equation}
by $b_n(z) = a_n(z_n+\delta_nz) - a_n(z_n+2\delta_n)$ and $v_n(z) = u_n(z_n+\delta_nz)$. Up to the choice of a subsequence we find a finite set $\Gamma_1 \subset \C$ and a finite-energy $\jtil_k$-holomorphic map $\vtil =(b,v):\C\setminus \Gamma_1 \to \R\times S^3$ such that $$ \vtil_n \to \vtil \text{ in } C^\infty_{loc}(\C\setminus \Gamma_1). $$ Clearly $E(\vtil) \leq \sup_n E(\vtil_n) \leq \sup_n E(\util_n) \leq C$. We may assume $\Gamma_1$ consists of negative punctures and $E(\vtil)>0$ if $\Gamma_1 \neq \emptyset$. If $\Gamma_1 = \emptyset$ then $$ \bar T \sim \int_{B_\epsilon(0)} u_n^*d\lambda_k = \int_{B_\epsilon(0)\setminus B_{\delta_n}(0)} u_n^*d\lambda_k + \int_{B_{\delta_n}(0)} u_n^*d\lambda_k \to \sigma(C) + \int_\D v^*d\lambda_k $$ which proves $E(\vtil)>0$ in this case as well. By Theorem~\ref{partial_asymptotics}, $\vtil$ is asymptotic to some closed $\lambda_k$-Reeb orbit at the (unique) positive puncture $\infty$, and using results of finite-energy cylinders with low $d\lambda_k$-area one can prove that this asymptotic orbit is $\bar P$. Here it is crucial that $\lambda_k$ is non-degenerate.

\begin{lemma}\label{dlambda_non_zero}
$\int_{\C\setminus \Gamma_1} v^*d\lambda_k >0$.
\end{lemma}

\begin{proof}
If not then, by Theorem~\ref{zera_dlambda_theorem}, there is a non-constant polynomial $p:\C\to \C$ such that $\Gamma_1 = p^{-1}(0)$ and $\vtil = F_{P}\circ p$. Here $P = (x,T)$ is some simply covered $\lambda_k$-Reeb orbit and $F_P : \C\setminus \{0\} \to \R\times S^3$ is the map $z=\est \mapsto (Ts,x(Tt))$. It is easy to see that $0\in \Gamma_1$ since $b_n(0) = \inf_{B_{R_n}(0)} b_n$ and points of $\Gamma_1$ are negative punctures. Since $\bar P$ is simply covered we have $\bar P = P$ and $p(z) = Az$ for some $A \neq 0$. We can estimate as above $$ \begin{aligned} \frac{\sigma(C)}{2} &\geq \int_{B_\epsilon(0)} u_n^*d\lambda_k - \bar T = \int_{B_\epsilon(0)\setminus B_{\delta_n}(0)} u_n^*d\lambda_k + \int_{B_{\delta_n}(0)} u_n^*d\lambda_k -\bar T \\ &\to \sigma(C) + \int_{\partial\D} v^*\lambda_k - \bar T = \sigma(C) + \bar T - \bar T = \sigma(C) \end{aligned} $$ which is a contradiction.
\end{proof}

Let us fix a non-vanishing global section
\begin{equation}\label{global_section}
  Z : S^3 \to \xi.
\end{equation}

\begin{lemma}\label{non-vanishing_pi_du}
Consider the projection $\pi_k : TS^3 \to \xi$ along the Reeb direction $\R R_k$. Then $\pi_k \cdot du_n$ is nowhere vanishing over $\D$ when $n$ is large enough.
\end{lemma}

\begin{proof}
If $n$ is large then, in view of Theorem~\ref{consequence_bishop_family}, $u_n(\partial \D) \subset \mathcal O$, where $\mathcal O \subset \DD_k$ is a neighborhood of $\partial \DD_k = \bar x(\R)$ satisfying $p\in \mathcal O \setminus \bar x(\R) \Rightarrow R_k|_p \not\in T_p\DD_k$. In the following we introduce polar coordinates $(r,\theta) \in (0,1]\times \R/2\pi\Z \simeq \D\setminus \{0\}$ and denote by $\partial_\theta u_n$ and $\partial_r u_n$ the corresponding partial derivatives. If $z\in\partial \D$ and $\pi_k \cdot du_n(z) = 0$ then $0 = \lambda_k \cdot \partial_\theta u_n(z) = \partial_r a_n(z)$, contradicting the strong maximum principle ($\Delta a\geq 0$).

Let $V$ be a vector field on $\DD_k$ parametrizing the characteristic foliation of $\DD_k$ which points out of $\DD_k$ at the boundary. Then $V$ has a unique non-degenerate source at $e$. We claim that $\pi_k\cdot \partial_\theta u_n(z)$ and $V\circ u_n(z)$ are linearly independent in $\xi|_{u_n(z)}$ for every $z\in\partial \D$ and $n$ large enough. If not we find $c_1,c_2\in \R$, $z\in\partial \D$ and $n$ so that $u_n(\partial\D) \subset \mathcal O$ and $$ 0= c_1 \pi_k\cdot \partial_\theta u_n(z) + c_2 V(u_n(z)) = \pi_k \cdot (c_1\partial_\theta u_n(z) + c_2 V(u_n(z))). $$ Consequently $c_1\partial_\theta u_n(z) + c_2 V(u_n(z)) \in \R R_k|_{u_n(z)} \cap T_{u_n(z)}\DD_k = 0$ which implies $c_1\partial_\theta u_n(z) + c_2 V(u_n(z))=0$ since $u_n(z) \in \mathcal O\setminus \bar x(\R), \ \forall z\in \partial \D$. If $c_1=0$ then $c_2=0$ because $V$ does not vanish on $u_n(\partial \D)$. If not then $\lambda_k$ vanishes on $\partial_\theta u_n(z)$ contradicting the strong maximum principle. In particular we showed $$ \wind(\pi_k \cdot \partial_\theta u_n(e^{i2\pi t}), V \circ u_n(e^{i2\pi t})) = 0. $$

By Theorem~\ref{consequence_bishop_family} the embedded loops $t\mapsto u_n(e^{i2\pi t})$ are transverse to $\xi$ and bound a disk $F_n$ in $\DD_k$ containing $e$, when $n$ is large. Since $e$ is the only (positive) zero of the section $V$ of $\xi|_{F_n}$ standard degree theory gives $$ \wind(V \circ u_n(e^{i2\pi t}),Z \circ u_n(e^{i2\pi t}))=1. $$

If $x+iy$ are standard holomorphic coordinates on $\D$ then, since $\pi_k \cdot du_n$ does not vanish over $\partial\D$ when $n\gg1$, we get $$ \wind (\pi_k \cdot \partial_x u_n(e^{i2\pi t}),\pi_k \cdot \partial_\theta u_n(e^{i2\pi t})) = \wind (1,ie^{i2\pi t}) =-1. $$

All the above winding numbers were computed endowing $\xi$ with the orientation induced by $d\lambda_k$, see Remark~\ref{rmk_winding_numbers}. Finally we compute
\[
  \begin{aligned}
    \wind(\pi_k \cdot \partial_x u_n(e^{i2\pi t}),Z \circ u_n(e^{i2\pi t})) &= \wind(\pi_k \cdot \partial_x u_n(e^{i2\pi t}),\pi_k \cdot \partial_\theta u_n(e^{i2\pi t})) \\
    & +\wind(\pi_k \cdot \partial_\theta u_n(e^{i2\pi t})),V \circ u_n(e^{i2\pi t})) \\
    & +\wind(V \circ u_n(e^{i2\pi t}),Z \circ u_n(e^{i2\pi t})) \\
    & = -1 + 0 + 1 = 0.
  \end{aligned}
\]
Thus $\pi_k \cdot \partial_x u_n$ does not vanish when $n$ is large enough.
\end{proof}

We will now use lemmas~\ref{dlambda_non_zero} and~\ref{non-vanishing_pi_du} to show that $\vtil$ is an embedded fast finite-energy plane. In view of Theorem~\ref{partial_asymptotics} and Definition~\ref{non_deg_puncture} there exists $R_0>1$ such that $\pi_k \cdot dv(z) \neq 0$ if $|z|\geq R_0$. If $R\geq R_0$ is fixed then
\begin{equation}\label{wind_infty_v_calculated}
  \begin{aligned}
    & \wind(\pi_k \cdot \partial_r v(Re^{i2\pi t}), Z\circ v(Re^{i2\pi t})) \\
    & = \lim_{n\to\infty} \wind(\pi_k \cdot \partial_r v_n(Re^{i2\pi t}), Z\circ v_n(Re^{i2\pi t})) \\
    & = \lim_{n\to\infty} \wind \left( \pi_k \cdot \left.\frac{d}{d\rho}\right|_{\rho=1} u_n(z_n+ \rho \delta_nR e^{i2\pi t}), Z\circ u_n(z_n+R\delta_ne^{i2\pi t}) \right) \\
    & = \lim_{n\to\infty} \wind(\pi_k \cdot \partial_ru_n(e^{i2\pi t}),Z\circ u_n(e^{i2\pi t})) \\
    & = \lim_{n\to\infty} \wind(\pi_k \cdot \partial_ru_n(e^{i2\pi t}),\pi_k \cdot \partial_xu_n(e^{i2\pi t})) \\
    & + \lim_{n\to\infty} \wind(\pi_k \cdot \partial_xu_n(e^{i2\pi t}), Z\circ u_n(e^{i2\pi t})) \\
    & = 1+0 = 1.
  \end{aligned}
\end{equation}
In the third and fourth equalities we strongly used Lemma~\ref{non-vanishing_pi_du}. $Z\circ v$ is a non-vanishing section of $v^*\xi$ and~\eqref{wind_infty_v_calculated} implies $$ \wind_\infty(\vtil,\infty,Z\circ v) = +1. $$ See section~\ref{section_alg_inv} for the definitions. Note that $\vtil$ is asymptotic to $\lambda_k$-Reeb orbits $\{P_z\}_{z\in\Gamma}$ with action $\leq E(\vtil) \leq C$. Consequently, since we assumed $k$ is large, $\mu_{CZ}(P_z)\geq 3$ for every $z\in\Gamma$. By Theorem~\ref{precise_asymptotics}, $\wind_\infty(\vtil,z,Z\circ v) = \wind(e,Z)$ for some eigenvector $e$ of the asymptotic operator $A_{P_z}$ satisfying $A_{P_z}e = \nu e$ with $\nu>0$. The inequality $\mu_{CZ}(P_z)\geq 3$ implies $$ \wind_\infty(\vtil,z,Z\circ v) \geq 2. $$ It follows from Lemma~\ref{lemma_wind_pi_infty} that
\begin{equation}\label{}
  \begin{aligned}
    0 &\leq \wind_\pi(\vtil) = \wind_\infty(\vtil) - 2 + \#\Gamma_1 + 1 \\
    &\leq 1- 2\#\Gamma_1 - 2 + \#\Gamma_1 + 1 = -\#\Gamma_1
  \end{aligned}
\end{equation}
which proves $\Gamma_1 = \emptyset$ and $\vtil$ is a fast plane. It only remains to show that $\vtil$ is an embedding. It is an immersion since $\wind_\pi(\vtil) = 0$ implies $\pi_k \cdot dv$ has no zeros. Since $\vtil$ is somewhere injective ($\bar P$ is simply covered), self-intersections are isolated. But if $\vtil$ has self-intersection points then positivity and stability of self-intersections of pseudo-holomorphic immersions implies that the disks $\util_n$ have self-intersections. This would be a contradiction since the $\util_n$ are embeddings. The proof of Theorem~\ref{existence_fast_non_deg_1} is complete.

\subsection{Obtaining fast planes in the degenerate case}

We still consider an arbitrary sequence $h_k$ as in~\eqref{functions_h_k}, and the contact forms $\lambda_k = h_k\lambda = h_kf\lambda_0$. Then $\bar P=(\bar x,\bar T)$ is a prime closed orbit of the Reeb vector field $R_k$ associated to $\lambda_k$. We let $\sigma_k>0$ be a constant defined as in the previous subsection, satisfying:
\begin{itemize}
  \item Every closed $\lambda_k$-Reeb orbit has action larger than $\sigma_k$.
  \item If $P'=(x',T')$ and $P''=(x'',T'')$ are closed $\lambda_k$-Reeb orbits satisfying $T',T''\leq \bar T$ and $T'\neq T''$ then $|T'-T''| > \sigma_k$.
\end{itemize}

Let $J \in \J$, $k$ be a fixed large integer and $H \subset \R\times S^3$ be a compact set. Consider the set $\Lambda(k,J,H)$ of embedded fast finite-energy $\jtil_k$-holomorphic planes $\util = (a,u) : \C\to \R\times S^3$ asymptotic to $\bar P$ satisfying the normalization conditions
\begin{equation}\label{}
  \begin{array}{ccc}
    \util(0) \in H, &  & \int_\D u^*d\lambda_k = \bar T - \sigma_k.
  \end{array}
\end{equation}
Here $\jtil_k$ is the almost-complex structure defined as in~\eqref{alm_cpx_str_def} using the vector field $R_k$ and the given complex structure $J:\xi\to \xi$. The following statement is Theorem 2.2 from~\cite{hry} applied to our particular situation. The arguments are implicitly contained in Appendix~\ref{appA}.

\begin{theorem}\label{comp_thm_1}
Suppose every $\lambda_k$-Reeb orbit $\hat P = (\hat x,\hat T)$ with $\hat T\leq \bar T$ satisfies $\mu_{CZ}(\hat P) \geq 3$. If $H \cap \R\times \bar x(\R) = \emptyset$ then $\Lambda(k,J,H)$ is $C^\infty_{loc}$-compact.
\end{theorem}

It is important to note that the choice of $J\in \J$ in the statement above is arbitrary. The next statement is a direct application of Theorem 2.5 from~\cite{hry} to our case.

\begin{theorem}\label{openbook}
Let $J \in \J$ and $k$ be such that there exists an embedded fast finite-energy $\jtil_k$-holomorphic plane asymptotic to $\bar P$. Suppose that $\Lambda(k,J,H)$ is $C^\infty_{loc}$-compact for every compact subset $H \subset \R\times S^3$ with $H \cap \R\times \bar x(\R) = \emptyset$. Then for every $l\geq 1$ one finds a $C^l$-map
\[
 \util = (a,u) : \R/\Z \times \C \rightarrow \R \times S^3
\]
satisfying:
\begin{enumerate}
 \item Each $\util(\vartheta,\cdot)$ is an embedded fast finite-energy $\jtil_k$-holomorphic plane asymptotic to $\bar P$.
 \item $u(\vartheta,\C) \cap \bar x(\R) = \emptyset \ \forall \vartheta\in \R/\Z$ and the map $u : \R/\Z \times \C \rightarrow S^3 \setminus \bar x(\R)$ is an orientation preserving $C^l$-diffeomorphism.
 \item Each $\cl{u(\vartheta,\C)}$ is a global surface of section for the $\lambda_k$-Reeb flow.
\end{enumerate}
\end{theorem}

There is an important consequence.

\begin{lemma}\label{all_linked}
Every periodic $\lambda$-Reeb orbit geometrically different from $\bar P$ links non-trivially with $\bar P$.
\end{lemma}

\begin{proof}
It is important to note that, as our arguments show so far, the conclusions of Theorem~\ref{existence_fast_non_deg_1} and Theorem~\ref{comp_thm_1} are true for a sequence of contact forms $g_k\lambda$, where $\{g_k\}$ is an arbitrary sequence of functions as in~\eqref{functions_h_k}. Let us suppose, by contradiction, that there exists some $\lambda$-Reeb orbit $\hat P = (\hat x,\hat T)$ satisfying $\hat x(\R) \neq \bar x(\R)$ which is not linked to $\bar P$. We find a sequence $g_k : S^3\to\R$ of smooth functions satisfying $g_k \to 1$ in $C^\infty$, $g_k \equiv 1$ and $dg_k \equiv0$ on $\hat x(\R) \cup \bar x(\R)$, such that $g_k\lambda$ are non-degenerate contact forms. In particular, $\bar P$ and $\hat P$ are periodic orbits for the Reeb flow of $g_k\lambda$. If $k$ is large enough the conclusions of theorems~\ref{existence_fast_non_deg_1} and~\ref{comp_thm_1} are true for the contact forms $g_k\lambda$. Then, by Theorem~\ref{openbook} every periodic Reeb orbit for the Reeb flow of $g_k\lambda$ (geometrically distinct of $\bar P$) is linked with $\bar P$, but this contradicts the fact that $\hat P$ is not linked to $\bar P$.
\end{proof}

Our goal is now to produce fast planes in the degenerate case. The first step in our construction is the following lemma asserting that if fast planes exist then they exist abundantly. The proof is postponed to the appendix.

\begin{lemma}\label{lemma_anyJ}
Let $M$ be a closed 3-manifold equipped with a non-degenerate contact form $\alpha$ such that $c_1(\xi)|_{\pi_2(M)} \equiv 0$, where $\xi = \ker\alpha$ is the induced contact structure. Suppose $P = (x,T)$ is a prime closed $\alpha$-Reeb orbit and let $\J_{\text{fast}}(P)$ be the set of $d\alpha$-compatible complex structures $J:\xi \to \xi$ such that there exist embedded fast finite-energy $\jtil$-holomorphic planes asymptotic to~$P$. Here $\jtil$ is defined as in~\eqref{alm_cpx_str_def} using $J$ and the Reeb vector field of~$\alpha$. If every contractible $\alpha$-Reeb orbit $P'=(x',T')$ with $T'\leq T$ satisfies $\mu_{CZ}(P')\geq 3$ then either $\J_{\text{fast}}(P) = \J(\xi,d\alpha)$ or $\J_{\text{fast}}(P) = \emptyset$.
\end{lemma}

The second step is to prove

\begin{lemma}\label{lemma_fast_plane_lambda}
For every $p\in S^3\setminus \bar x(\R)$ and every $J \in \J(\xi,d\lambda)$ there exists an embedded fast finite-energy $\jtil$-holomorphic plane $\util = (a,u): \C \to \R\times S^3$ asymptotic to $\bar P$ and satisfying $p\in u(\C)$.
\end{lemma}

The remaining of this subsection is devoted to the proof of Lemma~\ref{lemma_fast_plane_lambda}. Let us fix $p\in S^3 \setminus \bar x(\R)$ and $J \in \J(\xi,d\lambda)$ arbitrarily. Since $\lambda$ is dynamically convex we find that if $k$ is large enough then the conclusions of Theorem~\ref{existence_fast_non_deg_1} and Theorem~\ref{comp_thm_1} are simultaneously satisfied for the contact form $\lambda_k$. In view of Lemma~\ref{lemma_anyJ} and Theorem~\ref{openbook}
we obtain, for $k$ large, an embedded fast $\jtil_k$-holomorphic plane $\util_k=(a_k,u_k)$ asymptotic to $\bar P$ satisfying $p\in u_k(\C)$. Here $\jtil_k$ is the $\R$-invariant almost complex structure on $\R\times S^3$ given by $\jtil_k \cdot \partial_a = R_k$, $\jtil_k|_\xi \equiv J$ where $R_k$ is the Reeb vector field associated to $\lambda_k$. Setting $\jtil$ by $\jtil \cdot \partial_a = R$, $\jtil|_\xi \equiv J$, then $\jtil_k\to\jtil$ in $C^\infty$. We would like to examine now the limiting behavior of the sequence $\util_k$.

Let $\gamma>0$ be a number so that $\gamma < T'$ for every closed Reeb orbit $P' = (x',T')$ of $\lambda$ and of $\lambda_k$, for every $k$. The number $\gamma$ exists since, otherwise, we would find $k_j\to\infty$ and $T_j$-periodic $\lambda_{k_j}$-Reeb orbits with $T_j\to0$. Then the Arzel\`a-Ascoli theorem would provide a subsequence of these orbits which converge to a rest point of the $\lambda$-Reeb flow, a contradiction. After reparametrizing we may assume
\begin{equation}\label{normalization_conditions}
  \begin{array}{ccc}
    \int_\D u_k^*d\lambda = \bar T-\gamma & \text{and} & a_k(0) = \inf_\C a_k.
  \end{array}
\end{equation}
As usual, consider
\begin{equation}\label{}
  \Gamma = \{ z\in \C \mid \exists z_j \to z \text{ and } k_j \text{ such that } |d\util_{k_j}(z_j)| \to \infty \}.
\end{equation}
By Lemma~\ref{beforeclaim} we may assume, up to the choice of a subsequence, that $\#\Gamma<\infty$ and \eqref{normalization_conditions} implies $\Gamma \subset \D$. Choosing a further subsequence still denoted by $\util_k$, there is no loss of generality, perhaps after discarding a few points of $\Gamma$ and translating in the $\R$-coordinate, to assume that $\forall z\in\Gamma$ $\exists z_k\to z$ such that $|d\util_k(z_k)| \to \infty$, and that we find a smooth $\jtil$-holomorphic map $\util = (a,u) :\C\setminus \Gamma \to \R\times S^3$ satisfying $\util_k \to \util$ in $C^\infty_{loc}(\C\setminus \Gamma)$ and $E(\util) \leq \bar T$. We claim that $\util$ is non-constant. This is true if $\Gamma\neq\emptyset$ since every $z\in\Gamma$ is a bubbling-off point and, consequently, must have positive mass~\eqref{mass_defn}. If $\Gamma = \emptyset$ then this follows from $\int_\D u^*d\lambda = \bar T-\gamma > 0$.

The following statement is Lemma 6.24 from~\cite{hry}. Its original proof is contained in~\cite{props3}.

\begin{lemma}\label{lemma_nice_intersection}
If $\alpha$ is a contact form on $S^3$ and $\wtil$ is an embedded fast finite-energy plane in $\R\times S^3$ asymptotic to a closed (necessarily prime) Reeb orbit $\tilde P=(\tilde x,\tilde T)$ satisfying $\mu_{CZ}(\tilde P) \geq 3$, then $\wtil(\C) \cap \R\times \tilde x(\R) = \emptyset$.
\end{lemma}

\begin{lemma}\label{lemma_asymp_limit_util}
$u(e^{2\pi(s+it)}) \to \bar x(\bar Tt)$ in $C^\infty(\mathbb R/\mathbb Z,S^3) \mod \ \mathbb R/\mathbb Z$ as $s\to+\infty$.
\end{lemma}

\begin{proof}
As explained in~\cite{93}, see~\ref{fe_curves_section}, it is always possible to define the mass at a puncture $z$ of $\util$ by
\[
  m = - \lim_{\rho \to 0^+} \int_{\partial B_\rho(z)} u^*\lambda.
\]
$m=0 \Leftrightarrow z \text{ is removable}$. Each $z\in\Gamma$ is non-removable and negative, from where we conclude $\infty$ is a positive puncture. By the fundamental results from~\cite{93}, for every sequence $R_n \to \infty$ one finds a subsequence $R_{n_j}$, a real constant $c$ and a periodic $\lambda$-Reeb orbit $\hat P=(\hat x,\hat T)$ such that $u(R_{n_j}e^{i2\pi t}) \to \hat x(\hat Tt)$ as $j\to\infty$. Arguing by contradiction, assume $\hat P \neq \bar P$. If $\hat P$ and $\bar P$ are geometrically distinct then, when $j$ is large, $u(R_{n_j}e^{i2\pi t})$ is linked to $\bar P$ in view of Lemma~\ref{all_linked}. Since for every $j$ one has $u_k(R_{n_j}e^{i2\pi t}) \to u(R_{n_j}e^{i2\pi t})$ then $u_k(R_{n_j}e^{i2\pi t})$ is also linked to $\bar P$ when $j$ and $k$ are large. Therefore $\util_k(\C) \cap \R\times \bar x(\R){\neq \emptyset}$, contradicting Lemma~\ref{lemma_nice_intersection}. If $\hat P$ and $\bar P$ are not geometrically distinct then $\hat P = (\bar x,N\bar T)$, for some $N\geq 2$. But this implies $E(\util) = N\bar T> \bar T$, again a contradiction.
\end{proof}

\begin{lemma}\label{lemma_props_util}
We have that $\Gamma = \emptyset$, $\pi \cdot du$ is nowhere vanishing and $$ \lim_{s\to +\infty} {\rm wind}(t\mapsto \pi \cdot \partial_su(s,t), t\mapsto Z\circ u(s,t)) = 1 $$ where $u(s,t) = u(\est)$ and $Z$ is the global non-vanishing section~\eqref{global_section}.
\end{lemma}

\begin{proof}
If $\Gamma = \emptyset$ then $\int_\D u^*d\lambda = \bar T-\gamma > 0$. Assume $\Gamma\neq\emptyset$ and $\pi \cdot du$ vanishes identically. In view of Lemma~\ref{lemma_asymp_limit_util} and of Theorem~\ref{zera_dlambda_theorem} we have $\#\Gamma=1$ and, moreover, there are complex constants $A\neq0,D$ such that $\util(z) = F_{\bar P}(Az+D)$, where $F_{\bar P}(\est) = (\bar Ts,\bar x(\bar Tt))$. If $|D/A| <1$ then $$ \bar T = \int_{\partial \D} u^*\lambda = \lim_k \int_\D u_k^*d\lambda = \bar T - \gamma, $$ an absurdity. Thus $\Gamma = \{-D/A\} \subset \partial \D$ and we find $z'$ close to $-D/A$ such that $a(z') < a(0)-1$, implying $a_k(z') < a_k(0)$ for $k$ large and contradicting~\eqref{normalization_conditions}. We showed $\Gamma \neq\emptyset \Rightarrow \int u^*d\lambda > 0$.

As observed in the proof of Lemma~\ref{lemma_asymp_limit_util}, any given $z\in\Gamma$ is a negative puncture and we can find $\rho_n \to 0^+$ and a periodic $\lambda$-Reeb orbit $P' = (x',T')$ such that $u(z+ \rho_n e^{i2\pi t}) \to x'(T't+c)$ for some $c\in \R$. Since $\int u^*d\lambda>0$, we have $T'<\bar T$, and this implies $P'$ and $\bar P$ are geometrically distinct since $\bar P$ is prime. Fixing $n$ large we conclude that $u(z+ \rho_n e^{i2\pi t})$ is linked to $\bar P$ by Lemma~\ref{all_linked}, and hence so is $u_k(z+ \rho_n e^{i2\pi t})$ when $k$ is large. Consequently $\util_k(\C) \cap \R\times \bar x(\R){\neq \emptyset}$, contradicting Lemma~\ref{lemma_nice_intersection}. So far we showed $\Gamma = \emptyset$ and $\int u^*d\lambda>0$.

By the similarity principle there exists $\rho_n\to+\infty$ such that $\pi \cdot du(\rho_n e^{i2\pi t}) \neq 0, \ \forall t$. Since $\wind_\pi(\util_k) = 0$ for each $k$, we compute for any fixed $n$:
\[
  \begin{aligned}
    \wind(t\mapsto &\pi\cdot\partial_s u(\rho_ne^{i2\pi t}), t\mapsto Z\circ u(\rho_ne^{i2\pi t})) \\
    &= \lim_k \wind(t\mapsto \pi\cdot\partial_s u_k(\rho_ne^{i2\pi t}), t\mapsto Z\circ u_k(\rho_ne^{i2\pi t})) \\
    &= \lim_k \left( \{\text{algebraic count of zeros of } \pi\cdot du_k \text{ on } B_{\rho_n} \} + 1\right) = 1
  \end{aligned}
\]
Taking the limit as $n\to +\infty$, we use standard degree theory to conclude that the algebraic count of zeros of $\pi\cdot du$ on $\C$ is $1-1=0$. The proof is complete since all zeros count positively.
\end{proof}

So far we have proved that $\util$ is a $\jtil$-holomorphic finite-energy plane asymptotic to $\bar P$ only in the sense of Lemma~\ref{lemma_asymp_limit_util}. If we can show $\infty$ is a non-degenerate puncture of $\util$ (in the sense of Definition~\ref{non_deg_puncture}) then it will follow from Lemma~\ref{lemma_props_util} that $\util$ is a fast plane.

The following statement is a combination of the non-trivial lemmas 8.1 and 8.3 from~\cite{convex} with a slightly modified analysis of Appendix B.

\begin{lemma}\label{lemmas_convex}
Given any Martinet tube $(U,\Psi)$ of $\bar P$, there exists $R_0>0$ such that $\bigcup_k u_k(\C\setminus B_{R_0}) \subset U$ and, moreover, if we write $$ (id_\R\times\Psi) \circ \util_k(\est) = (a_k(s,t),\theta_k(s,t),x_k(s,t),y_k(s,t)) \text{ for } e^{2\pi s}\geq R_0 $$ then one finds constants $c_k,d_k\in \R$ such that 
\[
\begin{aligned}  
&\lim_{s\to+\infty} \sup_{k,t} \left( |D^\beta x_k(s,t)| + |D^\beta y_k(s,t)| \right) = 0 \\ 
&\lim_{s\to+\infty} \sup_{k,t} \left( |D^\beta [a_k(s,t) - \bar Ts - d_k]| + |D^\beta [\theta_k(s,t) - t - c_k]| \right) = 0  
\end{aligned}
\]
for every partial derivative $D^\beta = \partial^{\beta_1}_s\partial^{\beta_2}_t$. 
\end{lemma}

We choose a Martinet tube $(U,\Psi)$ for $\lambda$ and $\bar P$, as described in Definition~\ref{defn_martinet}, so that the section $t\mapsto (\Psi^*\partial_x)\circ \bar x(\bar Tt)$ of $\xi_{\bar P}$ extends as a non-vanishing section for every disk spanning $\bar P$. Here $B \subset \C$ is an open ball centered at the origin and $\Psi : U \to \R/\Z \times B$ is a diffeomorphism satisfying the following properties.
\begin{enumerate}
  \item $\Psi(\bar x(\bar T\theta)) = (\theta,0,0), \ \forall \theta \in \R/\Z$,
  \item $\Psi^*(g(d\theta+xdy)) = \lambda$ for some smooth $g:\R/\Z\times B\to\R^+$ satisfying $g\equiv\bar T$ and $dg \equiv 0$ on $\R/\Z\times 0$.
\end{enumerate}

It follows from the particular form of the $\lambda_k$ that
\[
 \begin{array}{ccc}
   \Psi_*\lambda_k = g^k(d\theta+xdy), & g^k \to g \text{ in } C^\infty_{loc}, & g^k\equiv\bar T \text{ and } dg^k \equiv0 \text{ on } \R/\Z\times \{0\}.
 \end{array}
\]
We shall write $\partial_\theta,\partial_x,\partial_y$ for the vectors $\Psi^*\partial_\theta, \Psi^*\partial_x, \Psi^*\partial_y$. In the frame $\{\partial_x,-x\partial_\theta+\partial_y\}$ the almost complex structure $J$ is represented by a matrix
\[
  j = \begin{pmatrix} j_{11} & j_{12} \\ j_{21} & j_{22} \end{pmatrix}
\]
The Reeb vector $R_k$ associated to $\lambda_k$ is represented in these coordinates by
\[
 R_k = (g^k)^{-2}(g^k+xg^k_x, g^k_y-xg^k_\theta, -g^k_x) = (R^1_k, R^2_k, R^3_k).
\]
An analogous formula, with $g$ in the place of $g^k$, holds for the Reeb vector $R=(R^1,R^2,R^3)$ associated to $\lambda$. Set
\[
 \begin{array}{cc}
   Y_k = (R^2_k,R^3_k), & Y = (R^2,R^3)
 \end{array}
\]
and
\[
 \begin{array}{cc}
   D^k(\theta,z) = \int_0^1 Y_k'(\theta,\tau z)d\tau, & D(\theta,z) = \int_0^1 Y'(\theta,\tau z)d\tau
 \end{array}
\]
where the prime denotes a derivative with respect to the variable $z=(x,y)$. By Lemma~\ref{lemmas_convex} we can find $s_0 \gg1$ so that the functions
\begin{equation}\label{}
  \begin{aligned}
    & z_k(s,t) = (x_k(s,t),y_k(s,t)),\ \ z(s,t) = (x(s,t),y(s,t)) \\
    & D^k(s,t) = D^k(\theta_k(s,t),z_k(s,t)), \ \ D(s,t) = D(\theta(s,t),z(s,t)) \\
    & j^k(s,t) = j(\theta_k(s,t),z_k(s,t)), \ \ j(s,t) =j(\theta(s,t),z(s,t))
  \end{aligned}
\end{equation}
are all well-defined. Here $a,a_k,\theta,\theta_k,x,x_k,y$ and $y_k$ are defined as in the statement of Lemma~\ref{lemmas_convex}. Then Cauchy-Riemann equations can be written as
\begin{equation}\label{cr_real_part}
  \begin{array}{c}
    \partial_sa_k - g^k(\theta_k,z_k)(\partial_t\theta_k + x_k\partial_ty_k) = 0 \\
    \partial_ta_k + g^k(\theta_k,z_k)(\partial_s\theta_k + x_k\partial_sy_k) = 0
  \end{array}
\end{equation}
and
\begin{equation}\label{cr_z_k}
  \partial_sz_k + j^k(s,t)\partial_tz_k + S_k(s,t)z_k = 0
\end{equation}
where $S_k(s,t) = [\partial_ta_kI - \partial_sa_kj^k(s,t)]D^k(s,t)$. The map $S(s,t)$ is defined analogously. Similar equations hold for $a(s,t)$, $\theta(s,t)$ and $z(s,t)$.

We can assume $c_k\to c$ for some $c\in\mathbb R$. Let us consider
\begin{equation}\label{def_Sk}
 \begin{aligned} S^\infty_k(t) & := -\bar Tj(t+c_k,0,0)Y'_k(t+c_k,0,0) \\ S^\infty(t) & := -\bar Tj(t+c,0,0)Y'(t+c,0,0) \end{aligned}
\end{equation}
where the $c_k$ are given by Lemma~3.18. This lemma and the properties of the functions $g_k$ and $g$ readily imply that
\begin{equation}\label{nice_Sk}
  \begin{aligned}
  & \lim_{s\to+\infty} \sup_{t,k} |D^\gamma [S_k(s,t)-S^\infty_k(t)]| = 0  \ \forall \gamma, \\
  & \lim_{s\to+\infty} \sup_t |D^\gamma [S(s,t)-S^\infty(t)]| = 0  \ \forall \gamma. \\
  & \lim_{k\to +\infty} \sup_t |D^m[S^\infty_k(t) - S^\infty(t)]| = 0 \ \forall m.
 \end{aligned}
\end{equation}

Note that $-j(t+c_k,0,0)\partial_t - S^\infty_k(t)$ is the representation of the asymptotic operator at $\bar P$ of the contact form $\lambda_k$ with respect to the $d\lambda_k$-symplectic frame $\{ e^k_1=\partial_x/\sqrt{g^k}, e^k_2 = (-x\partial_\theta + \partial_y)/\sqrt{g^k} \}$. Analogously, $-j(t+c,0,0)\partial_t - S^\infty(t)$ represents the asymptotic operator at $\bar P$ of the form $\lambda$ with respect to the $d\lambda$-symplectic frame $\{ e_1=\partial_x/\sqrt{g}, e_2 = (-x\partial_\theta + \partial_y)/\sqrt{g} \}$.

Take a smooth function $M : \R/\Z\times B \to \sp(1)$ satisfying
\begin{equation}\label{Mk}
  M(\theta,x,y)j(\theta,x,y) = J_0M(\theta,x,y) \text{ for all } (\theta,x,y)\in \R/\Z\times B.
\end{equation}
$M$ exists since $J$ is compatible with $d\lambda$ and $d\lambda_k$ $\forall k$. Setting
\begin{equation}\label{zeta_k}
\begin{array}{c}
  \begin{array}{cc}
    A_k(s,t) = M(\theta_k(s,t),x_k(s,t),y_k(s,t)), & \zeta_k(s,t) = A_k(s,t) z_k(s,t),
  \end{array} \\
  \Lambda_k(s,t) = (A_kS_k - \partial_sA_k - J_0\partial_tA_k)A_k^{-1}
\end{array}
\end{equation}
where $S_k = S_k(s,t)$ and $A_k = A_k(s,t)$ then, in view of~\eqref{cr_z_k}
\begin{equation}\label{equation_normalized_k}
  \partial_s\zeta_k + J_0 \partial_t \zeta_k + \Lambda_k \zeta_k = 0.
\end{equation}
Defining functions $A$, $\Lambda$ and $\zeta$ of $(s,t)$ analogously using $M$, $S$ and $z$ we obtain the corresponding equation
\begin{equation}\label{equation_normalized}
  \partial_s\zeta + J_0 \partial_t \zeta + \Lambda \zeta = 0.
\end{equation}
Defining
\begin{equation}\label{}
  \begin{array}{c}
    C_k^\infty(t) = (M(t+c_k,0,0)S_k^\infty(t) - J_0\partial_tM)M(t+c_k,0,0)^{-1} \\
    C(t) = (M(t+c,0,0)S^\infty(t) - J_0\partial_tM)M(t+c,0,0)^{-1}
  \end{array}
\end{equation}
then Lemma~\ref{lemmas_convex} and~\eqref{nice_Sk} together imply
\begin{equation}\label{nice_Lambda_k}
  \begin{aligned}
  & \lim_{s\to+\infty} \sup_{t,k} |D^\gamma [\Lambda_k(s,t)-C^\infty_k(t)]| = 0  \ \forall \gamma, \\
  & \lim_{s\to+\infty} \sup_t |D^\gamma [\Lambda(s,t)-C^\infty(t)]| = 0  \ \forall \gamma. \\
  & \lim_{k\to +\infty} \sup_t |D^m[C^\infty_k(t) - C^\infty(t)]| = 0 \ \forall m.
 \end{aligned}
\end{equation}
One checks that $-J_0\partial_t-C^\infty_k(t)$ and $-J_0\partial_t - C^\infty(t)$ are the representations of the corresponding asymptotic operators at $\bar P$ with respect to other symplectic frames.

\begin{remark}[$C^{l,\alpha,\delta}_0$-topology]\label{topology}
If $l\in\Z^+$, $\alpha\in(0,1)$, $\delta<0$ and $s_0 \in\R$ consider the space $C^{l,\alpha,\delta}_0([s_0,\infty)\times \R/\Z,\R^m)$ defined in~\cite{props3} as the set of maps $f(s,t)$ satisfying
\begin{itemize}
  \item $e^{-\delta s}D^\beta f(s,t) \in C^{0,\alpha}$ on $\left[s_0,\infty\right)\times \R/\Z$, $\forall |\beta|\leq l$,
  \item $\lim_{R\to +\infty} \Vert e^{-\delta s} D^\beta f(s,t)\Vert_{C^{0,\alpha}(\left[R,+\infty\right)\times
     S^1)}=0 \ \forall \ |\beta| \leq l$.
\end{itemize}
Then $C^{l,\alpha,\delta}_0([s_0,\infty)\times \R/\Z,\R^m)$ is a Banach space with the norm
\begin{equation}\label{norma_alpha_delta}
  \Vert f \Vert_{l,\alpha,\delta} = \Vert e^{-\delta s}f(s,t) \Vert_{C^{l,\alpha}(\left[s_0,\infty\right)\times\R/\Z)}.
\end{equation}
\end{remark}

The following proposition should be seen as some kind of uniform asymptotic analysis. The proof is postponed to the appendix, the arguments are essentially found in~\cite{props1}.

\begin{propn}\label{propcomp}
Suppose the maps $K_k : [0,+\infty)\times \R/\Z \rightarrow \R^{2n\times 2n} \ (k\geq 1)$ and $K^\infty_k,K^\infty : \R/\Z \rightarrow \R^{2n\times 2n}$ are smooth and satisfy:
\begin{enumerate}
  \item $K^\infty_k(t),K^\infty(t)$ are symmetric $\forall t$, and $K^\infty_k \to K^\infty$ in $C^\infty$ as $k\to +\infty$.
  \item $\lim_{s\rightarrow+\infty} \sup_{t,k} \left|D^\gamma[K_k(s,t)-K_k^\infty(t)]\right| = 0\ \forall \gamma$.
\end{enumerate}
Consider the unbounded self-adjoint operator $L$ on $L^2(\R/\Z,\R^{2n})$ defined by
\[
  Le = -J_0\dot{e} - K^\infty e.
\]
Denote $E = C^{l,\alpha,\delta}_0 \left( \left[ 0,+\infty \right) \times \R/\Z,\R^{2n} \right)$ for some $\delta < 0$ and $l\geq1$. Suppose $\{X_k\} \subset E$ are smooth maps satisfying
\begin{equation}\label{eqns_X_k}
  \partial_s X_k + J_0\partial_t X_k + K_k X_k = 0  \ \forall k.
\end{equation}
If $\delta\not\in\sigma(L)$ and $X_k$ is $C^\infty_{loc}$-bounded then $\{X_k\}$ has a convergent subsequence in $E$.
\end{propn}

The section $\kappa = Z \circ \bar x(\bar Tt)$ is a non-vanishing section of the bundle $(\bar x_{\bar T})^*\xi$ which extends to any disk spanning $\bar P$. Since the contact forms $\lambda_k$ are non-degenerate and $C^\infty$-convergent to $\lambda$, and the Conley-Zehnder index of $\bar P$ viewed as $\lambda_k$-Reeb orbit is $\geq 3$ if $k$ is large enough, there exists $\delta<0$ such that the eigenvalues with winding $1$ with respect to $\kappa$ of the asymptotic operators at $\bar P$ corresponding to $\lambda_k$ are less than $2\delta$. This follows from the description of the Conley-Zehnder index via self-adjoint operators described in~\ref{cz_operators_description}, and from the continuity of the spectrum of the asymptotic operators with respect to small perturbations of the contact form. Also, $\delta$ can be chosen so that it does not lie on the spectrum of the asymptotic operator at $\bar P$ corresponding to $\lambda$. By the identities $\wind_\infty(\util_k)=1$ we know that, as explained in Remark~\ref{asymp_evalue_wind_infty}, the asymptotic eigenvalue of $\util_k$ is $\leq 2\delta$ for every~$k$.

Using the number $\delta$ and an arbitrary choice of $l\geq 1$ to define $E$ as in the statement of Proposition~\ref{propcomp}, we conclude from Theorem~\ref{precise_asymptotics} that all the $z_k(s,t)$, and consequently also the $\zeta_k(s,t)$, belong to $E$. Note that $\zeta_k \to \zeta$ in $C^\infty_{\rm loc}$. In view of~\eqref{nice_Lambda_k} we can apply the above proposition to conclude that, up to selection of a subsequence, $\zeta_k$ converges in $E$. This implies that $\zeta \in E$ and $\zeta_k \to \zeta$ in $E$. In particular, since $l\geq 1$ was arbitrary, one finds $0<r<-\delta/2$ such that
\begin{equation}\label{upper2}
  \limsup_{s\to+\infty} e^{rs}|D^\gamma z(s,t)| = 0, \ \forall \gamma.
\end{equation}

Let us denote $$ F_k(s,t) = \begin{pmatrix} a_k(s,t) -\bar Ts - d_k \\ \theta_k(s,t) -t -c_k \end{pmatrix}. $$ Then, using~\eqref{upper2} and~\eqref{cr_real_part}, we have an equation $$ \partial_sF_k + \begin{pmatrix} 0 & -\bar T \\ \bar T^{-1} & 0 \end{pmatrix} \partial_tF_k + \Delta_k = 0 $$ for some sequence of functions $\Delta_k(s,t)$ satisfying $$ \lim_{s\to+\infty} \sup_{k,t} e^{rs}|D^\gamma \Delta_k(s,t)| = 0. $$

One can now argue as in~\cite{props1} to get $r>0$ (independent of $k\gg 1$) satisfying $$ \lim_{s\to+\infty} \sup_{k,t} e^{rs}|D^\gamma F_k(s,t)| = 0, \ \forall \gamma. $$ In particular, taking the limit as $k\to\infty$, we conclude that
\begin{equation}\label{upper3}
\begin{array}{ccc} e^{rs}\left( |\nabla [a(s,t)-Ts]| + |\nabla[\theta(s,t)-t]| \right) \to 0 & \text{as} & s\to+\infty. \end{array}
\end{equation}

Equations~\eqref{upper2} and~\eqref{upper3}, together with Lemma~\ref{lemma_asymp_limit_util} and Lemma~\ref{lemma_props_util}, imply that $\util$ has a non-degenerate puncture at $\infty$, in the sense of Definition~\ref{non_deg_puncture}. This shows $\util$ is a fast plane.

It remains to check that $\util$ is embedded. By Lemma~\ref{lemma_props_util}, $\util$ is an immersion. Self-intersections of $\util$ must be isolated since, otherwise, one could argue as in~\cite{props2} using Carleman's similarity principle to conclude that $\util$ is a non-trivial cover of some (somewhere injective) finite-energy plane. This would imply that $\bar P$ is not simply covered. Thus, self-intersections of $\util$ must be isolated. Stability and positivity of intersections would provide self-intersections of the $\util_k$ for $k$ large enough, a contradiction. Consequently, $\util$ has no self-intersections. If $w_k \in \C$ satisfies $u_k(w_k) = p$ then $\sup_k |w_k|<\infty$ by Lemma~\ref{lemmas_convex}, so we may assume $\{w_k\}$ is convergent in $\C$. It follows that $p \in u(\C)$, and Lemma~\ref{lemma_fast_plane_lambda} is proved.

\subsection{Constructing the global sections}

The construction of an open book decomposition with disk-like pages and binding $\bar P$ follows arguments from~\cite{convex} which are standard by now. Applying Theorem 2.3 from~\cite{hry} to our situation we get

\begin{lemma}\label{lemmafredholm}
Let $J \in \J(\xi,d\lambda)$ and $\util= (a,u)$ be an embedded $\jtil$-holomorphic fast finite-energy plane asymptotic to $\bar P$. Then $u(\C) \cap \bar x(\R) = \emptyset$ and $u:\C \to S^3\setminus \bar x(\R)$ is a proper embedding. Moreover for any $l$ there exists a $C^l$-embedding $f: \C\times B_r(0) \to \R\times S^3$, where $B_r(0) \subset \R^2$ denotes the ball of radius $r$ centered at the origin, satisfying:
\begin{enumerate}
  \item Each $f(\cdot,\tau)$, $\tau \in B_r(0)$, is an embedded fast finite-energy $\jtil$-holomorphic plane asymptotic to $\bar P$ and $f(\cdot,0) = \util$.
  \item If $\tau_0 \in B_r(0)$ and $\util_k$ are embedded $\jtil$-holomorphic fast finite-energy planes asymptotic to $\bar P$ satisfying $\util_k \to f(\cdot,\tau_0)$ in $C^\infty_{loc}$, then there exist $\tau_k \to \tau_0$ and complex numbers $A_k \to 1$, $B_k \to 0$ such that $\util_k(A_kz+B_k) = f(z,\tau_k)$, $\forall z\in \C$.
\end{enumerate}
\end{lemma}

Note that we need no generic assumptions on $\lambda$, like being non-degenerate or Morse-Bott. We refer to section 6 from~\cite{hry} for details.

We also need an uniqueness statement for fast planes which is lemmas 6.23 and 6.25 from~\cite{hry} combined. The proof is exactly the same of that of Theorem 4.11 from~\cite{props2}.

\begin{lemma}\label{equal_planes}
Let $\util=(a,u)$ and $\vtil=(b,v)$ be embedded fast finite-energy planes asymptotic to $\bar P$. Then either $u(\C) = v(\C)$ or $u(\C) \cap v(\C) = \emptyset$. In the first case one finds complex constants $B$ and $A\neq0$ and a real constant $c$ such that $\util(Az+B) = (b(z)+c,v(z)) \ \forall z\in \C$.
\end{lemma}

Let $p_0 \in S^3\setminus \bar x(\R)$ be a recurrent point for the Reeb flow $\phi_t$. Using Lemma~\ref{lemma_fast_plane_lambda} we find an embedded fast finite-energy plane $\util^t = (a^t,u^t)$ asymptotic to $\bar P$ satisfying $u^t(0) = \phi_t(p_0)$ and $a^t(0)=0$, for every $t\in \R$. The Reeb vector field $R$ is transverse to $u^0(\C)$ since $\wind_\pi(\util^0) = 0$. Thus $$ C = \{ t\in \R \mid \phi_t(p_0) \in u^0(\C) \} $$ is discrete, non-empty and closed. Our assumption on $p_0$ implies $C\setminus \{0\} \neq \emptyset$ and we will denote $T_0 = \inf C\cap \R^+ > 0$. Then $u^{T_0}(\C) = u^0(\C)$ by Lemma~\ref{equal_planes}.

\begin{lemma}\label{lemma_A}
If $0\leq t' < t'' < T_0$ then $u^{t'}(\C) \cap u^{t''}(\C) = \emptyset$.
\end{lemma}

\begin{proof}
We argue by contradiction and assume $u^{t'}(\C) \cap u^{t''}(\C) \neq \emptyset$. By Lemma~\ref{equal_planes} we must have $\Pi := u^{t'}(\C) = u^{t''}(\C)$. Thus $t'>0$ by the definition of $T_0$ since $t''<T_0$. We can concatenate the path $t\in [t',t''] \mapsto \phi_t(p_0)$ to a path contained in $\Pi$ to obtain a closed loop $\Gamma$ in $S^3\setminus \bar x(\R)$. $\Gamma$ intersects each plane $\{u^t(\C)\mid t'<t<t'' \}$ transversally and positively, thus $\Gamma$ is linked with $\bar P$. Consequently $\Gamma$ must also intersect $u^0(\C)$, again a contradiction to the definition of $T_0$.
\end{proof}

As before we consider a constant $$ 0< \gamma<\inf \{ T \mid P=(x,T) \text{ is a closed Reeb orbit of } \lambda \}. $$

\begin{lemma}\label{tail}
If $\vtil_k=(b_k,v_k)$ and $\vtil=(b,v)$ are fast embedded finite-energy planes asymptotic to $\bar P$ satisfying $\vtil_k \to \vtil$ in $C^\infty_{loc}$, and $\int_\D v_k^*d\lambda = \bar T - \gamma$ then for every neighborhood $U$ of $\bar x(\R)$ there exists $R_0>0$ such that $|z|\geq R_0 \Rightarrow v_k(z) \in U$.
\end{lemma}

\begin{proof}
This proof is essentially the proof of Lemma 8.1 from~\cite{convex}. Let $U$ be any given neighborhood of $\bar x(\R)$. Assuming the lemma is false we find $z_k \in \C$ satisfying $|z_k|\to+\infty$ and $v_k(z_k) \in \partial U$. Defining $(s_k,t_k)$ by $z_k = e^{2\pi (s_k+it_k)}$, we may consider a sequence $\wtil_k = (d_k,w_k) : \R\times \R/\Z \to \R\times S^3$ of finite-energy $\jtil$-holomorphic cylinders defined by
\[
\begin{array}{cc}
  d_k(s,t) = b_k(e^{2\pi((s+s_k)+i(t+t_k))}) - b_k(z_k), & w_k(s,t) = v_k(e^{2\pi((s+s_k)+i(t+t_k))}).
\end{array}
\]
For any fixed $L>0$ our assumptions on the $\vtil_k$ imply that $$ \int_{[-L,L]\times \R/\Z} w_k^*d\lambda \leq \gamma $$ when $k$ is large. Thus $|d\wtil_k|$ is $C^0_{loc}$-bounded by Lemma~\ref{beforeclaim}. Since $\wtil_k(0,0) \in 0\times S^3$ $\forall k$ then, up to the choice of a subsequence, we may assume $\wtil_k \to \wtil = (d,w)$ in $C^\infty_{loc}$ where $\wtil:\R\times \R/\Z \to \R\times S^3$ is a finite-energy $\jtil$-holomorphic map. Moreover, $\int w^*d\lambda \leq \gamma$ and we may estimate
\[
  \int_{\R/\Z} \bar x_{\bar T}^*\lambda - \int_{\R/\Z} w_k(s,\cdot)^*\lambda \leq \int_{\C\setminus \D} v_k^*d\lambda \leq \gamma, \ \forall s\in \R.
\]
This implies
\[
  \bar T \geq \int_{\R/\Z} w(s,\cdot)^*\lambda \geq \bar T - \gamma, \ \forall s\in \R
\]
Thus $\wtil$ is not a constant map, and if we identify $\R\times \R/\Z \simeq \C P^1 \setminus \{0,\infty\}$ via $z\simeq \est$ then $0$ is a non-removable negative puncture. Let us choose an arbitrary sequence $s_n \to -\infty$. By results from~\cite{93} $s_n$ has a subsequence, still denoted $s_n$, for which the following holds: there is a periodic Reeb orbit $P^- = (x^-,T^-)$ and a constant $c\in \R$ such that $w(s_n,t) \to x^-(T^-t + c)$ in $C^\infty(\R/\Z,S^3)$ as $n\to \infty$.

We claim $P^-$ is geometrically distinct of $\bar P$. If $\int w^*d\lambda=0$ then, by Theorem~\ref{zera_dlambda_theorem}, we know that $\wtil(\R\times \R/\Z) \subset \R\times x^-(\R)$. Consequently $x^-(\R)$ contains a point in $\partial U$, proving our claim in this case. If $\int w^*d\lambda >0$ then, by Stokes Theorem, $T^- < \bar T$, again implying our claim since $\bar T$ is the minimal period of $\bar x$.

Consequently, $P^-$ is linked to $\bar P$ and, fixing $n$ large enough, so is the loop $t\mapsto w(s_n,t)$. Now fixing $k$ large we conclude that $t\mapsto v_k(R_n e^{i2\pi t})$ is linked to $\bar P$ where $R_n = e^{2\pi s_n}$. This implies intersections of $\vtil_k$ with $\R\times \bar x(\R)$, contradicting Lemma~\ref{lemmafredholm}.
\end{proof}

\begin{lemma}\label{local_embeddings}
For every $t_0 \in \R$ and $l\geq 1$ there is a $C^l$-map $F = (H,G):\C\times I_\epsilon(t_0) \to \R\times S^3$ such that each $F(\cdot,t)$ is an embedded fast finite-energy plane asymptotic to $\bar P$ and, moreover, that $G: \C \times I_\epsilon(t_0) \to S^3\setminus \bar x(\R)$ is an embedding satisfying $G(t,0) = \phi_t(p_0)$ and $G(\C,t) = u^t(\C)$ $\forall t\in I_\epsilon(t_0) = (t_0-\epsilon,t_0+\epsilon)$.
\end{lemma}

Note that, by Lemma~\ref{equal_planes}, $G(\C,t) = u^t(\C)$ follows from $G(t,0) = \phi_t(p_0)$. Note also that the maps $G(\cdot,t)$ and $u^t(\cdot)$ might be different.

\begin{proof}
Consider the $C^l$-embedding $f:\C\times B_r(0) \to \R\times S^3$ given by Lemma~\ref{lemmafredholm} satisfying $f(\cdot,0) = \util^{t_0}(\cdot)$. Since the curve $t \mapsto (0,\phi_t(p_0))$ is transverse to $\util^{t_0}(\C)$ at $t = t_0$, we find a small interval $I$ around $t_0$ and unique $C^l$-smooth curves $\gamma : I \to B_r(0)$ and $\zeta : I \to \C$ satisfying $\gamma(t_0) = 0$, $\gamma' \neq 0$ on $I$, $\zeta(t_0) = 0$ and $f(\zeta(t),\gamma(t)) = (0,\phi_t(p_0))$. Let $t \in I \mapsto A(t) \in \C\setminus\{0\}$ be a $C^l$-smooth function such that the planes $z\mapsto f(\zeta(t)+A(t)z,\gamma(t))$ have $d\lambda$-area on $\D$ equal to $\bar T-\gamma$. Then $$ F(z,t) = (H(z,t),G(z,t)) = f(\zeta(t)+A(t)z,\gamma(t)) $$ is a $C^l$ map defined on $\C\times I$. Note that $DG(0,t_0)$ is non-singular since $R$ is transverse to $u^{t_0}(\C)$. Then, possibly after making $I$ smaller, we can find $\rho>0$ such that $G|_{B_\rho(0) \times I}$ is an embedding of $B_\rho(0) \times I$ into $S^3\setminus \bar x(\R)$.

Now we show that $G$ is an immersion. If not we find $(z_1,t_1) \in \C\times I$ and $(\delta z,\delta t) \neq 0$ such that $DG(z_1,t_1)\cdot (\delta z,\delta t) = 0$. Then $\delta t \neq0$ since $G(\cdot,t_1):\C\to S^3\setminus\bar x(\R)$ is an embedding by Lemma~\ref{lemmafredholm}, and we can assume $\delta t=1$. Plugging the formula for $G$ and differentiating we get
\begin{equation}\label{diff_tau}
  Df|_{(\zeta(t_1)+A(t_1)z_1,\gamma(t_1))} \cdot (A(t_1)\delta z + \zeta'(t_1) + A'(t_1)z_1,\gamma'(t_1)) = (c,0)
\end{equation}
where $(c,0) \in T_{\util^{t_1}(z_1)}(\R\times S^3)$, $c\neq 0$. Denoting $f=(h,g)$, we define $$ \Phi(s,z) = (h(\zeta(t_1)+A(t_1)z,\gamma(t_1))+cs,g(\zeta(t_1)+A(t_1)z,\gamma(t_1))). $$ Since $f$ is an embedding, we find unique $C^l$-smooth maps $w(s),\alpha(s) \in \C$, $\Gamma(s) \in B_r(0)$ satisfying $w(0) = \zeta(t_1)$, $\alpha(0) = A(t_1)$, $\Gamma(0) = \gamma(t_1)$, $\Gamma'(0)\neq 0$ and $$ \Phi(s,z) = f(w(s)+\alpha(s)z,\Gamma(s)). $$ Fix $z=z_1$ and differentiate in $s$ at $s=0$
\begin{equation}\label{diff_s}
  Df|_{(w(0)+\alpha(0)z_1,\Gamma(0))} \cdot (w'(0)+\alpha'(0)z_1,\Gamma'(0)) = (c,0).
\end{equation}
Subtracting~\eqref{diff_tau} from~\eqref{diff_s} we get $\Gamma'(0) = \gamma'(t_1)$. Here we strongly used that $Df$ is everywhere non-singular. Denoting $Q = (w(0),\Gamma(0)) = (\zeta(t_1),\gamma(t_1))$ we compute
\[
  \begin{aligned}
    (c,-R(\phi_{t_0}(p_0))) & = \left.\frac{d}{ds}\right|_{s=0} f(w(s),\Gamma(s)) - \left.\frac{d}{dt}\right|_{t=t_1} f(\zeta(t),\gamma(t)) \\
    & = Df|_{Q} \cdot (w'(0)-\zeta'(t_1),\Gamma'(0)-\gamma'(t_1)) \\
    & = Df|_{Q} \cdot (w'(0)-\zeta'(t_1),0) \in T_{\util^{t_1}(0)}\util^{t_1}(\C)
  \end{aligned}
\]
which implies $\wind_\pi(f(\cdot,\gamma(t_1)) > 0$, a contradiction. Thus $G$ is a $C^l$-immersion.

Note that $u^t(\C) = G(\C,t)$ for every $t\in I$ by Lemma~\ref{equal_planes}. We now claim that if $I$ is small enough then $G$ is injective on $\C\times I$. Suppose not. Then, again by Lemma~\ref{equal_planes}, there exist sequences $t^n_1,t^n_2 \to t_0$ such that $G(\C,t^n_1) = G(\C,t^n_2) \ \forall n$. We find unique $z_n \in \C$ satisfying $G(z_n,t^n_1) = G(0,t^n_2) = \phi_{t^n_2}(p_0)$. We must have $\liminf |z_n| \geq \rho>0$ because $G$ is 1-1 on $B_{\rho}(0) \times I$. Also, $\{\phi_t(p_0) \mid t\in \overline I\}$ does not meet a sufficiently small neighborhood $U$ of $\bar x(\R)$. Then Lemma~\ref{tail} implies $\limsup |z_n| < +\infty$. Consequently we can assume $z_n \to z^* \not= 0$. This implies $G(z^*,t_0) = G(0,t_0)$, a contradiction since $G(\cdot,t_0)$ is 1-1 by Lemma~\ref{lemmafredholm}. As remarked before, Lemma~\ref{equal_planes} implies $G(\C,t) = u^t(\C)$ as claimed.
\end{proof}

\begin{lemma}\label{lemma_B}
$\bigcup_{0\leq t < T_0} u^t(\C) = S^3\setminus \bar x(\R)$.
\end{lemma}

\begin{proof}
We will prove that $$ \Omega = \bigcup_{0<t<T_0} u^t(\C) \subset S^3\setminus \overline{u^0(\C)} $$ is open and closed in the connected set $S^3\setminus \overline{u^0(\C)}$. Let $q_n \in \Omega$, $q_n \to q\in S^3\setminus \overline{u^0(\C)}$. We use Lemma~\ref{local_embeddings} to find $0\leq \tau_1< \dots < \tau_N \leq T_0$, $\epsilon>0$ and smooth maps $F^i=(H^i,G^i) : \C\times I_\epsilon(\tau_i) \to \R\times S^3$ such that $\cup_{i=1}^N I_\epsilon(\tau_i) \supset [0,T_0]$ and $G^i(\C,t) = u^t(\C)$ $\forall t\in I_\epsilon(\tau_i)$. Moreover, $G^i:\C\times I_\epsilon(\tau_i) \to S^3\setminus \bar x(\R)$ is an embedding for each $i$. We find $t_n \in (0,T_0)$ such that $q_n \in u^{t_n}(\C)$. Up to the choice of a subsequence, $t_n \to t^* \in [0,T_0]$. If $t^* \in I_\epsilon(\tau_i)$ then $\exists z_n \in \C$ such that $q_n = G^i(z_n,t_n)$, and $\sup_n|z_n|<\infty$ by Lemma~\ref{tail}. So, up to the choice of a further subsequence, we may assume $z_n \to z^* \in \C$. Consequently $q = \lim_n G^i(z_n,t_n) = G^i(z^*,t^*)  \in G^i(\C,t^*) = u^{t^*}(\C)$. Hence $t^* \in (0,T_0)$ and $q\in \Omega$. Thus $\Omega$ is closed in $S^3\setminus \overline{u^0(\C)}$. By Lemma~\ref{local_embeddings} $\Omega$ is open in $S^3\setminus \overline{u^0(\C)}$.
\end{proof}

Lemmas~\ref{lemma_A},~\ref{local_embeddings} and~\ref{lemma_B} imply that the embedded planes $\{u^t(\C)\mid t\in[0,T_0]\}$ are pages of an open book decomposition of $S^3$ with binding $\bar P$. We still need to show that each page is a global surface of section for the Reeb flow.

\begin{lemma}\label{not_omega_limit}
Let $y$ be a Reeb trajectory in $S^3\setminus \bar x(\R)$ and assume that its $\omega$-limit set does not intersect $\bar x(\R)$. Then $y([a,+\infty)) \cap u^s(\D) \neq\emptyset$ for every $a\in \R$ and $s\in[0,T_0)$. Analogously, if we assume the $\alpha$-limit set of $y$ does not meet $\bar x(\R)$ then $y((-\infty,a]) \cap u^s(\D) \neq\emptyset$ for every $a\in \R$ and $s\in[0,T_0)$.
\end{lemma}

\begin{proof}
We only prove the lemma assuming the $\omega$-limit set of $y$ does not meet $\bar x(\R)$, the other case is analogous. Fix any $a\in \R$ and assume $y([a,+\infty)) \cap u^{t_0}(\C) \neq\emptyset$ for some $t_0 \in [0,T_0)$. Let
\[
  D = \{ t\in (t_0,T_0) \mid y([a,+\infty)) \cap u^r(\C) \neq\emptyset \ \forall r\in (t_0,t) \}.
\]
First we show $D \neq \emptyset$. Consider the embedding $G:\C\times I_\epsilon(t_0) \to S^3\setminus \bar x(\R)$ given by Lemma~\ref{local_embeddings} satisfying $G(\C,t) = u^t(\C) \ \forall t\in I_\epsilon(t_0)$, and take $\tau_0\geq a$ such that $y(\tau_0)\in u^{t_0}(\C)$. One can find a smooth curve $\gamma(\tau)$ such that $\gamma(\tau_0) = t_0$ and $y(\tau) \in G(\C,\gamma(\tau))$. Since $G(\C,t_0)$ is transverse to the Reeb vector we have $\gamma'(\tau_0) \neq 0$ and it is easy to check that $\gamma'(\tau_0)>0$.  Thus $\exists \delta>0$ such that $(t_0,t_0+\delta) \subset D$.

Now consider $t^* = \sup D$ and suppose, by contradiction, that $t^*<T_0$. By Lemma \ref{local_embeddings}, consider an embedding $G^*:\C\times I_\epsilon(t^*) \to S^3\setminus \bar x(\R)$ satisfying $G^*(\C,t) = u^t(\C) \ \forall t\in I_\epsilon(t^*)$. Take $\tau_n\geq a$ and $t_n <t^*$ such that $t_n\to t^*$ and $y(\tau_n) \in u^{t_n}(\D)$. Thus $\exists$ unique $z_n\in \C$ satisfying $G^*(z_n,t_n) = y(\tau_n)$. Using our hypothesis on the $\omega$-limit set of $y$ we find $q\in S^3\setminus \bar x(\R)$ such that $y(\tau_n) \to q$, after selecting a subsequence. Thus $\sup_n |z_n| < \infty$ by Lemma~\ref{tail}, so we may also assume $\lim_n z_n = z^*$ exists. This gives $G^*(z^*,t^*) = q$. If $\tau_n$ contains a bounded subsequence then $q \in y([a,+\infty))$ which proves $y([a,+\infty)) \cap u^{t^*}(\C) \neq \emptyset$. If $\tau_n \to +\infty$ then, by transversality of the Reeb vector with $u^{t^*}(\C)$, we again find that $y([a,+\infty)) \cap u^{t^*}(\C) \neq \emptyset$. Arguing as in the proof of $D\neq\emptyset$ we find $\delta>0$ such that $(t^*,t^*+\delta) \subset D$, contradicting the definition of $t^*$. Thus $t^* = T_0$.

So far we showed $T_0 =\sup D$ and $y([a,+\infty) \cap u^0(\D) = y([a,+\infty) \cap u^{T_0}(\D) \neq \emptyset$. Repeating the above arguments with $0$ in the place of $t_0$ we conclude that $y([a,+\infty) \cap u^s(\D) \neq \emptyset$ for every $s\in [0,T_0)$ as desired.
\end{proof}

Finally we need to recall the following lemma extracted from~\cite{hry} and proved by arguments from~\cite{convex}.

\begin{lemma}\label{omegalimit}
Let $\util = (a,u)$ be a fast finite-energy plane asymptotic to $\bar P = (\bar x,\bar T)$. If the $\omega$-limit ($\alpha$-limit) set of some Reeb trajectory $y(t)$ in $S^3\setminus \bar x(\R)$ intersects $\bar x(\R)$ then $\forall a\in\R$ $\exists t>a$ ($t<a$) such that $y(t) \in u(\C)$.
\end{lemma}

Lemmas \ref{not_omega_limit} and \ref{omegalimit} imply that each page $u^s(\C)$ is a global surface of section, for any $s\in [0,T_0)$. The proof of Theorem~\ref{main} is complete.

\section{Computing the self-linking number at a fixed point}\label{section_sl}

Suppose $\lambda$ is a (tight) dynamically convex contact form on $S^3$, and let $P = (x,T)$ be the boundary of the global disk section $D_0$ obtained by Theorem~\ref{theo_convex}, where $T$ is the minimal period of $x$. Let $P_1=(x_1,T_1)$ be an orbit given by a fixed point of the first return map to $D_0$. In particular, $P_1$ is unknotted and $T_1$ is the minimal period of $x_1$. In this section we prove

\begin{propn}\label{propn_sl}
$sl(P_1)=-1$.
\end{propn}

Orienting $S^3$ by $\lambda\wedge d\lambda$ one computes the linking number\footnote{The linking number $\link(c_1,c_2)$ of two oriented embedded circles $c_1$ and $c_2$ inside $S^3$ is defined as the oriented intersection number of $c_1$ with an oriented disk $D$ satisfying $\partial D = c_2$ (orientations included).} $\link(P_1,P)=+1$ by noting that $P_1$ intersects the disk-like global section once and positively. The Reeb vector field will be denoted by $R$.

\subsection{Outline of proof of Proposition~\ref{propn_sl}}

Below we construct a $C^0$ embedding $\varphi :\D\to S^3$ satisfying $\varphi(e^{i2\pi t}) = x_1(T_1t)$ that restricts to a smooth embedding $\varphi : \D\setminus\{0\} \to S^3$. Moreover, if $(r,\vartheta) \simeq re^{i\vartheta}$ are polar coordinates on $\D^* = \D\setminus\{0\}$ then the radial derivative $\varphi_r$ is never tangent to the Reeb vector. If $\varphi$ as just described exists then
\begin{equation}\label{section_W}
  W := \pi \cdot \varphi_r
\end{equation}
is a non-vanishing section of $\varphi^*\xi|_{\D^*}$, where $\pi$ is the projection~\eqref{projection}. Let $$ Z:S^3\to\xi $$ be a smooth global non-vanishing section and set
\begin{equation}\label{winding_global}
  m = \wind(W(e^{i2\pi t}),Z\circ \varphi(e^{i2\pi t}))
\end{equation}
where $\xi$ is oriented by $d\lambda$.

\begin{lemma}
If $\varphi$ as described above exists then $sl(P_1)=-m$.
\end{lemma}

\begin{proof}
Let $F_0 = \varphi(\D)$. For every neighborhood $\mathcal V$ of $\varphi(0)$ in $S^3$ there exists a smooth embedded disk $F \hookrightarrow S^3$ such that $F=F_0$ on $S^3\setminus \mathcal V$. Then $W$ defines a smooth section of $\xi|_{F\setminus \mathcal V}$ that can be extended smoothly to $F$. This extension will still be denoted by $W$. We claim that $W|_{\partial F}$ pushes $\partial F$ to a loop with zero intersection number with $F$. This can be seen by considering the homotopy $W^s = (1-s)W + s\varphi_r$ of vector fields along $\partial F$. The crucial fact is that $W^s$ is never tangent to $\partial F$, for every $s\in[0,1]$, which follows easily from $T\partial F = \R R|_{\partial F}$. Then $sl(P_1) = 0 - \wind_{\partial F}(W,Z|_{\partial F}) = -\wind_{\partial F_0}(W,Z|_{\partial F_0}) = -m$.
\end{proof}

The remaining of this section is devoted to the construction of $\varphi$ and to the proof of $m=1$.

\subsection{Constructing $\varphi$ and computing $m$}\label{construction_varphi}

We denote by $\xi|_P$ the bundle $(x_T)^*\xi \to \R/\Z$, where $x_T:\R/\Z\to M$ is the map $t\mapsto x(Tt)$. Consider coordinates $(\theta,x,y) \in \R/\Z \times \R^2$ and the contact form $\alpha_0 = d\theta + xdy$. According to~\cite{props1}, one can find a small open ball $B \subset \R^2$ centered at the origin, an open neighborhood $U$ of $x(\R)$ in $S^3$ and a Martinet tube $\Psi:U \stackrel{\sim}{\to} \R/\Z \times B$ as in Definition~\ref{defn_martinet}. We have $\Psi_*\lambda = f\alpha_0$ and $\Psi(x(Tt)) = (t,0,0)$, where $f$ satisfies $f|_{\R/\Z \times 0} \equiv T$ and $df|_{\R/\Z \times 0} \equiv 0$. It follows that $\Psi_*\xi = 0\times \R^2$ on $x(\R)$. According to Remark~\ref{rmk_martinet} we can and will assume that $\Psi^*\partial_x = Z$ over points of $P$.

Now consider the space $X_{l,\alpha,\delta}(P)$ of $C^{l,\alpha}_{loc}$-maps $F = (a,f) : \C\to \R\times S^3$ satisfying
\begin{itemize}
  \item $\exists d\in \R$ such that $a(s,t)-Ts-d \in C^{l,\alpha,\delta}_0([0,\infty)\times \R/\Z)$, where $a(s,t)=a(\est)$,
  \item $\exists R>0$ such that $f(\C\setminus B_R(0)) \subset U$,
  \item If $s_0\gg1$ and we define $(\theta(s,t),x(s,t),y(s,t)) = \Psi \circ f(\est)$ then $\exists c\in \R$ such that $(\theta-t-c,x,y) \in C^{l,\alpha,\delta}_0([s_0,\infty)\times \R/\Z,\R^3)$.
\end{itemize}
One can see that $X_{l,\alpha,\delta}(P)$ does not depend on the choice of the Martinet tube $(U,\Psi)$, and carries the structure of a separable smooth Banach manifold, but we will not make use of this fact. Here we only note that its topology is described by saying that $F_n=(a_n,f_n) \to F=(a,f)$ if, and only if,
\begin{itemize}
  \item $F_n \to F$ in $C^{l,\alpha}_{loc}(\C,\R\times S^3)$,
  \item If $s\gg 1$ then $f_n(\est),f(\est) \in U$ $\forall n$,
  \item Defining the functions $(a_n,\theta_n,x_n,y_n) = (id_\R\times \Psi) \circ F_n(\est)$ and also $(a,\theta,x,y) = (id_\R\times \Psi) \circ F(\est)$ for $s_0 \gg 1$ then the sequences $\theta_n-\theta$, $x_n-x$, $y_n-y$ and $a_n-a$ converge to zero in $C^{l,\alpha,\delta}_0([s_0,\infty)\times \R/\Z)$.
\end{itemize}

We may slightly perturb $\lambda$ in the $C^\infty$-topology to a non-degenerate contact form so that $P$ and $P_1$ are still Reeb orbits. As shown in the previous section, the conclusions of Theorem~\ref{existence_fast_non_deg_1} and Theorem~\ref{comp_thm_1} hold for this small perturbation, so that Theorem~\ref{openbook} can be applied. The conclusions of Theorem~\ref{openbook} are obtained in~\cite{hry} by the construction of a continuous map 
\begin{equation}\label{map_util}
  \util = (a,u) : \R/\Z \to X_{l,\alpha,\delta}(P)
\end{equation}
such that $(\vartheta,z) \mapsto u(\vartheta,z) := u(\vartheta)(z)$ is a diffeomorphism $\R/\Z \times \C \simeq S^3\setminus x(\R)$. Each page $\overline {u(\vartheta,\C)}$ is a global surface of section for the Reeb flow. The details of the proof of Theorem~\ref{openbook} are technical and follow closely the ideas from~\cite{char1,char2,convex}.

By the implicit function theorem we can find a smooth map $\vartheta \in \R/\Z \mapsto \gamma(\vartheta) \in \C$ such that $x_1(\R) \cap u(\vartheta,\C) = u(\vartheta,\gamma(\vartheta))$. Replacing $u$ by $u(\vartheta,z-\gamma(\vartheta))$ we still have a continuous map $\vartheta \mapsto u(\vartheta,\cdot) \in X_{l,\alpha,\delta}(P)$. By reparametrizing the $\R/\Z$-coordinate we can also assume, without loss of generality, that $u(\vartheta,0) = x_1(T_1\vartheta)$. We take the opportunity to state a useful lemma that follows from the above discussion.

\begin{lemma}\label{lemma_est}
Let $\util=(a,u)$ be the map~\eqref{map_util}. Composing with a smooth family of rotations we can assume $u(\vartheta,\est) \to x(Tt)$ as $s\to\infty$, $\forall \vartheta$. Moreover, $\forall \epsilon>0 \ \exists R_0 = e^{2\pi s_0} \gg 1$ such that $u(\vartheta,\C\setminus B_{R_0}(0)) \subset U \ \forall \vartheta$ and if we define $(\theta,x,y)(\vartheta,s,t) = \Psi \circ u(\vartheta,\est)$ for $s\geq s_0$ and $\vartheta,t\in\R/\Z$, then $$ e^{-\delta s}|D^\beta[(\theta,x,y)(\vartheta,s,t)-(t,0,0)]| \leq \epsilon $$ uniformly in $|\beta|\leq l$, $(\vartheta,t) \in \R/\Z\times \R/\Z$ and $s\geq s_0$.
\end{lemma}

For simplicity we assume $T_1 = 1$. Fix $0<\theta_1<1$, $0<\epsilon_0 <\min\{\theta_1/3,1/2\}$ and $s_0\gg1$ so that $s\geq s_0$ implies $$ w(\vartheta,s,t) = (\theta,x,y) = \Psi \circ u(\vartheta,\est) $$ is defined and
\begin{equation}\label{}
  \max \left\{ |w(\vartheta,s,t)-(t,0,0)|,\ |w_t(\vartheta,s,t)-(1,0,0)|,\ |(2s)^{3/2}w_s(\vartheta,s,t)| \right\} \leq \epsilon_0.
\end{equation}
This is possible in view of Lemma~\ref{lemma_est}.

Choose $0<r_0\ll1/2$ so small that $1/2r_0^2 \gg s_0$ and $2\theta_1/3r_0>11$. This implies $s := 1/2r^2 \gg s_0$ if $r\leq r_0$. Consider a smooth function $t:(0,1]\to\R$ satisfying $t^{\prime}(r)\geq0$, $t(r)\equiv0\text{ for } r_0\leq r\leq 1$ and $t(r)\equiv\frac{4\theta_1}{3r_0}r-\theta_1\text { for } 0 < r \leq\frac{r_0}{2}$.

Let $\rho : (0,1] \to [0,\infty)$ be a smooth function satisfying $\rho'<0$, $\rho(1)=0$ and $\rho(r) = \exp(\pi r^{-2})$ for $r\in(0,r_0]$. Now define for $\vartheta \in \R/\Z$ and $r \in (0,1]$
\begin{equation}\label{funcs}
  \begin{array}{ccc}
  z(r) = \rho(r) e^{i2\pi t(r)} &  & \varphi(r,\vartheta) = u(\vartheta,z(r)).
\end{array}
\end{equation}
Note that $$ \Psi \circ \varphi(r,\vartheta) = w(\vartheta,1/2r^2,t(r)) $$ is well-defined for $r\leq r_0$.

Note that $\varphi$ is a proper embedding of $(0,1]\times \R/\Z$ into $S^3 \setminus P$, and $\varphi(1,\vartheta) = u(\vartheta,0) = x_1(\vartheta)$. If we regard $\varphi$ as a map $\D^* \to S^3 \setminus P$ by $\varphi(re^{i2\pi\vartheta}) \simeq \varphi(r,\vartheta)$ then it can be extended to a continuous embedding of $\mathbb{D}$ into $S^3$ by setting $\varphi(0) = x(-\theta_1)$. The radial derivative $\varphi_r(\vartheta,r)$ over $\D^*$ is a non-vanishing vector of $Tu(\vartheta,\C)$, so it is never a multiple of $R$ since $u(\vartheta,\C)$ is an embedding transversal to the Reeb flow. This concludes the construction of the $C^0$-embedding $\varphi$ with the desired properties.

Regarding $(\theta,x,y) \in \R/\Z\times B$ as coordinate functions, then
\begin{equation}\label{claim_1}
  \theta(\Psi\circ \varphi(r,\vartheta)) \geq -\frac{2\theta_1}{3} \ \forall\vartheta \in \R/\Z, r\in[r_0/2,r_0].
\end{equation}
In fact, $\frac{r_0}{2}\leq r\leq r_0\Rightarrow s=(1/2\pi)\log e^{\pi/r^2}=1/2r^2\gg s_0$. Hence $\frac{r_0}{2}\leq r\leq r_0$ implies $|\Psi\circ \varphi(r,\vartheta)-(t(r),0,0)| \leq\epsilon_0<\theta_1/3$ and $$ t(r) \geq \frac{4\theta_1}{3r_0} \frac{r_0}{2} - \theta_1 = -\frac{\theta_1}{3}. $$ These two estimates together and the triangle inequality prove~\eqref{claim_1}.

Choose $\theta_2\in(-\theta_1,-2\theta_1/3)$ a regular value for the function $(\theta\circ\Psi \circ \varphi)|_{(0,r_0]\times \R/\Z}$ and set $$ S_0 = (\theta\circ \Psi\circ \varphi)^{-1}(\theta_2). $$ By~\eqref{claim_1} we have $(r,\vartheta) \in S_0 \Rightarrow r<r_0/2$. Also $S_0$ is bounded away from $\{0\}\times \R/\Z$ since $\lim_{r\rightarrow0} (\theta\circ \Psi\circ \varphi)(r,\vartheta) = -\theta_1<\theta_2$. Hence $S_0$ is a non-empty finite collection of embedded circles inside $(0,r_0/2)\times \R/\Z$. Now we claim that
\begin{equation}\label{claim_2}
  d\theta \cdot (\Psi\circ \varphi)_r > 10 \text{ on } S_0.
\end{equation}
To prove this note that $r<r_0/2$ implies
\[
(\Psi\circ \varphi)_r(r,\vartheta) = -(2s)^{3/2}w_s+\frac{4\theta_1}{3r_0}w_t
\]
Now we estimate
\[
  \begin{aligned}
    d\theta \cdot (\Psi\circ \varphi)_r &\geq \frac{4\theta_1}{3r_0}d\theta \cdot w_t - (2s)^{3/2}|w_s|   \\
    &\geq \frac{4\theta_1}{3r_0} \left[ d\theta \cdot (1,0,0) - |d\theta \cdot (w_t-(1,0,0))| \right] - (2s)^{3/2}|w_s| \\
    &\geq \frac{4\theta_1}{3r_0}(1-\epsilon_0)-\epsilon_0 \geq \frac{4\theta_1}{3r_0}\frac{1}{2}-\frac{1}{2} \geq 11-\frac{1}{2}>10
  \end{aligned}
\]
and establish~\eqref{claim_2}.

We also claim that $S_0$ has precisely one component which is a generator of $\pi_1((0,1]\times \R/\Z)$. In fact~\eqref{claim_2} implies that every component of $S_0$ is the image of a graph where $r$ is a function of $\vartheta$. Thus no component of $S_0$ is homotopically trivial in $(0,1]\times \R/\Z$ and every component of $S_0$ is a generator of $\pi_1((0,1]\times \R/\Z)$. Moreover, if there are two or more components then $d\theta \cdot (\Psi\circ \varphi)_r < 0$ along one of them, again in contradiction to~\eqref{claim_2}. From now on we might write $r(\vartheta) \in (0,r_0/2)$ to denote the unique smooth function such that $S_0 = \{ (r(\vartheta),\vartheta) : \vartheta \in \R/\Z\}$.

Consider the flat surface $\Sigma = \Psi^{-1}(\{ \theta = \theta_2 \})$ inside $U$. Then $$  \{u(\vartheta,z(r(\vartheta)))\} = \{ \varphi(r(\vartheta),\vartheta) \} = \varphi(S_0) \subset \Sigma \cap \text{image}(\varphi). $$ Note that $u(\vartheta,\C) \pitchfork \Sigma$ at the points of $\varphi(S_0)$ since there we have $|w_t-(1,0,0)|\leq\epsilon_0<1/2$ which implies $|d\theta \cdot w_t|\geq 1/2$. Denoting by $u(\vartheta)$ the map $z\mapsto u(\vartheta,z)$ we find, for each $\vartheta \in \R/\Z$, some $\beta(\vartheta) \in \R/\Z$ satisfying
\begin{equation}\label{nice_relation}
  d\theta \cdot d(\Psi\circ u(\vartheta))|_{z(r(\vartheta))} \cdot e^{i2\pi\beta(\vartheta)} = 0.
\end{equation}
By the above mentioned transversality and the implicit function theorem we can choose $\beta$ smoothly. Set
\[
  \sigma(\vartheta) = du(\vartheta)|_{z(r(\vartheta))} \cdot e^{i2\pi\beta(\vartheta)}
\]
and note that both $\sigma(\vartheta)$ and $\varphi_r(r(\vartheta),\vartheta)$ are vectors tangent to $u(\vartheta,\C)$ at the point $u(\vartheta)(z(r(\vartheta))) = \varphi(r(\vartheta),\vartheta)$, so $\pi \cdot \sigma(\vartheta)$ and $\pi \cdot \varphi_r(r(\vartheta),\vartheta)$ are non-vanishing vectors on $\xi$. Here we used that $u(\vartheta,\C)$ is transverse to the Reeb vector and $\pi$ is the projection~\eqref{projection}. It follows from~\eqref{claim_2}-\eqref{nice_relation} and from the definition of $\sigma$ that $\pi \cdot \sigma(\vartheta)$ and $\pi \cdot \varphi_r(r(\vartheta),\vartheta)$ are linearly independent $\forall \vartheta$ (it was used that $\pi \cdot du(\vartheta) : (\C,i) \to (\xi,J)$ is a complex isomorphism). We compute
\begin{equation}\label{wind1}
  \begin{aligned}
    m &= \wind(\pi \cdot \varphi_r(1,\vartheta), Z\circ \varphi(1,\vartheta)) = \wind(\pi \cdot \varphi_r(r(\vartheta),\vartheta),Z\circ \varphi(r(\vartheta),\vartheta)) \\
    &= \wind(\pi\cdot \varphi_r(r(\vartheta),\vartheta),\pi \cdot \sigma(\vartheta)) + \wind(\pi \cdot \sigma(\vartheta),Z\circ \varphi(r(\vartheta),\vartheta)) \\
    &= \wind(\pi \cdot \sigma(\vartheta),Z\circ \varphi(r(\vartheta),\vartheta)) \\
    &= \wind(\pi \cdot \sigma(\vartheta), \pi \cdot (\Psi^*\partial_x)\circ \varphi(r(\vartheta),\vartheta)).
  \end{aligned}
\end{equation}
The windings are computed as $\vartheta$ goes from $0$ to $1$ positively. Both vectors $\sigma(\vartheta)$ and $(\Psi^*\partial_x)\circ \varphi(r(\vartheta),\vartheta)$ belong to $\Psi^*(0\times \R^2) \subset TU$. Naturally in $\Psi^*(0\times \R^2)$ we have the complex structure $\Psi^*i$, where $i$ is the usual positive rotation by 90 degrees on $\R^2$. Moreover, the projection $\pi : \Psi^*(0\times \R^2) \to \xi$ is a vector bundle isomorphism since $\Psi^*(0\times \R^2)|_P = \xi|_P$ and $U$ is small. Clearly if $J$ is a $d\lambda$-compatible complex structure on $\xi$ then $\pi^*J$ is a complex structure on $\Psi^*(0\times \R^2)$ which is homotopic through complex structures to  $\Psi^*i$. This follows since $\Psi^*(d\theta\wedge dx\wedge dy)$ is a positive multiple of $\lambda \wedge d\lambda$ on $U$. Consequently we have
\begin{equation}\label{wind2}
  \wind(\pi \cdot \sigma(\vartheta), \pi \cdot (\Psi^*\partial_x)\circ \varphi(r(\vartheta),\vartheta)) = \wind (\Psi_*\sigma(\vartheta), \partial_x|_{\Psi\circ \varphi(r(\vartheta),\vartheta)})
\end{equation}
where the last winding is computed as a winding of non-vanishing vectors on $\R^2$ with respect to $i$.

We claim $\Psi_*\sigma(\vartheta)$ and $\frac{d}{d\vartheta} [ (\Psi\circ \varphi)(r(\vartheta),\vartheta) ]$ are everywhere linearly independent. To see this note that
\[
  \Psi_*\sigma(\vartheta) = d(\Psi \circ u(\vartheta))|_{z(r(\vartheta))} \cdot e^{i2\pi\beta(\vartheta)}
\]
and
\[
  \frac{d}{d\vartheta} [ (\Psi\circ \varphi)(r(\vartheta),\vartheta) ] = \partial_\vartheta (\Psi \circ u)|_{(\vartheta,z(r(\vartheta)))} \cdot 1 + d(\Psi\circ u(\vartheta))|_{z(r(\vartheta))} \cdot \frac{d}{d\vartheta}(z\circ r(\vartheta))
\]
so if these vectors are linearly dependent for some value of $\vartheta^*$ we get a contradiction to the fact that the derivative of the diffeomorphism $(\vartheta,z) \mapsto u(\vartheta,z)$ is non-singular at the point $(\vartheta^*,z(r(\vartheta^*)))$. We obtain
\begin{equation}\label{wind3}
  \wind (\Psi_*\sigma(\vartheta), \partial_x|_{\Psi\circ \varphi(r(\vartheta),\vartheta)}) = \wind \left( \frac{d}{d\vartheta} [ (\Psi\circ \varphi)(r(\vartheta),\vartheta) ], \partial_x|_{\Psi\circ \varphi(r(\vartheta),\vartheta)} \right).
\end{equation}

We claim that the last winding above equals $+1$. The image $\Gamma$ of the map $\vartheta \mapsto (\Psi\circ \varphi)(r(\vartheta),\vartheta)$ is an embedded circle in the plane $\{\theta = \theta_2\} \simeq \R^2$ and it winds around the origin. Otherwise $P_1$ would be homologically trivial in $S^3\setminus P$ because $\{ \varphi(r,\vartheta) : r(\vartheta)\leq r\leq 1\}$ is an embedded annulus inside $S^3\setminus P$ with one component equal to $\{\varphi(r(\vartheta),\vartheta)\}$ and the other equal to $P_1$. Moreover, $\Gamma$ bounds a disk $D \subset \{\theta = \theta_2\}$. Let us give $\Gamma$ the orientation induced from $\R/\Z$ by the embedding $\vartheta \mapsto (\Psi\circ \varphi)(r(\vartheta),\vartheta)$. By homotopy invariance the linking number $\link(P_1,P) = +1$ coincides with $\link(\Gamma,\theta\text{-axis})$, where the $\theta$-axis is oriented by $\partial_\theta$. Since the ambient space is $3$-dimensional we obtain $+1 = \link(\theta\text{-axis},\Gamma)$ which equals the intersection number of the $\theta$-axis with $D$, when $D$ is equipped with the orientation $o$ induced from the identity $\Gamma = \partial D$. Thus $o$ coincides with the orientation of $D$ induced by $dx\wedge dy$. This discussion shows that $\Gamma$ winds positively around the origin when $\R^2$ is oriented by $dx\wedge dy$. Now note that $\frac{d}{d\vartheta} [ (\Psi\circ \varphi)(r(\vartheta),\vartheta) ]$ is precisely a positive non-vanishing vector field tangent to $\Gamma$. So its winding coincides with the winding of $\Gamma$, proving our claim. Combining~\eqref{wind1}, \eqref{wind2} and \eqref{wind3} we obtain $m=1$, concluding the proof of Proposition~\ref{propn_sl}.

\appendix

\section{Proof of Lemma~\ref{lemma_anyJ}}\label{appA}

Let $M$, $\alpha$ and $P=(x,T)$ be as in the statement of the lemma. We write $\J = \J(\xi,d\alpha)$ and denote by $\J_{\text{fast}}(P)$ the set of $J \in \J$ for which there exists some embedded fast finite-energy $\jtil$-holomorphic plane asymptotic to $P$. Here $\jtil$ is defined in terms of $J$ and $\alpha$ by~\eqref{alm_cpx_str_def}. Since $\J$ is connected we may split the proof in the following two steps. \\

\noindent {\bf Step 1:} $\J_{\text{fast}}$ is closed. \\

\noindent {\bf Step 2:} $\J_{\text{fast}}$ is open.

\subsection{Proof of Step 1}

We follow~\cite{hry} closely. Let $H\subset M$ be a compact set with the following properties: $H \cap x(\R) = \emptyset$ and every continuous disk-map spanning $P$ passes through $H$. For every $J \in \J$ consider the set $\Lambda(J)$ consisting of fast embedded finite-energy planes $\util=(a,u)$ asymptotic to $P$ satisfying
\begin{equation}\label{norm_conditions}
  \begin{array}{cccc}
    \int_\D u^*d\alpha = T-\gamma, & u(0) \in H & \text{ and } & a(2) = 0.
  \end{array}
\end{equation}
Here $\gamma>0$ satisfies $\gamma <T'$ for every $P'=(x',T') \in \P$, and $\gamma < |T_0-T_1|$, for every $P_i=(x_i,T_i)$ such that $T_0\neq T_1$ and $T_0,T_1 \leq T$.

Consider $J_k \in \J_{\text{fast}}$ such that $J_k \to J$, and a sequence $\util_k=(a_k,u_k)$ of embedded fast $\jtil_k$-holomorphic planes asymptotic to $P$. After reparametrizing and translating in the $\R$-direction we may assume $\util_k \in \Lambda(J_k)$. Following~\cite{93}, let
\[
 \Gamma = \{ z\in \C \mid \exists z_j \to z \text{ and } k_j \text{ such that } |d\util_{k_j}(z_j)| \to \infty \}.
\]
It follows directly from Lemma~\ref{beforeclaim} and~\eqref{norm_conditions} that $\Gamma \subset \D$ and, up to the choice of a subsequence, we may assume $\#\Gamma<\infty$. Elliptic estimates show that, up to selection of a further subsequence, we may find a smooth $\jtil$-holomorphic map $\util =(a,u): \C\setminus \Gamma \to \R\times M$ such that $\util_k \to \util$ in $C^\infty_{loc}(\C\setminus \Gamma)$. It can be assumed that $\Gamma$ consists of negative or removable punctures and, consequently, $\infty$ must be positive. There is no loss of generality to assume that $\Gamma$ consists of negative punctures. Using results on cylinders with small $d\alpha$-area, as in~\cite{long} or~\cite{fols}, one proves that $\util$ is asymptotic to $P$ at $\infty$.

\begin{lemma}
$\int_{\C\setminus \Gamma} u^*d\alpha > 0$.
\end{lemma}

\begin{proof}
If $\Gamma=\emptyset$ then $\int_\D u^*d\alpha = \lim_k \int_\D u_k^*d\alpha = T-\gamma > 0$. If $\Gamma \neq\emptyset$ then clearly $\util$ is not constant. Arguing indirectly, assume $\Gamma\neq \emptyset$ and $\int_{\C\setminus \Gamma} u^*d\alpha = 0$. By Theorem~\ref{zera_dlambda_theorem}, $P$ cannot be simply covered if $\#\Gamma \geq 2$, proving $\#\Gamma=1$. In particular, it follows that there are complex constants $A\neq 0$ and $D$ such that $\util(z) = F_P(Az+D)$ where $F_P$ is the map $\est\in \C\setminus 0 \mapsto (Ts,x(Tt)) \in \R\times M$. In this case we must have $\Gamma = \{-D/A\} \subset \partial \D$ since, otherwise, we would get $$ T = \int_{\partial \D} u^*\alpha = \lim_k \int_{\partial \D} u_k^*\alpha = T-\gamma. $$ But this implies $u(0) \in H\cap x(\R)$, a contradiction.
\end{proof}

Let us enumerate $\Gamma = \{z_1,\dots,z_N\}$ and consider, for each $i$, the asymptotic orbit $P_i=(x_i,T_i)$ of $\util$ at $z_i$. Each $P_i$ is contractible. In fact, if $\rho$ is small and fixed, the contractible loop $t\mapsto u_k(z_j + \rho e^{i2\pi t})$ converges to $t\mapsto u(z_j + \rho e^{i2\pi t})$ which is homotopic to $x_i(T_it)$ in view of Theorem~\ref{partial_asymptotics}. Denoting ${x_i}_{T_i}(t) = x_i(T_it)$, let $\kappa_i$ be a non-vanishing section of ${x_i}_{T_i}^*\xi$ which extends to a non-vanishing section of $\xi|_{\mathcal D_i}$ for some (and hence any) disk $\mathcal D_i$ spanning $P_i$. Consider also $\kappa$ a non-vanishing section of $x_T^*\xi$ which extends over $\xi|_{\mathcal D}$ with no zeros, for spanning disks $\mathcal D$ for $P$. Let us compactify $\C\setminus \Gamma$ and obtain a surface with boundary $\Sigma$ by adding a circle at each point of $\Gamma$ and at $\infty$. By Theorem~\ref{partial_asymptotics} the map $u$ extends to a continuous map $\bar u:\Sigma\to M$ sending each boundary component onto the corresponding orbit $P_1,\dots,P_N,P$. We claim that the section $\kappa_1,\dots,\kappa_N,\kappa$ can be extended to a non-vanishing section $Z$ of $\bar u^*\xi$. In fact, there is no restriction to extending  $\kappa_1,\dots,\kappa_N$ over $\Sigma$ as a non-vanishing section $Z^*$. Disks $\mathcal D_i$ spanning $P_i$ can be glued to $\Sigma$ at $z_i$ to obtain a spanning disk $\mathcal D = u\#\mathcal D_1\#\dots \#\mathcal D_N$ for $P$. By our choice of the $\kappa_i$, $Z^*$ can be extended over $\mathcal D$ as a non-vanishing section. This shows that $\kappa$ does not wind with respect to $Z^*$, so $Z^*$ can be modified in a neighborhood of $\infty$ to a non-vanishing section $Z$ coinciding with $\kappa$ at $P$. This special non-vanishing section of $u^*\xi$ will be denoted by $B$.

Since $P$ is simply covered, there exists  Martinet tube $(U,\Psi)$ for $P$ such that $\kappa = (\Psi^*\partial_x) \circ x_T$. There is no loss of generality if we assume $B(\est) = (\Psi^*\partial_x) \circ u(\est)$ for large values of $s$. Again using results on cylinders with small $d\alpha$-area, one finds $R_0\gg1$ such that $|z|\geq R_0 \Rightarrow u_k(z) \subset U$, for every $k$. By Theorem~\ref{precise_asymptotics} we can also assume $|z|\geq R_0 \Rightarrow u(z) \subset U$ and $\pi \cdot du(z) \neq 0$. Here we used~\eqref{norm_conditions}.

\begin{lemma}
$\wind_\infty(\util,\infty,B) \leq 1$.
\end{lemma}

\begin{proof}
Since the maps $u_k$ provide capping disks for $P$, the section $(\Psi^*\partial_x)\circ u_k$ (defined only for $|z|\geq R_0$) can be extended to a non-vanishing section $B_k$ of $u_k^*\xi$. This follows from our particular choice of Martinet tube $(U,\Psi)$ made above, and implies $1=\wind_\infty(\util_k) = \wind_\infty(\util_k,\infty,B_k)$. We write $u_k(s,t) = u_k(\est)$ and $u(s,t) = u(\est)$. Let $s_k \to +\infty$ be so that
\[
\begin{array}{ccc}
  \pi\cdot \partial_s u_k(s_k,t) \neq0 \ \forall t & \text{and} & \wind(t\mapsto \pi \cdot \partial_su_k(s_k,t), B_k(e^{2\pi(s_k+it)})) = 1.
\end{array}
\]
For any fixed $s\geq (2\pi)^{-1} \log R_0$ the vectors $\pi \cdot \partial_s u(s,t)$ do not vanish, so that $\pi \cdot \partial_s u_k(s,t)$ do not vanish when $k\gg1$. We can estimate using standard degree theory
\[
 \begin{aligned}
   1 &- \wind(t\mapsto \pi \cdot\partial_s u(s,t), t\mapsto B(\est)) \\
   & = 1 - \wind(t\mapsto \pi \cdot\partial_s u(s,t), t\mapsto (\Psi^*\partial_x)\circ u(s,t)) \\
   & = \lim_{k\to \infty} 1 - \wind(t\mapsto \pi \cdot\partial_s u_k(s,t), t\mapsto (\Psi^*\partial_x)\circ u_k(s,t)) \\
   &= \lim_{k\to \infty} \left\{ \begin{array}{c}
                              \wind(t\mapsto \pi \cdot\partial_s u_k(s_k,t), t\mapsto B_k(s_k,t)) \\
                              - \wind(t\mapsto \pi \cdot\partial_s u_k(s,t), t\mapsto B_k(s,t))
                            \end{array} \right\} \\
   &= \lim_{k\to \infty} \# \left\{\text{zeros of } \pi \cdot du_k \text{ on } \{ e^{2\pi s} \leq |z| \leq e^{2\pi s_k} \} \right\}
 \end{aligned}
\]
Since the $\pi\cdot du_k$ satisfy a Cauchy-Riemann type equation, the (algebraic) count of zeros in the last line above is non-negative. The lemma follows by taking the limit as $s\to+\infty$.
\end{proof}

We will now argue indirectly to show that $\Gamma = \emptyset$. Fix $i$ and let $(s,t) \in \R^- \times \R/\Z$ be negative holomorphic cylindrical coordinates centered at $z_i$, as explained in Remark~\ref{cyl_coord}, and write $B(s,t)$ for the value of $B$ at the point corresponding to $(s,t)$. In particular, $B(s,t) \to \kappa_i(t)$ uniformly in $t$ as $s\to -\infty$. By Theorem~\ref{precise_asymptotics} there exists a smooth non-vanishing function $f(s,t)$ such that $f(s,t) \pi \cdot \partial_s u(s,t) \to e(t)$ in $C^0(\R/\Z,\xi)$ as $s\to-\infty$, where $e(t)$ is an eigenvector of $A_{P_i}$ corresponding to a positive eigenvalue. Consequently we can estimate
\[
  \begin{aligned}
    \wind_\infty(\util,z_i,B) &= \lim_{s\to -\infty} \wind (t\mapsto \pi\cdot \partial_su(s,t),t\mapsto B(s,t)) \\
    &= \wind(t\mapsto e(t),t\mapsto \kappa_i(t)) \geq 2
  \end{aligned}
\]
where in the last inequality we used $\mu(P_i)\geq 3$. Assuming $\Gamma\neq \emptyset$ we obtain the following contradiction
\[
  \begin{aligned}
    0 \leq \wind_\pi(\util) = \wind_\infty(\util) -2 + \#\Gamma + 1 \leq 1- 2\#\Gamma + \#\Gamma -1 = -\#\Gamma.
  \end{aligned}
\]
We used Lemma~\ref{lemma_wind_pi_infty}. This shows $\Gamma=\emptyset$ and proves that $\util$ is a non-constant fast finite-energy $\jtil$-holomorphic plane asymptotic to $P$.

To conclude Step 1 it remains only to prove that $\util$ is embedded. The identity $\wind_\pi(\util)=0$ shows that $\util$ is immersed. Self-intersection points of $\util$ must be isolated since, otherwise, one could use Carleman's similarity principle to show that $\util$ is a non-trivial cover of a somewhere injective finite-energy plane, see~\cite{props2}. In particular, this would imply that $P$ is not simply covered, a contradiction. Isolated self-intersections of $\util$ would imply self-intersections of $\util_k$ for $k\gg 1$, by stability and positivity of intersections of pseudo-holomorphic immersions. This shows $\util$ has no self-intersections and, consequently, $\util \in \Lambda(J)$ and $J \in \J_{\text{fast}}$.

\subsection{Proof of Step 2}

We need to revisit the functional analytic set-up for the Fredholm theory of embedded finite-energy surfaces constructed in~\cite{props3}, however, we will slightly modify the discussion there. From now on fix $J_0 \in \J_{\text{fast}}$ and an embedded fast finite-energy $\tilde J_0$-holomorphic plane $\util=(a,u):\C\to \R\times M$ asymptotic to $P$ at $\infty$.

With $l\geq 1$ fixed let $\mathcal K^l$ be the set of $C^l$-sections $K$ of $\mathcal L_\R(\xi)$ satisfying
\[
  \begin{array}{ccc}
    J_0K + KJ_0 = 0 & \text{and} & d\alpha(u,Kv) + d\alpha(Ku,v) = 0 \ \ \forall u,v\in \xi_p,p\in M.
  \end{array}
\]
Then $\mathcal K^l$ becomes a Banach space with the $C^l$-norm $\Vert\cdot\Vert_{C^l}$. When $\Delta>0$ is small enough then every $K \in \mathcal K^l$ satisfying $\Vert K \Vert_{C^l}<\Delta$ induces some $d\lambda$-compatible complex structure $J$ on $\xi$ of class $C^l$ by $J = J_0\exp(-J_0K)$. The set of $J$ arising in this way will be denoted by $\mathcal U_\Delta$, and the bijective correspondence between $\mathcal U_\Delta$ and a $\Delta$-ball in $\mathcal K^l$ gives $\mathcal U_\Delta$ the structure of a (trivial) Banach manifold. If we can prove that $J_0$ has an open neighborhood $\mathcal O$ in $\mathcal U_\Delta$ such that $\mathcal O\cap \J \subset \J_{\text{fast}}$ then it will follow that $J_0$ is an interior point of $\J_{\text{fast}}$ because $\mathcal O\cap \J$ is open in $\J$.

Let $U$ be a small neighborhood of $x(\R)$ in $M$ and $$ \Psi:U \stackrel{\sim}{\to} \R/\Z\times B $$ be a Martinet tube as described in~\ref{asymptotic_behavior_section}. Thus $\Psi^*(f(d\theta+xdy)) = \alpha$ and $\Psi(x(T\theta)) = (\theta,0,0)$ where $(\theta,z=(x,y)) \in B$ are standard coordinates in $\R/\Z\times \R^2$, and $f>0$ satisfies $f\equiv T$, $df\equiv0$ on $\R/\Z\times 0$. We have an induced frame
\begin{equation}\label{sections_tube}
  \{e_1 \simeq \partial_x, e_2 \simeq -x\partial_\theta + \partial_y\}
\end{equation}
of $\xi|_U$, and we assume $\{e_1\circ F,e_2\circ F\}$ extends to a frame of $F^*\xi$ for some (and hence any) continuous disk map $F:\D \to M$ satisfying $F(e^{i2\pi t}) = x(Tt)$.

Let $g$ be any $\R$-invariant Riemannian metric on $\R\times M$ such that $\jtil_0$ is an isometry, and denote by $N\util$ the normal bundle of $\util(\C)$. The asymptotic behavior of $\util$ explained in~\ref{asymptotic_behavior_section} (see~\cite{props1}) implies that $\tilde\xi_{u(z)} \cap d\util_z(T_z\C) = 0$ when $|z|$ is large, where $\tilde\xi=\pi_M^*\xi$. Here $\pi_M:\R\times M \to M$ denotes the projection onto the second factor. We may consider a $\jtil_0$-invariant subbundle $L\subset \util^*T(\R\times M)$ that coincides with $\util^*\tilde\xi$ on $\C\setminus B_{R_0}(0)$, where $R_0\gg1$, and coincides with $N\util$ on $B_{R_0-1}(0)$. Possibly after making $R_0$ larger, we have also $|z|\geq R_0 \Rightarrow u(z) \in  U$. One finds a $(d\alpha,J_0)$-unitary frame $\{n_1,n_2\}$ of $\xi|_U$ such that $\{\tilde n_1(\util(z)),\tilde n_2(\util(z))\}$ extends to a smooth $\jtil_0$-complex frame of $L$. Here $\tilde n_i(\util(z)) = n_i(u(z))$, $i=1,2$. The identity
\begin{equation}\label{winding_normal_xi}
  \wind(t\mapsto n_1 \circ u (Re^{i2\pi t}), t\mapsto e_1 \circ u (Re^{i2\pi t}))=1 \ \ (R\gg1)
\end{equation}
follows immediately from Theorem 1.8 in~\cite{props3}.

Let $\bar n_i$ be extensions of $\tilde n_i \circ \util$ to $\C$ such that $\{\bar n_1,\bar n_2\}$ is a $\jtil_0$-complex frame. Since $g$ is $\R$-invariant, its injectivity radius is uniformly bounded away from zero all over $\R\times M$. Thus we can find a small ball $B'\subset \C$ centered at $0$ such that the map
\[
\begin{array}{cc}
  \Phi : \C\times B' \to \R\times M, & (z,w) \mapsto \exp_{\util(z)} ( \Re(w) \bar n_1(z) + \Im(w) \bar n_2(z) )
\end{array}
\]
is an embedding onto a neighborhood of $\util(\C)$. Moreover, the image of any map $U:\C\to \R\times M$ sufficiently close to $\util$ in the strong $C^1$-topology coincides with $\{\Phi(z,v(z)) \mid z\in \C\}$, for some $v:\C\to B'$ small in the strong $C^1$-topology.

In the following we denote $\bar J = \Phi^*\jtil$ and write $$ \bar J(z,v) = \begin{bmatrix} j_1(z,v) & \Delta_1(z,v) \\ \Delta_2(z,v) & j_2(z,v) \end{bmatrix} $$ in $2\times 2$-blocks, for every $J \in \mathcal U_\Delta$. Let $\sigma(s,t) = \est$. After making $B'$ smaller and $R_0$ larger the map $$ F(s,t,v) = (id_\R\times \Psi) \circ \Phi \circ (\sigma\times id_{\R^2})(s,t,v) $$ is defined for $e^{2\pi s}\geq R_0$ and $v\in B'$. Using the asymptotic behavior of $\util$ from Theorem~\ref{precise_asymptotics} and the $\R$-invariance of $g$ one shows
\begin{equation}\label{map_F}
  \limsup_{s\to+\infty} \sup_{v\in K,t\in\R/\Z} |D^\beta F(s,t,v) | < \infty
\end{equation}
for every non-trivial multi-index $\beta=(\beta_1,\beta_2,\beta_3,\beta_4)$ and $K\subset B'$ compact. For each $J \in \mathcal U_\Delta$ consider the almost complex structure $\hat J$ on $\R\times \R/\Z\times B$ defined by $(id_\R\times \Psi)^*\hat J = \jtil$. Seen as a smooth $4\times 4$ matrix, $\hat J$ is independent of the $\R$-coordinate. If $\underline J = (\sigma\times id_{\R^2})^*\bar J = F^*\hat J$ then, in view of~\eqref{map_F}, we have
\begin{equation}\label{map_J_underline}
  \limsup_{s\to+\infty} \sup_{v\in K,t\in\R/\Z} |D^\beta \underline J(s,t,v) | < \infty
\end{equation}
for every $\beta$ and $K\subset B'$ compact.

Given $\delta<0$, $\gamma \in (0,1)$ and a finite-dimensional real vector space $V$ we can define the Banach space $C^{l,\gamma,\delta}_0(\C,V)$ of $C^{l,\gamma}$ functions $v: \C\to V$ such that the map $(s,t) \mapsto v\circ \sigma(s,t)$, $s\geq 0$, belongs to $C^{l,\gamma,\delta}_0([0,+\infty)\times \R/\Z,V)$, as explained in Remark~\ref{topology}. The norm is defined by $\Vert v\Vert_{C^{l,\gamma}(\D)} + \Vert v\circ\sigma \Vert_{l,\gamma,\delta}$, where $\Vert \cdot\Vert_{l,\gamma,\delta}$ is the norm~\eqref{norma_alpha_delta}. From now on we fix $l\geq 2$, $\gamma\in(0,1)$ as above, and $\delta<0$ will be fixed \emph{a posteriori}. Clearly the subset $\V \subset C^{l,\gamma,\delta}_0(\C,\C)$ of maps with image in $B'$ is an open set.

Let $Y\to \C$ be the vector bundle with fibers $Y_z = \{ \R\text{-linear maps } T_z\C \to \C \}$. The space $C^{l-1,\gamma,\delta}_0(Y)$ consists of sections $A:\C \to Y$ of class $C^{l-1,\gamma}$ such that $(s,t) \mapsto A(\sigma(s,t))\cdot \partial_t\sigma(s,t)$ belongs to $C^{l-1,\gamma,\delta}_0$ on $\R^+ \times \R/\Z$. The norm is $\|A\|_{C^{l-1,\gamma}(\D)} + \|A\circ\sigma \cdot\sigma_t\|_{l-1,\gamma,\delta}$.

Similarly to~\cite{props3} consider the smooth map $H:\V \times \mathcal U_\Delta \to C^{l-1,\gamma,\delta}(Y)$ given by
\begin{equation}\label{badH}
 H(v,J) = \Delta_2(z,v) + j_2(z,v)\cdot dv - dv \cdot j_1(z,v) - dv\cdot \Delta_1(z,v)\cdot dv.
\end{equation}
One has to make use of~\eqref{map_J_underline} to verify that $H$ takes values on $C^{l-1,\gamma,\delta}(Y)$. Standard arguments show that $H$ is smooth. The equation $H(v,J)=0$ holds if, and only if, the embedding $z\mapsto (z,v(z))$ has a $\bar J$-invariant tangent space. Writing $$ \bar J_0 = \begin{bmatrix} j^0_1 & \Delta^0_1 \\ \Delta^0_2 & j^0_2 \end{bmatrix} $$ in $2\times 2$-blocks, we have $j^0_1(z,0)=j^0_2(z,0)=i$ and $\Delta^0_1(z,0)=\Delta^0_2(z,0)=0$. Differentiating~\eqref{badH} at $v=0$ we get
\begin{equation}\label{badH_linear}
  D_1H(0,J_0)\zeta = i \cdot d\zeta - d\zeta \cdot i + C \cdot\zeta
\end{equation}
where $C(z) = D_2\Delta^0_2(z,0)$. Let $d\phi_{Tt}:\xi|_{x(0)} \to \xi_{x(Tt)}$ be represented as a path $\varphi(t) \in \sp(1)$ using the frame $\{n_1,n_2\}$, and let $S = -i\varphi'\varphi^{-1}$. Denoting by $\tilde C(s,t)$ the linear map $u\mapsto (C(\sigma(s,t)) \cdot u) \cdot \partial_t\sigma(s,t)$, the following asymptotic behavior is proved in~\cite{props3}
\begin{equation}\label{}
  \lim_{s\to\infty} \sup_t |D^\beta[\tilde C(s,t)-S(t)]| = 0, \ \forall |\beta|\leq l.
\end{equation}
There is a corresponding representation $L = -i\partial_t - S$ of the asymptotic operator at $P$. From now one we assume $\delta \in (-\infty,0) \setminus \sigma(L)$, and denote by $\mu_{CZ}^N(P,\delta)$ the $\delta$-weighted Conley-Zehnder of $P$ computed with respect to the frame $\{n_1,n_2\}$, which is defined as twice the winding of the largest eigenvalue of $L$ below $\delta$, plus $0$ or $1$ depending whether the winding jumps when we compare with the smallest eigenvalue above $\delta$.

There exists a bundle splitting $Y = Y^{1,0} \oplus Y^{0,1}$ and we may define a Banach space of sections $C^{l-1,\gamma,\delta}_0(Y^{0,1})$ analogously as above. Differentiating the identity $\bar J_0^2=-I$ one shows $C(z)\cdot \zeta$ is $\C$-anti-linear, so that the linear map $D_1H(0,J_0)$ may be viewed as an operator $D_0:C^{l,\gamma,\delta}_0(\C,\C) \to C^{l-1,\gamma,\delta}_0(Y^{0,1})$.

\begin{theorem}[Hofer, Wysocki and Zehnder]\label{fredholm_operator_thm}
The operator $D_0$ is Fredholm with index $\mu_{CZ}^N(P,\delta) + 1$. Moreover, there exists a smooth Banach bundle $\E\to \V \times\mathcal U_\Delta$, with fibers modeled on $C^{l-1,\gamma,\delta}_0(Y^{0,1})$, and a smooth section $\eta$ such that $\eta(v,J)=0 \Leftrightarrow H(v,J)=0$ and the partial vertical derivative $D_1\eta(0,J_0)$ coincides with $D_0$.
\end{theorem}

Since $\mu_{CZ}(P) \geq 3$ we have $\mu(\varphi) \geq 1$ by~\eqref{winding_normal_xi} and by our choice of Martinet tube. This allows room for us to place $\delta<0$ precisely on the spectral gap between the largest eigenvalue with winding $0$ and the smallest eigenvalue with winding $1$, so that $\mu_{CZ}^N(P,\delta) = 1$ and the Fredholm index of $D_1\eta(0,J_0)$ is $2$.

\begin{lemma}\label{D_1H_2_surj}
With this choice of $\delta$ the operator $D_0$ is surjective.
\end{lemma}

\begin{proof}
Let $a \in \ker D_0$ be non-zero and denote $a(s,t) = a(\sigma(s,t))$. Then $a_s + ia_t + \tilde C a=0$. By standard asymptotic analysis from~\cite{props1} (see section 6 in~\cite{hry} for a detailed account), one finds an eigenvector $e(t)$ of $L$ associated to some eigenvalue $\nu <0$ such that $$ a(s,t) = e^{\int_{s_0}^s h(\tau) d\tau} ( e(t) + R(s,t)) $$ where $h(\tau)-\nu$ and $R(s,t)$ converge to $0$ uniformly in $t$ as $s\to+\infty$. Since $a\in C^{l,\gamma,\delta}_0$ we must have $\nu < \delta$, so that $\wind(t\mapsto e(t)) \leq 0$. By standard degree theory the algebraic count of zeros of $a(z)$ on $\C$ is $\leq 0$. However, Carleman's Similarity principle implies that all zeros count positively. Hence, $a$ never vanishes. If the kernel had 3 or more linearly independent sections then a non-trivial linear combination of them would have to vanish at some point and, consequently, would vanish identically, providing a contradiction. This shows that $\dim \ker D_0\leq 2$ and the conclusion follows since the index is 2.
\end{proof}

Clearly $\ker D_0$ splits. By the implicit function theorem we find, possibly after shrinking $\V$, a neighborhood $\O$ of $J_0$ in $\mathcal U_\Delta$ so that the (local) universal moduli space defined by
\[
  \M = \{ (v,J) \in \V\times \O \mid H(v,J) = 0 \}
\]
is a smooth Banach submanifold with tangent space $\ker D_0$, and the projection $\Pi(v,J) = J$ restricts to a submersion of $\M$ onto $\O$.

Hence, for any fixed $J \in \O\cap \J$ there exists some $v\in \V$ so that the embedding $z\mapsto (z,v(z))$ has a $\bar J$-invariant tangent space. Following the appendix of~\cite{props3}, for a given $0<\epsilon<2\pi$ it is possible to find $\phi \in C^{l,\gamma,-\epsilon}_0(\C,\C)$ so that the map $\psi(z) = z+\phi(z)$ is a diffeomorphism and $\wtil(z) = \Phi(\psi(z),v\circ \psi(z))$ is an embedded finite-energy $\jtil$-holomorphic plane asymptotic to $P$. To conclude Step 2 it remains to show

\begin{lemma}
The plane $\wtil$ above is fast.
\end{lemma}

\begin{proof}
Consider $(S(s,t),T(s,t)) = \sigma^{-1} \circ\psi \circ \sigma(s,t) \in \R\times \R/\Z$. It follows from the formula $S+iT = (2\pi)^{-1}\log(\psi\circ\sigma (s,t))$ that $|(S-s,T-t)| + |D(S-s,T-t)| \to 0$ as $s\to \infty$, uniformly in $t\in\R/\Z$. Let us write
\[
  (a(s,t),\theta(s,t),x(s,t),y(s,t)) = (id_\R\times \Psi) \circ \util \circ \sigma(s,t) = F(s,t,0),
\]
and
\[
  \begin{aligned}
    (A(s,t),&\Theta(s,t),X(s,t),Y(s,t)) = (id_\R\times \Psi) \circ \wtil\circ \sigma(s,t) \\
    &= F(S,T,v\circ \sigma(S,T)).  \\
  \end{aligned}
\]
Using~\eqref{map_F}, and the fact that $|v\circ\sigma(S,T)|$ decay as fast as $e^{\delta s}$ when $s\to\infty$, we can estimate
\[
  \begin{aligned}
    & |(X(s,t),Y(s,t))-(x(S,T),y(S,T))| \\
    & \leq |F(S,T,v\circ \sigma(S,T)) - F(S,T,0)| \leq Ce^{\delta s}
  \end{aligned}
\]
for some $C>0$, if $s\gg1$. In view of Theorem~\ref{precise_asymptotics} we know that all partial derivatives of $(x(s,t),y(s,t))$ of all orders decay exponentially to zero faster than $e^{(\mu+r)s}$ for every $r>0$, where $\mu<\delta$ is the asymptotic eigenvalue\footnote{To see why $\mu<\delta$ note that the winding of the asymptotic eigenvalue of $\util$ with respect to $e_1 \simeq \partial_x$ is $\wind_\infty(\util)=1$, as explained in Remark~\ref{asymp_evalue_wind_infty}. This is so since the Martinet tube is aligned with a trivialization of $x_T^*\xi$ induced by some (and hence any) capping disk. But, by our choice of $\delta$, if an eigenvalue winds $\leq 0$ with respect to $n_1$ (or $\leq 1$ with respect to $e_1$) it must be smaller than $\delta$.}. Taking $r>0$ small we get decay faster than $e^{\delta s}$. Since the partial derivatives of $\phi(\est)$ with respect to $s$ and $t$ decay to $0$ uniformly (like $e^{-\epsilon s}$), we conclude that
\[
  |(x(S,T),y(S,T)) - (x(s,t),y(s,t))| \leq C'e^{\delta s} \text{ when }s\gg 1.
\]
Both estimates above give
\[
  |(X(s,t),Y(s,t))-(x(s,t),y(s,t))| \leq C'' e^{\delta s}
\]
which implies
\[
  \begin{aligned}
    |(X(s,t),Y(s,t))| &\leq |(x(s,t),y(s,t))| + |(X(s,t),Y(s,t))-(x(s,t),y(s,t))| \\
    &\leq C''' e^{\delta s}.
  \end{aligned}
\]
Now, in view of Theorem~\ref{precise_asymptotics}, $(X(s,t),Y(s,t))$ decays slower than $e^{(\nu-r)s}$ for some $\nu \in \sigma(L)\cap (-\infty,0)$ and any $r>0$. Here $\nu$ is the asymptotic eigenvalue of $\wtil$. Hence the estimate above shows that $\nu < \delta$, implying $\wind_\infty(\wtil)\leq 1$. Thus $\wind_\infty(\wtil)=1$ by Lemma~\ref{lemma_wind_pi_infty}.
\end{proof}

We showed that $\mathcal O\cap \J\subset \J_{\text{fast}}$.

\begin{remark}
It should be noted that the elegant work of Wendl~\cite{wendl1,wendl2} can be used to start a new proof of the above mentioned compactness and automatic transversality result for embedded fast finite-energy planes.
\end{remark}

\section{Proof of Proposition~\ref{propcomp}}

First we need two auxiliary lemmas.

\begin{lemma}\label{appB_lemma_1}
Suppose $K_k$, $K^\infty_k$, $K^\infty$ and $L$ are as in the statement of Proposition~\ref{propcomp}. For each $s\geq0$ and $k\in\Z^+$ consider the unbounded self-adjoint operator
\[
 \begin{array}{cc}
   L_k(s):W^{1,2}\subset L^2 \rightarrow L^2, & L_k(s)e = -J_0\dot{e} - S_k(s,\cdot)e
 \end{array}
\]
where $S_k(s,t) = \frac{1}{2} \left[ K_k(s,\cdot) + K_k(s,\cdot)^T \right]$. If $\delta\in(-\infty,0)\setminus\sigma(L)$ then we find $s_0\geq0$ and $c>0$ such that
\[
\norma{[L_k(s)-\delta]e}_{L^2} \geq c\norma{e}_{L^2} \ \forall s\geq s_0, \ \forall e\in W^{1,2}, \ \forall k\geq0.
\]
\end{lemma}

\begin{proof}
We proceed indirectly. Abbreviating $\norma{\cdot}_2=\norma{\cdot}_{L^2}$, suppose there are sequences $s_j \rightarrow+\infty$ and $\{e_j\}\subset W^{1,2}$ such that $\norma{e_j}_2\equiv1$ and $\norma{[L_{k_j}(s_j)-\delta]e_j}_2\rightarrow0$. The equation
\[
 \partial_t e_j = J_0[L_{k_j}(s_j)-\delta]e_j + J_0S_{k_j}(s_j,\cdot)e_j + \delta J_0e_j
\]
shows that $\{e_j\}$ is a bounded sequence in $W^{1,2}$. Consequently we can assume the existence of some $e\in L^2$ satisfying $\norma{e_j-e}_2\rightarrow0$. The above equation now shows that $\{e_j\}$ is a Cauchy sequence in $W^{1,2}$, which proves $e\in W^{1,2}$ and $e_j\rightarrow e$ in $W^{1,2}$. However, the asymptotic conditions imposed on the functions $K_{k_j}$ imply that $L_{k_j}(s_j)$ converges to $L$ in the strong operator topology of continuous linear maps from $W^{1,2}$ into $L^2$. Hence $[L-\delta]e=0$, contradicting $\delta\not\in\sigma(L)$ since $\norma{e}_2=1$.
\end{proof}

\begin{lemma}\label{asympdecaylemma}
Let $K_k$, $K^\infty_k$, $K^\infty$ and $L$ be as in the statement of Proposition~\ref{propcomp}. Suppose $\lambda<\delta<0$ are such that $\lambda\in\sigma(L)$ and $(\lambda,\delta]\cap\sigma(L)=\emptyset$. Then one can find $0<r<\delta-\lambda$ and $s_1\geq0$ such that
\[
\norma{X(s,\cdot)}_{L^2} \leq e^{(\delta-r)(s-s_1)}\norma{X(s_1,\cdot)}_{L^2} \ \ \forall s\geq s_1
\]
for any $k$ and for any smooth function $X :[0,+\infty)\times \R/\Z \to \R^{2n}$ satisfying
\begin{equation}\label{solution_X}
  \begin{array}{cc}
    \partial_s X + J_0\partial_t X + K_kX = 0, & \lim_{s\rightarrow+\infty} e^{-\delta s}\norma{X(s,\cdot)}_{L^2} = 0.
  \end{array}
\end{equation}
\end{lemma}

\begin{proof}
We abbreviate $\norma{\cdot}_2=\norma{\cdot}_{L^2}$ and assume $X \in C^\infty$ satisfies~\eqref{solution_X} for some $k$. Setting $Y=e^{-\delta s}X$, $S_k = \frac{1}{2}[K_k(s,t)+K_k(s,t)^T]$, $A_k(s,t) = \frac{1}{2}[K_k(s,t)-K_k(s,t)^T]$ we have
\begin{equation}\label{nice_decay_asymp_uniform}
  \limsup_{s,k\rightarrow+\infty} \sup_t \left( |D^\gamma[S_k(s,t)-K^\infty(t)]| + |D^\gamma A_k(s,t)| \right) = 0 \ \forall |\gamma| \leq l
\end{equation}
and
\begin{equation}\label{eqlocal2}
  Y_s + J_0 Y_t + K_kY + \delta Y = Y_s - [L_k(s)-\delta]Y + A_kY = 0
\end{equation}
where $L_k(s)$ is the operator defined in the statement of Lemma~\ref{appB_lemma_1}. Moreover, $\lim_{s\rightarrow+\infty} \norma{Y(s,\cdot)}_2 = 0$. If we set $g(s)=\frac{1}{2}\norma{Y(s,\cdot)}_2^2$ then
\[
  g^\prime(s) = \left< Y_s,Y \right> = \left< [L_k(s)-\delta]Y-A_kY,Y \right> = \left< [L_k(s)-\delta]Y,Y \right>
\]
and
\[
 \begin{aligned}
  g^{\prime\prime}(s) &= \left< -J_0Y_{ts} - S_kY_s - \delta Y_s - (\partial_sS_k) Y,Y \right> + \left< [L_k(s)-\delta]Y,Y_s \right> \\
  &= \left< [L_k(s)-\delta]Y_s,Y \right> - \left< (\partial_sS_k) Y,Y \right> + \left< [L_k(s)-\delta]Y,[L_k(s)-\delta]Y - A_kY \right> \\
  &= \left< [L_k(s)-\delta]Y - A_kY,[L_k(s)-\delta]Y \right> - \left< (\partial_sS_k) Y,Y \right> \\
  &+ \norma{[L_k(s)-\delta]Y}^2_2 - \left< [L_k(s)-\delta]Y,A_kY \right> \\
  &= 2\norma{[L_k(s)-\delta]Y}^2_2 - 2\left< [L_k(s)-\delta]Y,A_kY \right> - \left< (\partial_sS_k) Y,Y \right>.
 \end{aligned}
\]
We used that $L_k(s)$ is self-adjoint $\forall s$. Let $s_0$ and $c>0$ be given by the previous lemma. Then $s\geq s_0$ implies
\[
 \begin{aligned}
  g^{\prime\prime}(s) &\geq 2c^2\norma{Y}_2^2 - 2c\norma{Y}_2\norma{A_kY}_2 - \norma{Y}_2\norma{[\partial_sS_k]Y}_2 \\
  &\geq 4 g(s) \left( c^2 - c\norma{A_k}_{L^\infty(S^1)} - \frac{1}{2}\norma{\partial_sS_k}_{L^\infty(S^1)} \right)
 \end{aligned}
\]
In view of~\eqref{nice_decay_asymp_uniform}, for any $\epsilon>0$ small $\exists s_1\geq s_0$ be such that
\[
 s\geq s_1 \Rightarrow \left( c^2 - c\norma{A_k}_{L^\infty(S^1)} - \frac{1}{2}\norma{\partial_sS_k}_{L^\infty(S^1)} \right) \geq (c-\epsilon)^2 \ \forall k.
\]
Consequently we have an estimate $g^{\prime\prime}(s) \geq 4(c-\epsilon)^2 g(s)$ whenever $s\geq s_1$, independently of $k$ or $X$. Now we note that any positive $C^2$ function on $[s_1,+\infty)$ satisfying $g^{\prime\prime}(s) \geq 4(c-\epsilon)^2 g(s)$ and $\lim_{s\rightarrow +\infty} g(s) = 0$ must also satisfy $g(s)\leq g(s_1) e^{ -2(c-\epsilon)(s-s_1) }$ for all $s\geq s_1$. Again this is independent of $k$ or $X$. The proof is complete if we set $r<\min\{\delta-\lambda,c-\epsilon\}$ since $2g(s) = e^{-2\delta s}\norma{X(s,\cdot)}_2^2$.
\end{proof}

We now turn to the proof of Proposition~\ref{propcomp}. The definition of $E$ and the assumptions on $K_k$ together imply that
\begin{equation}\label{uniformdecayXn}
 \lim_{s\rightarrow+\infty} \sup_t e^{-\delta s}\left|D^\gamma X_k(s,t)\right| = 0 \ \forall \gamma
\end{equation}
holds for each fixed $k$. It is important here that $D^\gamma$ is a partial derivative of any order, and we used equations~\eqref{eqns_X_k} and the H\"older estimates for the operator $\partial_s+J_0\partial_t$. We fix a number $\lambda<\delta$ satisfying $(\lambda,\delta] \cap \sigma(L) = \emptyset$ and proceed in three steps.
\\

\noindent {\bf STEP 1:} $\forall m \geq 0$ one can find numbers $0<r_m<\delta-\lambda$, $s_m > 0$ such that
\[
 \sqrt{ \sum_{j=0}^m \norma{(\partial_s)^j X_k (s,\cdot)}_{L^2(S^1)}^2 } \leq e^{(\delta-r_m)(s-s_m)} \sqrt{ \sum_{j=0}^m \norma{(\partial_s)^j X_k (s_m,\cdot)}_{L^2(S^1)} } \ \forall s\geq s_m \ \forall k.
\]

\begin{proof}[Proof of STEP 1]
Fix $k\geq0$. We have equations
\[
 \left\{
  \begin{aligned}
   & \partial_sX_k + J_0\partial_tX_k + K_k X_k = 0 \\
   & \partial_{ss}^2X_k + J_0\partial_{ts}^2X_k + \partial_s[K_k X_k] = 0 \\
   & \cdots \\
   & (\partial_s)^{m+1}X_k + J_0\partial_t (\partial_s)^m X_k + \partial_s^m[K_k X_k] = 0
  \end{aligned}
 \right.
\]
Defining $Z_k(s,t) = \left( X_k , \partial_sX_k , \dots , \partial_s^mX_k \right) ^ T$ then $\partial_sZ_k + \hat J_0\partial_tZ_k + G_k Z_k = 0$ where
\[
 G_k(s,t) = \begin{bmatrix}
             K_k & 0 & \cdots & 0 \\
             H_k^{21} & K_k & \cdots & 0 \\
             \cdots \\
             H_k^{(m+1)1} & H_k^{(m+1)2} & \cdots & K_k  \\
            \end{bmatrix};\
 \hat J_0 = \begin{bmatrix}
             J_0 & 0 & \cdots & 0 \\
             0 & J_0 & \cdots & 0 \\
             \cdots \\
             0 & 0 & \cdots & J_0  \\
            \end{bmatrix}.
\]
$G_k$ is a lower triangular matrix of $2n\times2n$ blocks and every term $H_k^{ij}$ ($i>j$) below the diagonal satisfies
\[
 \lim_{s\rightarrow+\infty} \sup_{k,t} |D^\gamma H_k^{ij}(s,t)| = 0 \ \forall \gamma.
\]
This follows from the hypotheses of Proposition~\ref{propcomp}. These remarks and (\ref{uniformdecayXn}) show that we can apply Lemma~\ref{asympdecaylemma} to the $Z_k(s,t)$ and find numbers $s_m>0$ and $0<r_m<\delta-\lambda$ as desired. Note that $s_m$ and $r_m$ are independent of $k$.
\end{proof}

\noindent {\bf STEP 2:} $\forall l\geq0$ one can find $0<r_l<\delta-\lambda$, $s_l>0$ and $c_l>0$ such that
\begin{equation}
 \max_{\norma{\gamma}\leq l} \norma{D^\gamma X_k(s,\cdot)}_{L^\infty(S^1)} \leq c_l e^{(\delta-r_l)(s-s_l)} \sum_{|\beta|\leq l+1} \norma{D^\beta X_k(s_l,\cdot)}_{L^\infty(S^1)}
\end{equation}
for every $k\geq0$ and $s\geq s_l$.

\begin{proof}[Proof of STEP 2]
We prove by induction that $\forall m \geq 0 \ \exists C_m>0$, $0<r_m<\delta-\lambda$ and $s_m>0$ such that $s\geq s_m$ implies
\begin{equation}\label{indstep2again}
 \sqrt{ \sum_{|\gamma|\leq m} \norma{D^\gamma X_k (s,\cdot)}_{L^2(S^1)}^2 } \leq C_m e^{(\delta-r_m)(s-s_m)} \sqrt{ \sum_{|\gamma|\leq m} \norma{D^\gamma X_k (s_m,\cdot)}_{L^2(S^1)}^2 } \ \ \forall k.
\end{equation}
The assertion for $m=0$ is already proved in Step 1. Assuming (\ref{indstep2again}) is proved for $m$, we claim that it holds for $m+1$. Writing $D^\gamma = \partial_s^{i}\partial_t^{j}$ we proceed by induction on $j$ to show that
\begin{equation}\label{indstep2again_2}
 \sqrt{ \norma{\partial_s^{m+1-j}\partial_t^{j} X_k (s,\cdot)}_{L^2(S^1)}^2 } \leq C' e^{(\delta-r')(s-s')} \sqrt{ \widetilde \sum \norma{D^\gamma X_k (s',\cdot)}_{L^2(S^1)}^2 } \ \ \forall s\geq s'
\end{equation}
for every $0\leq j\leq m+1$, where $\widetilde\sum$ indicates a sum over all multi-indices $\gamma=(\gamma_1,\gamma_2)$ satisfying either $|\gamma|\leq m$ or $|\gamma|=m+1$ and $\gamma_2\leq j$, and the constants $C',r',s'$ are independent of $k$.

The case $j = 0$ follows from Step 1 for $m+1$ and from~\eqref{indstep2again} for $m$. Now fix $0< b \leq m+1$ and assume~\eqref{indstep2again_2} holds for $j=0,\dots,b-1$. Let $\beta = (m+1-b,b)$. Equation $\partial_s X_k + J_0\partial_t X_k + K_k X_k =0$ implies
\[
 \begin{aligned}
 \partial_s \left( \partial_s^{m+1-b}\partial_t^{b-1} X_k \right) + J_0 D^\beta X_k &= \partial_s^{m+1-b}\partial_t^{b-1} \left( \partial_s X_k + J_0 \partial_t X_k \right) \\
 &= - \partial_s^{m+1-b}\partial_t^{b-1} \left(K_kX_k\right).
 \end{aligned}
\]
Thus
\[
 D^\beta X_k = J_0 \left( \partial_s^{m+2-b}\partial_t^{b-1} X_k + \partial_s^{m+1-b}\partial_t^{b-1} \left(K_kX_k\right) \right)
\]
The asymptotic uniform bounds on derivatives of $K_k(s,t)$ in $s$ and $k$, and the induction hypothesis, imply~\eqref{indstep2again_2} holds for $j=0,\dots,b$. The induction step is complete and~\eqref{indstep2again_2} is proved for every $j\leq m+1$. This proves~\eqref{indstep2again} for $m+1$. We showed~\eqref{indstep2again} holds $\forall m$.

Using \eqref{indstep2again} for $m=l+1$ we obtain $\hat C_l>0$ independent of $n$ such that
\[
 \sqrt{ \sum_{|\gamma|\leq l} \norma{D^\gamma X_k (s,\cdot)}_{W^{1,2}(S^1)}^2 } \leq \hat C_l e^{(\delta-r_{l+1})(s-s_{l+1})} \sum_{|\gamma|\leq l+1} \norma{D^\gamma X_k (s_{l+1},\cdot)}_{L^\infty(S^1)}
\]
for every $s\geq s_{l+1}$ and $k\geq0$. The conclusion follows since $W^{1,2}(S^1)\hookrightarrow L^\infty(S^1)$.
\end{proof}

\noindent {\bf STEP 3:} There exists $X_\infty$ and a subsequence $X_{k_j}$ such that $X_{k_j} \rightarrow X_\infty$ in $C^{l,\alpha,\delta}_0$.

\begin{proof}[Proof of STEP 3]
Since $X_k$ is $C^\infty_{loc}$-bounded we can assume, up to selection of a subsequence, that $\exists X_\infty$ such that $X_k \rightarrow X_\infty$ in $C^\infty_{loc}$. Fix $\gamma$ and $\epsilon>0$. By the previous step $\exists s_1 \gg 0$ such that
\[
 s\geq s_1 \Rightarrow \sup_{k,t} e^{-\delta s} \left|D^\gamma X_k(s,t)\right| \leq \frac{\epsilon}{2}.
\]
This implies
\[
  \sup_{s\geq s_1,t\in S^1} e^{-\delta s} |D^\gamma X_\infty(s,t)| \leq \frac{\epsilon}{2}.
\]
We find $k_1 \in \Z^+$ such that
\[
 k \geq k_1 \Rightarrow \sup_{[0,s_1]\times \R/\Z} e^{-\delta s} \left| D^\gamma[X_k-X_\infty] \right| \leq \epsilon.
\]
Hence if $s\geq 0$ and $k\geq k_1$ then $e^{-\delta s} \left| D^\gamma[X_k-X_\infty] \right| \leq \epsilon$. We proved
\[
 \lim_{k\rightarrow \infty} \sup_{s,t} e^{-\delta s} \left| D^\gamma[X_k-X_\infty] \right| = 0
\]
for each fixed $\gamma$. If we fix $l\geq1$ then the above limits for $|\gamma|\leq l+1$ imply
\[
 \lim_{k\rightarrow+\infty} \max_{|\beta|\leq l} \norma{ e^{-\delta s} D^\beta[X_k-X_\infty] } _{C^{0,\alpha}([0,+\infty)\times \R/ \Z)} = 0.
\]
The conclusion follows.
\end{proof}

\end{document}